\documentclass[12pt]{amsart}
\usepackage{latexsym, amsmath, amssymb, amsthm}
\usepackage{epsfig}
\usepackage{amscd}
\usepackage{latexsym}
\usepackage{pstricks,pst-node,cite}

\newtheorem{thm}{Theorem}[section]
\newtheorem{lem}[thm]{Lemma}
\newtheorem{prop}[thm]{Proposition}
\newtheorem{cor}[thm]{Corollary}
\newtheorem{exa}[thm]{Example}

\theoremstyle{remark}

\newcommand{\tensor}{\otimes}

\newcommand{\Inv}{\mbox{\rm Inv}}

\newcounter{fignum}
\ifx\blackandwhite\undefined
  \psset{linecolor=blue} \else
  \psset{linecolor=blue}
\fi

\psset{linewidth=.6pt,dash=3pt 3pt,doublesep=.05,dotsize=1pt 5}
\SpecialCoor

\newcommand{\vv}{\begin{pspicture}[shift=-.4](-.6,-.5)(.6,.5)
\psbezier(.5;45)(.25;45)(.25;315)(.5;315)
\psbezier(.5;225)(.25;225)(.25;135)(.5;135)
\end{pspicture}}
\newcommand{\hh}{\begin{pspicture}[shift=-.4](-.6,-.5)(.6,.5)
\psbezier(.5;45)(.25;45)(.25;135)(.5;135)
\psbezier(.5;225)(.25;225)(.25;315)(.5;315)
\end{pspicture}}
\newcommand{\btwohvert}{\begin{pspicture}[shift=-.6](-.6,-.7)(.6,.7)
\qline(0, .25)(.35, .6)\qline(0, .25)(-.35, .6)
\qline(0,-.25)(.35,-.6)\qline(0,-.25)(-.35,-.6)
\psline[doubleline=true](0,-.25)(0,.25)
\end{pspicture}}
\newcommand{\btwohhoriz}{\begin{pspicture}[shift=-.4](-.8,-.5)(.8,.5)
\qline( .25,0)( .6,.35)\qline( .25,0)( .6,-.35)
\qline(-.25,0)(-.6,.35)\qline(-.25,0)(-.6,-.35)
\psline[doubleline=true](-.25,0)(.25,0)
\end{pspicture}}

\newcommand{\middlearrow}{\lput{:U}{\begin{pspicture}[shift=0](0,0)(0,0)
\psline[arrows=->,arrowscale=1.5](2.2pt,0)(2.3pt,0)\end{pspicture}}}
\newcommand{\mda}{\lput{:U}{\begin{pspicture}[shift=0](0,0)(0,0)
\psline[linecolor=black,arrows=->,arrowscale=1.7](3.2pt,0)(3.4pt,0)\end{pspicture}}}
\newcommand{\mdb}{\lput{:U}{\begin{pspicture}[shift=0](0,0)(0,0)
\psline[linecolor=lightgray,arrows=->,arrowscale=1.7](3.2pt,0)(3.4pt,0)\end{pspicture}}}

\begin{document}

\author{Dongseok KIM}

\address{Department of Mathematics \\ Kyungpook National University \\ Taegu 702-201 Korea}
\email{dongseok@knu.ac.kr}

\title[Jones-Wenzl idempotents]
{Jones-Wenzl idempotents For Rank 2 Simple Lie algebras}

\begin{abstract}
Temperley-Lieb algebras have been generalized to web spaces for rank
2 simple Lie algebras. Using these webs, we find a complete
description of the Jones-Wenzl idempotents for the quantum
$\mathfrak{sl}(3)$ and $\mathfrak{sp}(4)$ by single clasp
expansions. We discuss applications of these expansions.
\end{abstract}

\maketitle
\section{Introduction}

After the discovery of the Jones polynomial~\cite{Jones:subfactor,
Jones:poly}, its generalizations have been studied in many different
ways. Using the quantum $\mathfrak{sl}(2)$ representation theory,
the Jones polynomial can be seen as a polynomial invariant of a
colored link whose components are colored by the two dimensional
vector representation of the quantum $\mathfrak{sl}(2)$. By using
all irreducible representations of the quantum $\mathfrak{sl}(2)$,
one can find the colored Jones polynomial and it has been
extensively studied \cite{FK:canonical, Khovanov:colored,
kirbyMelvin:witten, Murakami:coloredjones, Vybornov:yang,
Witten:pathint}.

The other direction is to use the representation theory of other
complex simple Lie algebras from the original work of Reshetikhin
and Turaev \cite{RT:1,RT:2}. These quantized simple Lie algebras
invariants can be found by using the Jones-Wenzl idempotents and
fundamental representations. In this philosophy, Kuperberg
introduced web spaces of simple Lie algebras of rank $2$,
$\mathfrak{sl}(3)$, $\mathfrak{sp}(4)$ and $G_2$ as generalizations
of Temperley-Lieb algebras corresponding to $\mathfrak{sl}(2)$
\cite{Kuperberg:spiders}. Then he successively generalized the
result for $\mathfrak{sl}(2)$ \cite{RTW:gottingen1932} that the
dimension of the invariant subspace of the tensor of irreducible
representations of the quantum $\mathfrak{sl}(3)$ and
$\mathfrak{sp}(4)$ is equal to the dimension of web spaces of the
given boundary with respect to the relations in
Figure~\ref{relations} and Figure~\ref{b2rel} respectively
\cite{Kuperberg:spiders}. But there was no immediate generalization
to other Lie algebras until new results for
$\mathfrak{so}(7)$~\cite{math.QA/0601209} and $\mathfrak{sl}(4)$
\cite{KK:sl4}. The quantum $\mathfrak{sl}(3)$ invariants have many
interesting results \cite{chbili, chbili:qm, Khovanov:sl3,
KK:notdual, OY:quantum, Sokolov} also have been generalized to the
quantum $\mathfrak{sl}(n)$ \cite{JK:sln, KR:factor, MOY:Homfly,
Sikora:sln, yokota:skeinforn}. An excellent review can be found in
\cite{SL:review}.

Ohtsuki and Yamada generalized Jones-Wenzl idempotents (these were
called \emph{magic weaving elements}) for the quantum $\mathfrak
{sl}(3)$ web spaces by taking the expansions in Proposition
\ref{a2exppro} and \ref{a2abexpprop1} as a definition of clasps
\cite{OY:quantum}. On the other hand, Kuperberg abstractly proved
the existence of generalized Jones-Wenzl idempotents for other
simple Lie algebras of rank $2$, he called \emph{clasps}
\cite{Kuperberg:spiders}. In the recursive formula shown in
Figure~\ref{a12exp1}, the resulting webs have two (one with one
clasp) clasps, thus it is called a \emph{double clasps expansion} of
the clasp of weight $n$. There is an expansion for which each
resulting web has just one clasp as depicted in
Figure~\ref{a1exppic}~\cite{FK:canonical}. We called it a \emph{
single clasp expansion} of the clasp of weight $n$. These expansions
are very powerful tools for graphical calculus \cite{FK:canonical,
kim:thesis, Vybornov:yang}. We provide single clasp expansions of
all quantum $\mathfrak{sl}(3)$ clasps together with double,
quadruple clasps expansions of all quantum $\mathfrak{sl}(3)$
clasps. We also find single and double clasp expansions of some
quantum $\mathfrak{sp}(4)$ clasps.

Using expansions of clasps, Lickorish first found a quantum
$\mathfrak{sl}(2)$ invariants of 3-manifolds \cite{Lickorish:su,
Lickorish:skein}. Ohtsuki and Yamada did for the quantum $\mathfrak
{sl}(3)$ \cite{OY:quantum} and Yokota found for the quantum
$\mathfrak{sl}(n)$ \cite{yokota:skeinforn}. For applications of
single clasp expansions, first we provide a criterion which
determines the periodicity of a link using colored quantum
$\mathfrak{sl}(3)$ and $\mathfrak{sp}(4)$ link invariants. We
discuss a generalization of $3j, 6j$ symbols for the quantum
$\mathfrak{sl}(3)$ representation theory. At last, we review how
$\mathfrak{sl}(3)$ invariants extend for a special class of graphs.

The outline of this paper is as follows. In section~\ref{review}, we
review the original Jones-Wenzl idempotents and provide some
algebraic background of the representation theory of
$\mathfrak{sl}(3)$ and $\mathfrak{sp}(4)$. We provide single clasp
expansions of all clasps for the quantum $\mathfrak{sl}(3)$ in
section~\ref{singlesl3}. In section~\ref{singlesp4} we study single
clasp expansions of some clasps for the quantum $\mathfrak{sp}(4)$.
In section~\ref{application}, we will discuss some applications of
the quantum $\mathfrak{sl}(3)$ clasps and their single clasp
expansions. In section~\ref{lemmas}, we prove two key lemmas.

A part of the article is from the author's Ph. D. thesis. Precisely,
section~\ref{singlesl3} and \ref{lemmas} are from \cite[section
2.3]{kim:thesis} and section~\ref{singlesp4} is from \cite[section
2.4]{kim:thesis}.  \vskip 1cm

 \noindent{\bf Acknowledgements}
The author would like to thank Greg Kuperberg for introducing the
subject and advising the thesis, Mikhail Khovanov, Jaeun Lee and
Myungsoo Seo for their attention to this work. Also, the referee has
been very helpful and critical during refereing and revising. The
\TeX\, macro package PSTricks~\cite{PSTricks} was essential for
typesetting the equations and figures. The author was supported in
part by KRF Grant M02-2004-000-20044-0.

\section{Jones-Wenzl
idempotents and algebraic back ground}\label{review}

For standard terms and notations for representation theory, we refer
to~\cite{FultonHarris:gtm}.

\subsection{Jones-Wenzl idempotents}\label{jwback}

An explicit algebraic definition of Jones-Wenzl idempotents can be
found in~\cite{FK:canonical}. We will recall a definition of
Jones-Wenzl idempotents which can be generalized for other simple
Lie algebras. Let $T_n$ be the $n$-th Temperley-Lieb algebra with
generators, $1, e_1,e_2,\ldots , e_{n-1}$, and relations,

\begin{align*}
e_i^2&=-(q^{\frac{1}{2}}+q^{-\frac{1}{2}})e_i,\\
e_ie_j&=e_je_i \hskip 2cm \mathrm{if}\hskip .2cm |i-j|\ge 2,\\
e_i&=e_ie_{i\pm 1}e_i.
\end{align*}
For each $n$, the algebra $T_n$ has an idempotent $f_n$ such that
$f_nx=xf_n=\epsilon(x)f_n$ for all $x\in T_n$, where $\epsilon$ is
an augmentation. These idempotents were first discovered by
Jones~\cite{Jones:subfactor} and Wenzl~\cite{Wenzl:Proj}. They found
a recursive formula:

$$f_{n}=f_{n-1}+\frac{[n-1]}{[n]}f_{n-1}e_{n-1}f_{n-1},$$
as illustrated in Figure \ref{a12exp1} where we use a rectangular
box to represent $f_n$ and the quantum integers are defined as
$$
[n]=\frac{q^{\frac{n}{2}}-q^{-\frac{n}{2}}}{q^{\frac{1}{2}}-q^{-\frac{1}{2}}}.
$$ Thus, they are named \emph{Jones-Wenzl
idempotents(projectors)}. It has the following properties 1) it is
an idempotent 2) $f_ne_i=0=e_if_n$ where $e_i$ is a U-turn from the
$i$-th to the $i+1$-th string as shown in Figure \ref{claspaxiom}.
The second property is called \emph{the annihilation axiom}. We will
discuss the importance of Jones-Wenzl idempotents in
section~\ref{webdef}. In Figure~\ref{a1exppic}, $n$ stands for the
number of strings and $``i"$ stands for $i$-th string from the
right. We will use this convention for the rest of the article.

\begin{figure}
$$
\begin{pspicture}[shift=-.8](-.1,-.3)(1.3,1.3)
\rput[t](.5,-.1){$n$}\qline(.5,0)(.5,.4)
\psframe[linecolor=black](0,.4)(1,.6)
\qline(.5,.6)(.5,1)\rput[b](.5,1.1){$n$}
\end{pspicture}
= \begin{pspicture}[shift=-.8](-.3,-.3)(1.5,1.3)
\rput[t](.5,-0.1){$n-1$}\qline(.5,0)(.5,.4)
\psframe[linecolor=black](0,.4)(1,.6)
\qline(.5,.6)(.5,1)\rput[b](.5,1.1){$n-1$} \qline(1.3,0)(1.3,1)
\end{pspicture}
+ \frac{[n-1]}{[n]} \begin{pspicture}[shift=-1.3](-.3,-.3)(1.55,2.3)
\rput[t](.5,-.1){$n-1$} \qline(.5,0)(.5,.4)
\psframe[linecolor=black](0,.4)(1,.6)
\qline(.25,.6)(.25,1.4)\rput[l](.35,1){$n-2$}
\psframe[linecolor=black](0,1.4)(1,1.6) \qline(.5,1.6)(.5,2)
\rput[b](.5,2.1){$n-1$}
\psarc(1,.6){.2}{0}{180}\qline(1.2,0)(1.2,.6)
\psarc(1,1.4){.2}{180}{0}\qline(1.2,1.4)(1.2,2)
\end{pspicture}
$$
\caption{A double clasps expansion of the clasp of weight $n$.}
\label{a12exp1}
\end{figure}

\begin{figure}
$$
\begin{pspicture}[shift=-.4](-.3,0)(2.5,1)
\rput[r](-0.1,.5){$n$}\qline(0,.5)(.4,.5)
\psframe[linecolor=black](.4,0)(.6,1) \qline(.6,.5)(1.4,.5)
\rput[b](1,.6){$n$} \psframe[linecolor=black](1.4,0)(1.6,1)
\qline(1.6,.5)(2,.5)\rput[l](2.1,.5){$n$}
\end{pspicture}
= \begin{pspicture}[shift=-.4](-.5,0)(1.3,1)
\rput[r](-0.1,.5){$n$}\qline(0,.5)(.4,.5)
\psframe[linecolor=black](.4,0)(.6,1)
\qline(.6,.5)(1,.5)\rput[l](1.1,.5){$n$}
\end{pspicture}
\hskip 1cm ,\hskip 1cm \begin{pspicture}[shift=-.4](-.3,0)(3.3,1)
\qline(0,.5)(.4,.5)\rput[r](-.1,.5){$n$}
\psframe[linecolor=black](.4,0)(.6,1)
\qline(.6,.833)(1,.833)\rput[l](1.1,.833){$k$}
\psarc(.6,.5){.167}{-90}{90}
\qline(.6,.167)(1,.167)\rput[l](1.1,.167){$n-k-2$}
\end{pspicture} = 0
$$
\caption{Properties of the Jones-Wenzl idempotents.}
\label{claspaxiom}
\end{figure}

\begin{figure}
$$
\begin{pspicture}[shift=-.6](-.1,-.3)(2.3,2.3)
\rput[t](1,-.1){$n$}\qline(1,0)(1,.4)
\psframe[linecolor=black](0,.4)(2,.6)
\qline(1,.6)(1,2)\rput[b](1,2.1){$n$}
\end{pspicture}
=  \begin{pspicture}[shift=-.6](-.1,-.3)(2.3,2.3)
\rput[t](.8,-.1){$n-1$}\qline(.8,0)(.8,.4)\qline(2,0)(2,2)
\psframe[linecolor=black](0,.4)(1.6,.6)
\qline(.8,.6)(.8,2)\rput[b](.8,2.1){$n-1$}
\end{pspicture} + \sum_{i=2}^{n} \frac{[n+1-i]}{[n]}
\begin{pspicture}[shift=-.6](-.4,-.3)(2.35,2.5) \rput[t](1,-.1){$n-1$}
\qline(1,0)(1,.4) \psframe[linecolor=black](0,.4)(2,.6)
\qline(.2,.6)(.2,2) \qline(.4,.6)(.4,2) \qline(.6,.6)(.6,2)
\rput(0.75,2.3){$``i"$} \rput(2.18,2.3){$``1"$}
\qline(2.2,2)(2.2,1.5) \qline(2.2,1.5)(1.5,1) \qline(1.5,1)(1.5,.6)
\qline(2,2)(2,1.5) \qline(2,1.5)(1.3,1) \qline(1.3,1)(1.3,.6)
\qline(1.8,2)(1.8,1.5) \qline(1.8,1.5)(1.1,1) \qline(1.1,1)(1.1,.6)
\psarc(1.9,.6){.3}{0}{180}\qline(2.2,0)(2.2,.6)
\psarc(1,1.7){.2}{180}{0}\qline(1.2,1.7)(1.2,2) \qline(.8,1.7)(.8,2)
\end{pspicture}
$$
\caption{A single clasp expansion of the clasp of weight $n$.}
\label{a1exppic}
\end{figure}

\subsection{The representation theory of
$\mathfrak{sl}(3)$}\label{sl3back}

The Lie algebra $\mathfrak{sl}(3)$ is the set of all $3\times 3$
complex matrices with trace zero, which is an $8$ dimensional vector
space with the Lie bracket. Let $\lambda_i$ be a fundamental
dominant weight of $\mathfrak{sl}(3)$, $i=1, 2$. All finite
dimensional irreducible representation of $\mathfrak{sl}(3)$ are
determined by its highest weight $\lambda=a\lambda_1+b\lambda_2$,
denoted by $V_{\lambda}$ where $a$ and $b$ are all nonnegative
integers. We will abbreviate $V_{a\lambda_1+b\lambda_2}$ by $
V(a,b)$. The dimension and the quantum dimension of the fundamental
representation $V_{\lambda_1}\cong (V_{\lambda_2})^{*}$ of
$\mathfrak{sl}(3)$ are $3, [3]$. The weight space of a fundamental
representation $V(1,0)$ is $[1,0]$, $[-1,1]$ and $[0, -1]$. The
weight space of a fundamental representation $V(0,1)$ is $[0, 1]$,
$[1, -1]$ and $ [-1, 0]$. Thus, one can easily find the following
decomposition formula for a tensor product of a fundamental
representation and an irreducible representation,

\begin{align*}
V_{\lambda_1}\otimes V_{a\lambda_1+b\lambda_2} &\cong
V_{(a+1)\lambda_1+b\lambda_2}
\oplus V_{(a-1)\lambda_1+(b+1)\lambda_2} \oplus  V_{a\lambda_1+(b-1)\lambda_2},\\
V_{\lambda_2}\otimes V_{a\lambda_1+b\lambda_2} &\cong
V_{a\lambda_1+(b+1)\lambda_2} \oplus
V_{(a+1)\lambda_1+(b-1)\lambda_2} \oplus
V_{(a-1)\lambda_1+b\lambda_2},
\end{align*}
with a standard reflection rule, a refined version of the Brauer's
theorem~\cite[pp.142]{Humphreys:gtm}. Using these tensor rules, one
can find the following lemma.

\begin{lem} For integers $a , b\ge 1$,
$$\mathrm{dim}(\mathrm{Inv}(V_{\lambda_1}^{\otimes a}\otimes
V_{\lambda_2}^{\otimes b}\otimes V_{(b-1)\lambda_1+a\lambda_2}))=
ab. $$ \label{lema2dim}
\end{lem}

To compare the weight of cut paths and clasps, we recall the usual
partial ordering of the weight lattice of ${\mathfrak sl}(3)$ as
\begin{eqnarray*}
a\lambda_1 + b\lambda_2 & \succ & (a+1)\lambda_1 + (b-2)\lambda_2, \\
a\lambda_1 + b\lambda_2 & \succ & (a-2)\lambda_1 + (b+1)\lambda_2.
\end{eqnarray*}

\subsection{The representation theory of
$\mathfrak{sp}(4)$}\label{sp4back}

The Lie algebra $\mathfrak{sp}(4)$ is the set of all $4\times 4$
complex matrices of the following form,
$$\left[\begin{matrix} A & B \\ C & -^tA \end{matrix}\right], ~~~~
\mathrm{where} ~ ^tB=B, ^tC=C$$ which is a $10$ dimensional vector
space with the Lie bracket, where $A, B$ and $C$ are $2\times 2$
matrices. Let $\lambda_i$ be a fundamental dominant weight of
$\mathfrak{sp}(4)$, $i=1, 2$. All finite dimensional irreducible
representation of $\mathfrak{sp}(4)$ are determined by its highest
weight $\lambda=a\lambda_1+b\lambda_2$, denoted by $V_{\lambda}$
where $a$ and $b$ are all nonnegative integers. We will abbreviate
$V_{a\lambda_1+b\lambda_2}$ by $ V(a,b)$. The dimension and the
quantum dimension of the fundamental representation
$V_{\lambda_1}(V_{\lambda_2})$ of $\mathfrak{sp}(4)$ are $4, [4] (5,
[5]$, respectively). The weight space of a fundamental
representation $V(1,0)$ is $[1,0]$, $[-1,1]$, $[1,-1]$ and $[-1,
0]$. The weight space of a fundamental representation $V(0,1)$ is
$[0, 1]$, $[0, -1]$, $[2,-1]$, $[-2,1]$ and $ [0, 0]$. Thus, one can
easily find the following decomposition formula for a tensor product
of a fundamental representation and an irreducible representation,

\begin{align*}
V_{\lambda_1}\otimes V_{a\lambda_1+b\lambda_2} &\cong
V_{(a+1)\lambda_1+b\lambda_2}
\oplus V_{(a-1)\lambda_1+(b+1)\lambda_2} \oplus
V_{(a+1)\lambda_1+(b-1)\lambda_2} \oplus V_{(a-1)\lambda_1+b\lambda_2},\\
V_{\lambda_2}\otimes V_{a\lambda_1+b\lambda_2} &\cong
V_{a\lambda_1+(b+1)\lambda_2} \oplus V_{a\lambda_1+(b-1)\lambda_2}
\oplus V_{(a-2)\lambda_1+(b+1)\lambda_2}  \oplus
V_{(a+2)\lambda_1+(b-1)\lambda_2} \oplus V_{a\lambda_1+b\lambda_2},
\end{align*}
with a similar reflection rule. Using these tensor rules, one can
find the following two lemmas.

\begin{lem} For a positive integer $n$,
$$\mathrm{dim}(\mathrm{Inv}(V_{\lambda_1}^{\otimes n+1}\otimes
V_{(n-1)\lambda_1}))= \frac{n(n+1)}{2}.$$ \label{lemb2dim1}
\end{lem}

\begin{lem} For a positive integer $n$,
$$\mathrm{dim}(\mathrm{Inv}(V_{\lambda_2}^{\otimes n+1}\otimes
V_{(n-1)\lambda_2}))= \frac{n(n+1)}{2}.$$ \label{lemb2dim2}
\end{lem}

There is a natural partial ordering of the $\mathfrak{sp}(4)$ weight
lattice given by

\begin{eqnarray}
a\lambda_1 + b\lambda_2 &\succ &(a-2)\lambda_1 + (b+1)\lambda_2, \nonumber\\
a\lambda_1 + b\lambda_2 &\succ &(a+2)\lambda_1 + (b-2)\lambda_2.
\nonumber
\end{eqnarray}

\subsection{Invariant vector spaces and Web spaces}
\label{webdef}

In this subsection, we briefly review the web spaces, full details
can be found in \cite{Kuperberg:spiders}. Let $V_i$ be an
irreducible representation of complex simple Lie algebras
$\mathfrak{g}$. One of classical invariant problems is to
characterize the vector space of invariant tensors
$$\Inv(V_1\tensor V_2\tensor\ldots\tensor V_n),$$
together with algebraic structures such as tensor products, cyclic
permutations and contractions. The dimension of such a vector space
is given by Cartan-Weyl character theory; $\dim(\Inv(V_1\tensor
V_2\tensor\ldots\tensor V_n))$ is the number of copies of the
trivial representation in the decomposition of $V_1\tensor
V_2\tensor\ldots\tensor V_n$ into irreducible representations. For
this algebraic space, we look for a geometric counterpart which can
preserve the algebraic structure of the invariant spaces. The
discovery of quantum groups opens the door for the link between
invariant spaces and topological invariants of links and manifolds.
For quantum $\mathfrak{sl}(2)$, the dimension of the invariant
spaces\ of $V_1^{\otimes 2n}$ is the dimension of the $n-$th
Temperley-Lieb algebra as a vector space which is generated by chord
diagrams with $2n$ marked points on the boundary of the disk where
$V_1$ is the vector representation of $\mathfrak{sl}(2)$. In
particular, this space is free, $i.e.,$ there is no relation between
chord diagrams. To represent any irreducible representations other
than the vector representation, we use Jones-Wenzl idempotents as we
described in section~\ref{jwback}. Then all webs in the web space of
a tensor of irreducible representations $V_{i_1}\tensor
V_{i_2}\tensor\ldots\tensor V_{i_n}$ can be obtained from webs in
the web space of $V_{1}^{\otimes \sum_k i_k}$ and by attaching
Jones-Wenzl idempotents of weight $i_k$ along the boundary (some
webs become zero by the annihilation axiom, no longer a basis for
web space and the other nonzero webs are called \emph{basis webs}),
where $V_i$ is the irreducible representation of the quantum
$\mathfrak{sl}(2)$ of highest weight $i$ and $k=1, 2, \ldots, n$.
For example \cite{Kuperberg:spiders}, the web
$$
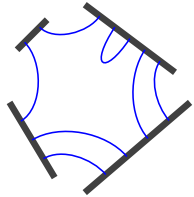
\begin{pspicture}[shift=-1.2](0,0)(2.4,2.5)
\psline[linecolor=darkgray,linewidth=2.5pt](0,1.2)(.6,.2)
\psline[linecolor=darkgray,linewidth=2.5pt](.1,1.9)(.5,2.3)
\psline[linecolor=darkgray,linewidth=2.5pt](1,2.5)(2.2,1.6)
\psline[linecolor=darkgray,linewidth=2.5pt](1,0)(2.4,1.2)
\pnode(1.28,0.24){a1}\pnode(1.56,0.48){a2}\pnode(1.84,0.72){a3}
\pnode(2.12,0.96){a4}
\pnode(2,1.75){b1}\pnode(1.8,1.9){b2}\pnode(1.6,2.05){b3}
\pnode(1.4,2.2){b4}\pnode(1.2,2.35){b5}
\pnode(.4,2.2){c1}\pnode(.2,2.0){c2}
\pnode(.15,.95){d1}\pnode(.3,.7){d2}\pnode(.45,.45){d3}
\nccurve[angleA=130.6,angleB=236.3]{a4}{b1}
\nccurve[angleA=130.6,angleB=236.3]{a3}{b2}
\nccurve[angleA=236.3,angleB=236.3,ncurv=5]{b3}{b4}
\nccurve[angleA=236.3,angleB=-45]{b5}{c1}
\nccurve[angleA=-45,angleB=31]{c2}{d1}
\nccurve[angleA=31,angleB=130.6]{d2}{a2}
\nccurve[angleA=31,angleB=130.6]{d3}{a1}
\end{pspicture}
$$
is not a basis web of $V_{2}\tensor V_{3}\tensor V_4\tensor V_{5}$,
which instead has basis
$$
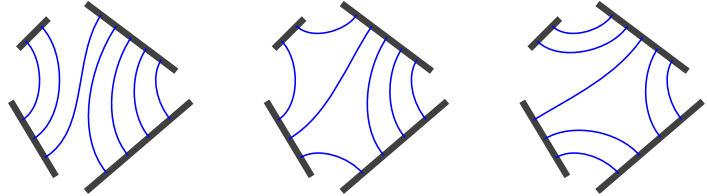
\begin{pspicture}[shift=-1.2](0,0)(2.4,2.5)
\psline[linecolor=darkgray,linewidth=2.5pt](0,1.2)(.6,.2)
\psline[linecolor=darkgray,linewidth=2.5pt](.1,1.9)(.5,2.3)
\psline[linecolor=darkgray,linewidth=2.5pt](1,2.5)(2.2,1.6)
\psline[linecolor=darkgray,linewidth=2.5pt](1,0)(2.4,1.2)
\pnode(1.28,0.24){a1}\pnode(1.56,0.48){a2}\pnode(1.84,0.72){a3}
\pnode(2.12,0.96){a4}
\pnode(2,1.75){b1}\pnode(1.8,1.9){b2}\pnode(1.6,2.05){b3}
\pnode(1.4,2.2){b4}\pnode(1.2,2.35){b5}
\pnode(.4,2.2){c1}\pnode(.2,2.0){c2}
\pnode(.15,.95){d1}\pnode(.3,.7){d2}\pnode(.45,.45){d3}
\nccurve[angleA=130.6,angleB=236.3]{a1}{b4}
\nccurve[angleA=130.6,angleB=236.3]{a2}{b3}
\nccurve[angleA=130.6,angleB=236.3]{a3}{b2}
\nccurve[angleA=130.6,angleB=236.3]{a4}{b1}
\nccurve[angleA=236.3,angleB=31]{b5}{d3}
\nccurve[angleA=-45,angleB=31]{c1}{d2}
\nccurve[angleA=-45,angleB=31]{c2}{d1}
\end{pspicture}
\hspace{1cm} \begin{pspicture}[shift=-1.2](0,0)(2.4,2.5)
\psline[linecolor=darkgray,linewidth=2.5pt](0,1.2)(.6,.2)
\psline[linecolor=darkgray,linewidth=2.5pt](.1,1.9)(.5,2.3)
\psline[linecolor=darkgray,linewidth=2.5pt](1,2.5)(2.2,1.6)
\psline[linecolor=darkgray,linewidth=2.5pt](1,0)(2.4,1.2)
\pnode(1.28,0.24){a1}\pnode(1.56,0.48){a2}\pnode(1.84,0.72){a3}
\pnode(2.12,0.96){a4}
\pnode(2,1.75){b1}\pnode(1.8,1.9){b2}\pnode(1.6,2.05){b3}
\pnode(1.4,2.2){b4}\pnode(1.2,2.35){b5}
\pnode(.4,2.2){c1}\pnode(.2,2.0){c2}
\pnode(.15,.95){d1}\pnode(.3,.7){d2}\pnode(.45,.45){d3}
\nccurve[angleA=130.6,angleB=31]{a1}{d3}
\nccurve[angleA=130.6,angleB=236.3]{a2}{b3}
\nccurve[angleA=130.6,angleB=236.3]{a3}{b2}
\nccurve[angleA=130.6,angleB=236.3]{a4}{b1}
\nccurve[angleA=236.3,angleB=31]{b4}{d2}
\nccurve[angleA=236.3,angleB=-45]{b5}{c1}
\nccurve[angleA=-45,angleB=31]{c2}{d1}
\end{pspicture}
\hspace{1cm} \begin{pspicture}[shift=-1.2](0,0)(2.4,2.5)
\psline[linecolor=darkgray,linewidth=2.5pt](0,1.2)(.6,.2)
\psline[linecolor=darkgray,linewidth=2.5pt](.1,1.9)(.5,2.3)
\psline[linecolor=darkgray,linewidth=2.5pt](1,2.5)(2.2,1.6)
\psline[linecolor=darkgray,linewidth=2.5pt](1,0)(2.4,1.2)
\pnode(1.28,0.24){a1}\pnode(1.56,0.48){a2}\pnode(1.84,0.72){a3}
\pnode(2.12,0.96){a4}
\pnode(2,1.75){b1}\pnode(1.8,1.9){b2}\pnode(1.6,2.05){b3}
\pnode(1.4,2.2){b4}\pnode(1.2,2.35){b5}
\pnode(.4,2.2){c1}\pnode(.2,2.0){c2}
\pnode(.15,.95){d1}\pnode(.3,.7){d2}\pnode(.45,.45){d3}
\nccurve[angleA=130.6,angleB=31]{a1}{d3}
\nccurve[angleA=130.6,angleB=31]{a2}{d2}
\nccurve[angleA=130.6,angleB=236.3]{a3}{b2}
\nccurve[angleA=130.6,angleB=236.3]{a4}{b1}
\nccurve[angleA=236.3,angleB=31]{b3}{d1}
\nccurve[angleA=236.3,angleB=-45]{b4}{c2}
\nccurve[angleA=236.3,angleB=-45]{b5}{c1}
\end{pspicture}
$$
where the Jones-Wenzl idempotents were presented by the thick gray
lines instead of boxes.

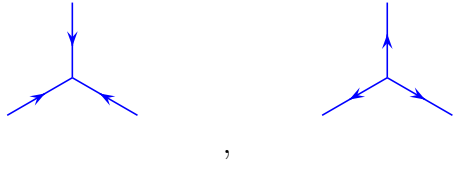
\begin{figure}
$$
\begin{pspicture}[shift=-.9](-1,-1)(1,1) \rput(0,1){\rnode{a1}{$$}}
\rput(-.866,-.5){\rnode{a2}{$$}} \rput(.866,-.5){\rnode{a3}{$$}}
\rput(0,0){\rnode{b1}{$$}} \ncline{a1}{b1}\middlearrow
\ncline{a2}{b1}\middlearrow \ncline{a3}{b1}\middlearrow
\end{pspicture}
\hskip 1cm , \hskip 1cm \begin{pspicture}[shift=-.9](-1,-1)(1,1)
\rput(0,1){\rnode{a1}{$$}} \rput(-.866,-.5){\rnode{a2}{$$}}
\rput(.866,-.5){\rnode{a3}{$$}} \rput(0,0){\rnode{b1}{$$}}
\ncline{b1}{a1}\middlearrow \ncline{b1}{a2}\middlearrow
\ncline{b1}{a3}\middlearrow
\end{pspicture}
$$
\caption{Generators of the quantum $\mathfrak{sl}(3)$ web space.}
\label{generator}
\end{figure}

\begin{figure}
\begin{eqnarray}
\begin{pspicture}[shift=-.4](-.6,-.5)(.6,.5)
\pscircle(0,0){.4}\psline[arrows=->,arrowscale=1.5](.1,.4)(.11,.4)
\end{pspicture}
& = & [3] \label{a2defr21}
\\  \begin{pspicture}[shift=-.7](0,.8)(0,.8)
 \begin{pspicture}[shift=-.7](-1.5,-.8)(1.5,.8) \pnode(.4;180){a2}
\pnode(.4;0){a3}
\qline(-1.2,0)(-.4,0)\psline[arrowscale=1.5]{->}(-.7,0)(-.9,0)
\nccurve[angleA=90,angleB=90,nodesep=1pt]{a2}{a3}\middlearrow
\nccurve[angleA=-90,angleB=-90,nodesep=1pt]{a2}{a3}\middlearrow
\qline(1.2,0)(.4,0)\psline[arrowscale=1.5]{->}(.9,0)(.7,0)
\end{pspicture}\end{pspicture}
& = & - [2] \begin{pspicture}[shift=-.7](-.8,-.8)(.8,.8) \qline(.6,0)(-.6,0)
\psline[arrowscale=1.5]{->}(.1,0)(-.1,0)
\end{pspicture}  \label{a2defr22}
\\
\begin{pspicture}[shift=-1](-1.1,-1.1)(1.1,1.1) \qline(1,1)(.5,.5)
\qline(-1,1)(-.5,.5) \qline(-1,-1)(-.5,-.5) \qline(1,-1)(.5,-.5)
\qline(.5,.5)(-.5,.5) \qline(-.5,.5)(-.5,-.5)
\qline(-.5,-.5)(.5,-.5) \qline(.5,-.5)(.5,.5)
\psline[arrowscale=1.5]{<-}(.1,.5)(-.1,.5)
\psline[arrowscale=1.5]{<-}(-.1,-.5)(.1,-.5)
\psline[arrowscale=1.5]{<-}(.5,.1)(.5,-.1)
\psline[arrowscale=1.5]{<-}(-.5,-.1)(-.5,.1)
\psline[arrowscale=1.5]{->}(.85,.85)(.65,.65)
\psline[arrowscale=1.5]{<-}(-.85,.85)(-.65,.65)
\psline[arrowscale=1.5]{->}(-.85,-.85)(-.65,-.65)
\psline[arrowscale=1.5]{<-}(.85,-.85)(.65,-.65)
\end{pspicture}
&= &   \begin{pspicture}[shift=-1](0,1.1)(0,1.1)\begin{pspicture}[shift=-1](-1.1,-1.1)(1.1,1.1)
\pnode(1;45){a1}\pnode(1;135){a2}\pnode(1;225){a3}\pnode(1;315){a4}
\nccurve[angleA=225,angleB=315]{a1}{a2}\middlearrow
\nccurve[angleA=45,angleB=135]{a3}{a4}\middlearrow
\end{pspicture}\end{pspicture} +  \begin{pspicture}[shift=-1](0,1.1)(0,1.1)\begin{pspicture}[shift=-1](-1.1,-1.1)(1.1,1.1)
\pnode(1;45){a1}\pnode(1;135){a2}\pnode(1;225){a3}\pnode(1;315){a4}
\nccurve[angleA=225,angleB=135]{a1}{a4}\middlearrow
\nccurve[angleA=45,angleB=315]{a3}{a2}\middlearrow
\end{pspicture} \end{pspicture}
\label{a2defr23}
\end{eqnarray}
\caption{Complete relations of the quantum $\mathfrak{sl}(3)$ web
space.} \label{relations}
\end{figure}

A first generalization of Temperley-Lieb algebras was made for
simple Lie algebras of rank 2, $\mathfrak{sl}(3)$,
$\mathfrak{sp}(4)$ and $G_2$ \cite{Kuperberg:spiders}. Each diagrams
appears in a geometric counterpart of the invariant vectors is
called a web, precisely a directed and weighted cubic planar graph.
Unfortunately, some of webs are no longer linearly independent for
simple Lie algebra other than $\mathfrak{sl}(2)$. For example, we
look at the web space of $\mathfrak{sl}(3)$ representations. Let
$V_{\lambda_1}$ be the vector representation of the quantum
$\mathfrak{sl}(3)$ and $V_{\lambda_2}$ be the dual representation of
$V_{\lambda_1}$.  The web space of a fixed boundary (a sequence of
$V_{\lambda_1}$ and $V_{\lambda_2}$) is a vector space spanned by
the all webs of the given boundary which is generated by the webs in
Figure \ref{generator} (as inward and outward arrows) modulo by the
subspace spanned by the equation of diagrams which are called a
complete set of the relations,
equations~\ref{a2defr21},~\ref{a2defr22} and \ref{a2defr23} as
illustrated in Figure~\ref{relations}. We have drawn a web in
Figure~\ref{webexa}. We might use the notation $+, -$ for
$V_{\lambda_1}, V_{\lambda_2}$ but it should be clear. For several
reasons, such as the positivity and the integrality
\cite{Le:integral}, we use $-[2]$ in relation~\ref{a2defr22} but one
can use a quantum integer $[2]$ and get an independent result. If
one uses $[2]$, one can rewrite all results in here by multiplying
each trivalent vertex by the complex number $i$.

\begin{figure}
$$
\begin{pspicture}[shift=-1.9](-1.7,-1.7)(1.7,2.2)
\pnode(.8;0){a1}\pnode(.8;60){a2}\pnode(.8;120){a3}
\pnode(.8;180){a4}\pnode(.8;240){a5}\pnode(.8;300){a6}
\rput(1.6;0){\rnode{b1}{$-$}}\pnode(1.6;60){b2}
\rput(1.6;120){\rnode{b3}{$-$}}\rput(1.6;180){\rnode{b4}{$+$}}
\rput(1.6;240){\rnode{b5}{$-$}}\rput(1.6;300){\rnode{b6}{$+$}}
\rput(1.6,1.386){\rnode{c1}{$-$}}\rput(.4,2.078){\rnode{c2}{$-$}}
\ncline{a2}{a1}\middlearrow\ncline{a2}{a3}\middlearrow
\ncline{a2}{b2}\middlearrow\ncline{a4}{a3}\middlearrow
\ncline{a4}{a5}\middlearrow\ncline[nodesepB=3pt]{a4}{b4}\middlearrow
\ncline{a6}{a5}\middlearrow\ncline{a6}{a1}\middlearrow
\ncline{a6}{b6}\middlearrow\ncline[nodesepA=3pt]{b1}{a1}\middlearrow
\ncline{b3}{a3}\middlearrow\ncline{b5}{a5}\middlearrow
\ncline[nodesepA=3pt]{c1}{b2}\middlearrow\ncline{c2}{b2}\middlearrow
\end{pspicture}
$$
\caption{An example of the webs with a boundary $(+ - + - - - -)$.}
\label{webexa}
\end{figure}
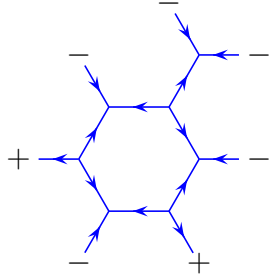

\begin{figure}
$$
\begin{pspicture}[shift=-.8](-.5,-.9)(2.2,.9)
\rput(-.2,-.5){$+$}\rput(-.2,0){$-$}\rput(-.2,.5){$+$}
\qline(0,-.5)(.4,-.5)\qline(0,0)(.4,0)\qline(0,.5)(.4,.5)
\psframe[linecolor=darkgray](.4,-.7)(.6,.7) \qline(.6,0)(1,0)
\pccurve[angleA=0,angleB=240](.6,-.5)(1.25,-.433)
\pccurve[angleA=0,angleB=120](.6,.5)(1.25,.433)
\qline(1,0)(1.25,.433)\qline(1,0)(1.25,-.433)
\qline(1.25,.433)(1.75,.433)\qline(1.25,-.433)(1.75,-.433)
\rput(1.95,.433){$+$}\rput(1.95,-.433){$+$}
\pccurve[arrows=*-*,linestyle=dashed,linecolor=darkgray,
    angleA=0,angleB=0,ncurv=1.7](.5,-.9)(.5,.9)
\end{pspicture} = 0
$$
\caption{An example of the annihilation axiom with a cut path.}
\label{webaxiom2}
\end{figure}
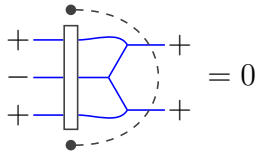

To define the generalization of Jones-Wenzl idempotents,
\emph{clasps}, we first generalize the annihilation axiom for other
web spaces. We need to introduce new concepts: a \emph{cut path} is
a path which is transverse to strings of a web, and the
\emph{weight} of a cut path is the sum of weights of all decorated
strings which intersect with the cut path. For example, the weight
of the clasp as depicted in Figure~\ref{webaxiom2} is
$2V_{\lambda_1}$ abbreviated by $(2,0)$. Then we can generalize the
annihilation axiom as follows: if we attach the clasp to a web which
has a cut path of a weight less than that of the clasp, then it is
zero. Since the weight of the clasp shown in Figure~\ref{webaxiom2}
is $(2,1)$ and there is a cut path of weight $(2,0)$, the web in
Figure~\ref{webaxiom2} is zero by the annihilation axiom. For
$\mathfrak{sl}(3)$, the clasp $\omega$ of weight $(a,b)$ is defined
to be the web in the web space of $V_{\lambda_1}^{\otimes a}\otimes
V_{\lambda_2}^{\otimes b}\otimes (V_{\lambda_1}^*)^{\otimes
a}\otimes (V_{\lambda_2}^*)^{\otimes b}$, say $W$, which satisfies
the annihilation axiom and the idempotent axiom ($\omega^2 =
\omega$). One can see the dimension of the web space of $W$ is one,
$i. e.$, all webs in the web space of $W$ are multiples of $\omega$.
However, the clasp of weight $(a,b)$ is unique by the idempotent
axiom (it is nonzero). An algebraic proof of the existence of clasps
for $\mathfrak{sl}(3)$ and $\mathfrak{sp}(4)$ is
given~\cite{Kuperberg:spiders}. On the other hand, the double clasps
expansion and the quadruple clasps expansion formulae
\cite{OY:quantum} do concretely show the existence of the
$\mathfrak{sl}(3)$ clasp. Using these expansions one can find
Example~\ref{webexa1} (we omit some of arrows on the edges of webs,
but it should be clear).

\begin{exa}\label{webexa1}  The complete
expansions of the clasps of weight $(2,0)$ and $(3,0)$ are
$$
\begin{pspicture}[shift=-.6](-.5,-.7)(.5,.7)
\psline(-.2,.5)(-.2,.1)\psline[arrowscale=1.5]{->}(-.2,.4)(-.2,.2)
\psline(-.2,-.1)(-.2,-.5)\psline[arrowscale=1.5]{->}(-.2,-.2)(-.2,-.4)
\psline(.2,.5)(.2,.1)\psline[arrowscale=1.5]{->}(.2,.4)(.2,.2)
\psline(.2,-.1)(.2,-.5)\psline[arrowscale=1.5]{->}(.2,-.2)(.2,-.4)
\psframe[linecolor=black](-.4,-.1)(.4,.1)
\end{pspicture}
 =   \begin{pspicture}[shift=-.6](-.5,-.7)(.5,.7)
\psline(-.2,.5)(-.2,-.5)\psline[arrowscale=1.5]{->}(-.2,.1)(-.2,-.1)
\psline(.2,.5)(.2,-.5)\psline[arrowscale=1.5]{->}(.2,.1)(.2,-.1)
\end{pspicture}  +
\frac{1}{[2]} \begin{pspicture}[shift=-.4](-.5,-.5)(.5,.5)
\psline(-.5,.5)(0,.2)\psline[arrowscale=1.5]{->}(-.4,.44)(-.2,.32)
\psline(.5,.5)(0,.2)\psline[arrowscale=1.5]{->}(.4,.44)(.2,.32)
\psline(-.5,-.5)(0,-.2)\psline[arrowscale=1.5]{->}(-.4,-.44)(-.2,-.32)
\psline(.5,-.5)(0,-.2)\psline[arrowscale=1.5]{->}(.4,-.44)(.2,-.32)
\psline(0,-.2)(0,.2)\psline[arrowscale=1.5]{->}(0,-.1)(0,.1)
\end{pspicture}
$$
$$ \begin{pspicture}[shift=-.9](-.6,-.7)(.6,1.2)
\psline(-.4,1)(-.4,.1)\psline[arrowscale=1.5]{->}(-.4,.4)(-.4,.2)
\psline(-.4,-.1)(-.4,-.5)\psline[arrowscale=1.5]{->}(-.4,-.2)(-.4,-.4)
\psline(.4,1)(.4,.1)\psline[arrowscale=1.5]{->}(.4,.4)(.4,.2)
\psline(.4,-.1)(.4,-.5)\psline[arrowscale=1.5]{->}(.4,-.2)(.4,-.4)
\psline(0,1)(0,.1)\psline[arrowscale=1.5]{->}(0,.4)(0,.2)
\psline(0,-.1)(0,-.5)\psline[arrowscale=1.5]{->}(0,-.2)(0,-.4)
\psframe[linecolor=black](-.6,-.1)(.6,.1)
\end{pspicture}
 =   \begin{pspicture}[shift=-.9](-.6,-.7)(.6,1.2)
\psline(-.4,1)(-.4,-.5)\psline[arrowscale=1.5]{->}(-.4,.1)(-.4,-.1)
\psline(.4,1)(.4,-.5)\psline[arrowscale=1.5]{->}(.4,.1)(.4,-.1)
\psline(0,1)(0,-.5)\psline[arrowscale=1.5]{->}(0,.1)(0,-.1)
\end{pspicture}  + \frac{[2]}{[3]} \left(
 \begin{pspicture}[shift=-.9](-.6,-.7)(.6,1.2)
\psline(-.4,1)(-.4,-.5)\psline[arrowscale=1.5]{->}(-.4,.1)(-.4,-.1)
\psline(.4,.1)(.4,-.5)\psline[arrowscale=1.5]{->}(.4,-.1)(.4,-.3)
\psline(0,1)(.2,.7)(.2,.4)(0,.1) \psline(.4,1)(.2,.7)(.2,.4)(.4,.1)
\psline(0,.1)(0,-.5)\psline[arrowscale=1.5]{->}(0,-.1)(0,-.3)
\end{pspicture} +
 \begin{pspicture}[shift=-.9](-.6,-.7)(.6,1.2)
\psline(.4,1)(.4,-.5)\psline[arrowscale=1.5]{->}(.4,.1)(.4,-.1)
\psline(-.4,.1)(-.4,-.5)\psline[arrowscale=1.5]{->}(-.4,-.1)(-.4,-.3)
\psline(0,1)(-.2,.7)(-.2,.4)(0,.1)
\psline(-.4,1)(-.2,.7)(-.2,.4)(-.4,.1)
\psline(0,.1)(0,-.5)\psline[arrowscale=1.5]{->}(0,-.1)(0,-.3)
\end{pspicture}
\right)+  \frac{[1]}{[3]}\left( \begin{pspicture}[shift=-.9](-.6,-.7)(.6,1.2)
\psline(-.4,1)(-.4,.25)(-.2,0)(-.2,-.2)(-.4,-.5)
\psline(0,1)(.2,.7)(.2,.5)(-.2,0)(-.2,-.2)(0,-.5)
\psline(.4,1)(.2,.7)(.2,.5)(.4,.25)(.4,-.5)
\psline[arrowscale=1.5]{->}(-.4,.5)(-.4,.3)
\psline[arrowscale=1.5]{->}(.4,.1)(.4,-.1)
\end{pspicture} + \begin{pspicture}[shift=-.9](-.6,-.7)(.6,1.2)
\psline(.4,1)(.4,.25)(.2,0)(.2,-.2)(.4,-.5)
\psline(0,1)(-.2,.7)(-.2,.5)(.2,0)(.2,-.2)(0,-.5)
\psline(-.4,1)(-.2,.7)(-.2,.5)(-.4,.25)(-.4,-.5)
\psline[arrowscale=1.5]{->}(.4,.5)(.4,.3)
\psline[arrowscale=1.5]{->}(-.4,.1)(-.4,-.1)
\end{pspicture}\right) + \frac{1}{[2][3]} \begin{pspicture}[shift=-1](0,.9)(0,.9)\begin{pspicture}[shift=-.8](-.6,-.7)(.6,1.2)
\pccurve[angleA=-90,angleB=240](-.4,1)(0,.5)\middlearrow
\pccurve[angleA=-90,angleB=-30](.4,1)(0,.5)\middlearrow
\pcline(0,1)(0,.5)\middlearrow
\pccurve[angleA=120,angleB=90](0,0)(-.4,-.5)\middlearrow
\pccurve[angleA=30,angleB=90](0,0)(.4,-.5)\middlearrow
\pcline(0,0)(0,-.5)\middlearrow
\end{pspicture}\end{pspicture}
$$
\end{exa}

\section{Single Clasp Expansions for the quantum $\mathfrak{sl}(3)$ clasps}
\label{singlesl3}

First we look at a single clasp expansion of the clasp of weight
$(n,0)$ where the weight $(a,b)$ stands for $a\lambda_1 +b\lambda_2$
in section~\ref{sl3no}. We can easily find a single clasp expansion
of the clasp of weight $(0,n)$ by reversing arrows in the equation
presented in the formula of Proposition \ref{a2exppro}. In
section~\ref{sl3ab}, we find a single clasp expansion of the clasp
of weight $(a,b)$ and double clasps expansions. Kuperberg showed
that for a fixed boundary, all webs of the given boundary are cut
outs from the hexagonal tiling of the plane with the given boundary
\cite{Kuperberg:spiders}.

\subsection{Single clasp expansions of a clasp of weight (n,0)}
\label{sl3no}

First we find a single clasp expansions of a clasp of weight (n,0)
in Proposition~\ref{a2exppro}. It is worth to mention that $i)$ this
expansion can be obtained from a complete expansion (linear
expansions of webs without any clasps) which can be found by using a
double clasps expansion iteratively \cite{OY:quantum} and then
attaching a clasp of weight $(n-1,0)$ to each web in the expansion;
$ii)$ the single clasp expansion in Proposition~\ref{a2exppro} holds
for any $\mathfrak{sl}(n)$ where $n\ge 4$ because $\mathfrak{sl}(3)$
is naturally embedded in $\mathfrak{sl}(n)$. By symmetries, there
are four different single clasp expansions depending on where the
clasp of weight $(n-1,0)$ is located. For equation~\ref{a2exp}, the
clasp is located at the southwest corner, which will be considered
the standard expansion, otherwise, we will state the location of the
clasp.

We demonstrate Proposition~\ref{a2exppro} for $n=2, 3$ directly
using the presentations of the clasps in Example~\ref{webexa1}. For
$n=2$, Proposition~\ref{a2exppro} is identical to the first formula
of Example \ref{webexa1}. For $n=3$, by attaching the clasp of
weight $(2,0)$ to the southwest corner of each web in the second
formula of Example \ref{webexa1},

$$ \begin{pspicture}[shift=-1.1](-.6,-1.2)(.6,1.2)
\psline(-.4,1)(-.4,.1)\psline[arrowscale=1.5]{<-}(-.4,.4)(-.4,.2)
\psline(-.4,-.1)(-.4,-.5)\psline[arrowscale=1.5]{<-}(-.4,-.2)(-.4,-.4)
\psline(.4,1)(.4,.1)\psline[arrowscale=1.5]{<-}(.4,.4)(.4,.2)
\psline(.4,-.1)(.4,-.5)\psline[arrowscale=1.5]{<-}(.4,-.2)(.4,-.4)
\psline(0,1)(0,.1)\psline[arrowscale=1.5]{<-}(0,.4)(0,.2)
\psline(0,-.1)(0,-.5)\psline[arrowscale=1.5]{<-}(0,-.2)(0,-.4)
\psframe[linecolor=black](-.6,-.1)(.6,.1)
\psframe[linecolor=black](-.6,-.7)(.2,-.5) \qline(.4,-1)(.4,-.5)
\psline(-.4,-1)(-.4,-.7)\psline[arrowscale=1.5]{<-}(-.4,-.75)(-.4,-.95)
\psline(0,-1)(0,-.7)\psline[arrowscale=1.5]{<-}(0,-.75)(0,-.95)
\end{pspicture} =   \begin{pspicture}[shift=-1.1](-.6,-1.2)(.6,1.2)
\psline(-.4,1)(-.4,-.5)\psline[arrowscale=1.5]{<-}(-.4,.1)(-.4,-.1)
\psline(.4,1)(.4,-.5)\psline[arrowscale=1.5]{<-}(.4,.1)(.4,-.1)
\psline(0,1)(0,-.5)\psline[arrowscale=1.5]{<-}(0,.1)(0,-.1)
\psframe[linecolor=black](-.6,-.7)(.2,-.5) \qline(.4,-1)(.4,-.5)
\psline(-.4,-1)(-.4,-.7)\psline[arrowscale=1.5]{<-}(-.4,-.75)(-.4,-.95)
\psline(0,-1)(0,-.7)\psline[arrowscale=1.5]{<-}(0,-.75)(0,-.95)
\end{pspicture}  + \frac{[2]}{[3]} \left(
 \begin{pspicture}[shift=-1.1](-.6,-1.2)(.6,1.2)
\psline(-.4,1)(-.4,-.5)\psline[arrowscale=1.5]{<-}(-.4,.1)(-.4,-.1)
\psline(.4,.1)(.4,-.5)\psline[arrowscale=1.5]{<-}(.4,-.1)(.4,-.3)
\psline(0,1)(.2,.7)(.2,.4)(0,.1) \psline(.4,1)(.2,.7)(.2,.4)(.4,.1)
\psline(0,.1)(0,-.5)\psline[arrowscale=1.5]{<-}(0,-.1)(0,-.3)
\psframe[linecolor=black](-.6,-.7)(.2,-.5) \qline(.4,-1)(.4,-.5)
\psline(-.4,-1)(-.4,-.7)\psline[arrowscale=1.5]{<-}(-.4,-.75)(-.4,-.95)
\psline(0,-1)(0,-.7)\psline[arrowscale=1.5]{<-}(0,-.75)(0,-.95)
\end{pspicture} +
 \begin{pspicture}[shift=-1.1](-.6,-1.2)(.6,1.2)
\psline(.4,1)(.4,-.5)\psline[arrowscale=1.5]{<-}(.4,.1)(.4,-.1)
\psline(-.4,.1)(-.4,-.5)\psline[arrowscale=1.5]{<-}(-.4,-.1)(-.4,-.3)
\psline(0,1)(-.2,.7)(-.2,.4)(0,.1)
\psline(-.4,1)(-.2,.7)(-.2,.4)(-.4,.1)
\psline(0,.1)(0,-.5)\psline[arrowscale=1.5]{<-}(0,-.1)(0,-.3)
\psframe[linecolor=black](-.6,-.7)(.2,-.5) \qline(.4,-1)(.4,-.5)
\psline(-.4,-1)(-.4,-.7)\psline[arrowscale=1.5]{<-}(-.4,-.75)(-.4,-.95)
\psline(0,-1)(0,-.7)\psline[arrowscale=1.5]{<-}(0,-.75)(0,-.95)
\end{pspicture}
\right)+  \frac{[1]}{[3]}\left(  \begin{pspicture}[shift=-1.1](-.6,-1.2)(.6,1.2)
\psline(-.4,1)(-.4,.25)(-.2,0)(-.2,-.2)(-.4,-.5)
\psline(0,1)(.2,.7)(.2,.5)(-.2,0)(-.2,-.2)(0,-.5)
\psline(.4,1)(.2,.7)(.2,.5)(.4,.25)(.4,-.5)
\psline[arrowscale=1.5]{<-}(-.4,.5)(-.4,.3)
\psline[arrowscale=1.5]{<-}(.4,.1)(.4,-.1)
\psframe[linecolor=black](-.6,-.7)(.2,-.5) \qline(.4,-1)(.4,-.5)
\psline(-.4,-1)(-.4,-.7)\psline[arrowscale=1.5]{<-}(-.4,-.75)(-.4,-.95)
\psline(0,-1)(0,-.7)\psline[arrowscale=1.5]{<-}(0,-.75)(0,-.95)
\end{pspicture} +  \begin{pspicture}[shift=-1.1](-.6,-1.2)(.6,1.2)
\psline(.4,1)(.4,.25)(.2,0)(.2,-.2)(.4,-.5)
\psline(0,1)(-.2,.7)(-.2,.5)(.2,0)(.2,-.2)(0,-.5)
\psline(-.4,1)(-.2,.7)(-.2,.5)(-.4,.25)(-.4,-.5)
\psline[arrowscale=1.5]{<-}(.4,.5)(.4,.3)
\psline[arrowscale=1.5]{<-}(-.4,.1)(-.4,-.1)
\psframe[linecolor=black](-.6,-.7)(.2,-.5) \qline(.4,-1)(.4,-.5)
\psline(-.4,-1)(-.4,-.7)\psline[arrowscale=1.5]{<-}(-.4,-.75)(-.4,-.95)
\psline(0,-1)(0,-.7)\psline[arrowscale=1.5]{<-}(0,-.75)(0,-.95)
\end{pspicture}\right)
+ \frac{1}{[2][3]}   \begin{pspicture}[shift=-1.1](0,1.2)(0,1.2)\begin{pspicture}[shift=-.9](-.6,-1.2)(.6,1.2)
\pccurve[angleA=-90,angleB=240](0,.5)(-.4,1)\middlearrow
\pccurve[angleA=-90,angleB=-30](0,.5)(.4,1)\middlearrow
\pcline(0,.5)(0,1)\middlearrow
\pccurve[angleA=120,angleB=90](-.4,-.5)(0,0)\middlearrow
\pccurve[angleA=30,angleB=90](.4,-.5)(0,0)\middlearrow
\pcline(0,-.5)(0,0)\middlearrow
\psframe[linecolor=black](-.6,-.7)(.2,-.5) \qline(.4,-1)(.4,-.5)
\psline(-.4,-1)(-.4,-.7)\psline[arrowscale=1.5]{<-}(-.4,-.75)(-.4,-.95)
\psline(0,-1)(0,-.7)\psline[arrowscale=1.5]{<-}(0,-.75)(0,-.95)
\end{pspicture}\end{pspicture}$$

Since $ \begin{pspicture}[shift=-.9](0,.9)(0,.9) \begin{pspicture}[shift=-.9](-.5,-.7)(.5,1.2)
\pcline(-.2,.1)(0,.5)\middlearrow \pcline(.2,.1)(0,.5)\middlearrow
\pcline(0,1)(0,.5)\middlearrow
\psline(-.2,-.1)(-.2,-.5)\psline[arrowscale=1.5]{->}(-.2,-.4)(-.2,-.2)
\psline(.2,-.1)(.2,-.5)\psline[arrowscale=1.5]{->}(.2,-.4)(.2,-.2)
\psframe[linecolor=black](-.4,-.1)(.4,.1)
\end{pspicture}\end{pspicture} =0$, we find

$$  \begin{pspicture}[shift=-1.1](-.6,-1.2)(.6,1.2)
\psline(-.4,1)(-.4,.1)\psline[arrowscale=1.5]{<-}(-.4,.4)(-.4,.2)
\psline(-.4,-.1)(-.4,-1)\psline[arrowscale=1.5]{<-}(-.4,-.2)(-.4,-.4)
\psline(.4,1)(.4,.1)\psline[arrowscale=1.5]{<-}(.4,.4)(.4,.2)
\psline(.4,-.1)(.4,-1)\psline[arrowscale=1.5]{<-}(.4,-.2)(.4,-.4)
\psline(0,1)(0,.1)\psline[arrowscale=1.5]{<-}(0,.4)(0,.2)
\psline(0,-.1)(0,-1)\psline[arrowscale=1.5]{<-}(0,-.2)(0,-.4)
\psframe[linecolor=black](-.6,-.1)(.6,.1)
\end{pspicture}
 =  \begin{pspicture}[shift=-1.1](-.8,-1.2)(.6,1.2)
\psline(-.4,1)(-.4,-.5)\psline[arrowscale=1.5]{<-}(-.4,.1)(-.4,-.1)
\psline(.4,1)(.4,-.5)\psline[arrowscale=1.5]{<-}(.4,.1)(.4,-.1)
\psline(0,1)(0,-.5)\psline[arrowscale=1.5]{<-}(0,.1)(0,-.1)\psframe[linecolor=black](-.6,-.7)(.2,-.5)
\qline(.4,-1)(.4,-.5)
\psline(-.4,-1)(-.4,-.7)\psline[arrowscale=1.5]{<-}(-.4,-.75)(-.4,-.95)
\psline(0,-1)(0,-.7)\psline[arrowscale=1.5]{<-}(0,-.75)(0,-.95)
\end{pspicture}  + \frac{[2]}{[3]}
  \begin{pspicture}[shift=-1.1](-.8,-1.2)(.6,1.2)
\psline(-.4,1)(-.4,-.5)\psline[arrowscale=1.5]{<-}(-.4,.1)(-.4,-.1)
\psline(.4,.1)(.4,-.5)\psline[arrowscale=1.5]{<-}(.4,-.1)(.4,-.3)
\psline(0,1)(.2,.7)(.2,.4)(0,.1) \psline(.4,1)(.2,.7)(.2,.4)(.4,.1)
\psline(0,.1)(0,-.5)\psline[arrowscale=1.5]{<-}(0,-.1)(0,-.3)
\psframe[linecolor=black](-.6,-.7)(.2,-.5) \qline(.4,-1)(.4,-.5)
\psline(-.4,-1)(-.4,-.7)\psline[arrowscale=1.5]{<-}(-.4,-.75)(-.4,-.95)
\psline(0,-1)(0,-.7)\psline[arrowscale=1.5]{<-}(0,-.75)(0,-.95)
\end{pspicture} +  \frac{[1]}{[3]}  \begin{pspicture}[shift=-1.1](-.8,-1.2)(.6,1.2)
\psline(.4,1)(.4,.25)(.2,0)(.2,-.2)(.4,-.5)
\psline(0,1)(-.2,.7)(-.2,.5)(.2,0)(.2,-.2)(0,-.5)
\psline(-.4,1)(-.2,.7)(-.2,.5)(-.4,.25)(-.4,-.5)
\psline[arrowscale=1.5]{<-}(.4,.5)(.4,.3)
\psline[arrowscale=1.5]{<-}(-.4,.1)(-.4,-.1)
\psframe[linecolor=black](-.6,-.7)(.2,-.5) \qline(.4,-1)(.4,-.5)
\psline(-.4,-1)(-.4,-.7)\psline[arrowscale=1.5]{<-}(-.4,-.75)(-.4,-.95)
\psline(0,-1)(0,-.7)\psline[arrowscale=1.5]{<-}(0,-.75)(0,-.95)
\end{pspicture}
$$

This verifies the $n=3$ case of Proposition~\ref{a2exppro}.

\begin{prop} For a positive integer $n$,
\begin{align}  \begin{pspicture}[shift=-1.2](-.1,-.3)(2.6,2.3)
\rput[t](1,-.1){$n$}\qline(1,0)(1,.4) \rput[c](1,1.2){$\ldots$}
\qline(.1,.6)(.1,2) \qline(.3,.6)(.3,2) \qline(.5,.6)(.5,2)
\qline(1.5,.6)(1.5,2) \qline(1.7,.6)(1.7,2) \qline(1.9,.6)(1.9,2)
\psline[arrowscale=1.5]{->}(.1,1)(.1,1.2)
\psline[arrowscale=1.5]{->}(.3,1)(.3,1.2)
\psline[arrowscale=1.5]{->}(.5,1)(.5,1.2)
\psline[arrowscale=1.5]{->}(1.5,1)(1.5,1.2)
\psline[arrowscale=1.5]{->}(1.7,1)(1.7,1.2)
\psline[arrowscale=1.5]{->}(1.9,1)(1.9,1.2)
\psframe[linecolor=black](0,.4)(2,.6)
\end{pspicture}
=  \begin{pspicture}[shift=-1.2](-.1,-.3)(2.6,2.3)
\rput[t](.9,-.1){$n-1$}\qline(1,0)(1,.4) \rput[c](1,1.2){$\ldots$}
\qline(.1,.6)(.1,2) \qline(.3,.6)(.3,2) \qline(.5,.6)(.5,2)
\qline(1.5,.6)(1.5,2) \qline(1.7,.6)(1.7,2) \qline(2.1,0)(2.1,2)
\psline[arrowscale=1.5]{->}(.1,1)(.1,1.2)
\psline[arrowscale=1.5]{->}(.3,1)(.3,1.2)
\psline[arrowscale=1.5]{->}(.5,1)(.5,1.2)
\psline[arrowscale=1.5]{->}(1.5,1)(1.5,1.2)
\psline[arrowscale=1.5]{->}(1.7,1)(1.7,1.2)
\psline[arrowscale=1.5]{->}(2.1,1)(2.1,1.2)
\psframe[linecolor=black](0,.4)(1.8,.6)
\end{pspicture} + \sum_{i=2}^{n} \frac{[n+1-i]}{[n]}  \begin{pspicture}[shift=-1.3](-.6,-.3)(3.35,2.5)
\rput[t](1,-.1){$n-1$} \rput[b](.4,1.8){$\ldots$} \psline(1,0)(1,.4)
\psline(.1,.6)(.1,2.2) \psline(.7,.6)(.7,2.2)
\psline(.9,1.8)(.9,2.2) \psline(.9,1.8)(1.1,1.6)
\psline(1.1,1.6)(1.1,1.4) \psline(1.1,1.4)(1.3,.9)
\psline(1.3,1.8)(1.1,1.6) \psline(1.3,1.8)(1.3,2.2)
\psline(1.3,.6)(1.3,.9) \psline(1.5,1.4)(1.3,.9)
\psline(1.5,1.4)(1.5,2.2) \psline(1.5,1.4)(1.7,.9)
\psline(2.2,1.4)(2,.9) \psline(2.2,1.4)(2.4,.9)
\psline(2.2,1.4)(2.2,2.2) \psline(2.4,.9)(2.6,1.4)
\psline(2.4,.6)(2.4,.9) \psline(2.6,1.4)(2.8,.9)
\psline(2.6,1.4)(2.6,2.2) \psline(2.8,.9)(3,1.4)
\psline(2.8,.9)(2.8,.6) \psline(3,1.4)(3.2,.9)
\psline(3.2,.9)(3.2,0) \psline[arrowscale=1.5]{->}(.1,1.3)(.1,1.5)
\psline[arrowscale=1.5]{->}(.7,1.3)(.7,1.5)
\psline[arrowscale=1.5]{->}(.9,1.9)(.9,2.1)
\psline[arrowscale=1.5]{->}(1.16,1.25)(1.24,1.05)
\psline[arrowscale=1.5]{->}(1.3,1.9)(1.3,2.1)
\psline[arrowscale=1.5]{->}(1.3,.65)(1.3,.85)
\psline[arrowscale=1.5]{->}(1.44,1.25)(1.36,1.05)
\psline[arrowscale=1.5]{->}(1.5,1.7)(1.5,1.9)
\psline[arrowscale=1.5]{->}(1.56,1.25)(1.64,1.05)
\psline[arrowscale=1.5]{->}(2.14,1.25)(2.06,1.05)
\psline[arrowscale=1.5]{->}(2.26,1.25)(2.34,1.05)
\psline[arrowscale=1.5]{->}(2.2,1.7)(2.2,1.9)
\psline[arrowscale=1.5]{<-}(2.46,1.05)(2.54,1.25)
\psline[arrowscale=1.5]{->}(2.4,.65)(2.4,.85)
\psline[arrowscale=1.5]{->}(2.66,1.25)(2.74,1.05)
\psline[arrowscale=1.5]{->}(2.6,1.7)(2.6,1.9)
\psline[arrowscale=1.5]{<-}(2.86,1.05)(2.94,1.25)
\psline[arrowscale=1.5]{->}(2.8,.65)(2.8,.85)
\psline[arrowscale=1.5]{<-}(3.2,.55)(3.2,.35)
\rput[b](.9,2.4){$``i"$} \rput[b](1.85,1.8){$\ldots$}
\rput[b](2.6,2.4){$``1"$} \psframe[linecolor=black](0,.4)(3,.6)
\end{pspicture} \label{a2exp}
\end{align} \label{a2exppro}
\end{prop}
\begin{proof}
We prove the linear independence of the webs in the right-hand side of the
equation~\ref{a2exp}. Suppose there exists a linear combination of
webs which is zero, let us denote $c_i$ be the coefficient of this
linear combination corresponding to the $i$-th web in the right-hand side of the
equation~\ref{a2exp}.

\begin{align*} 0
=  c_1 \begin{pspicture}[shift=-1.4](-.1,-.3)(2.6,2.7)
\rput[t](.9,-.1){$n-1$}\qline(1,0)(1,.4) \rput[c](1,1.2){$\ldots$}
\qline(.1,.6)(.1,2) \qline(.3,.6)(.3,2) \qline(.5,.6)(.5,2)
\qline(1.5,.6)(1.5,2) \qline(1.7,.6)(1.7,2) \qline(2.1,0)(2.1,2)
\psline[arrowscale=1.5]{->}(.1,1)(.1,1.2)
\psline[arrowscale=1.5]{->}(.3,1)(.3,1.2)
\psline[arrowscale=1.5]{->}(.5,1)(.5,1.2)
\psline[arrowscale=1.5]{->}(1.5,1)(1.5,1.2)
\psline[arrowscale=1.5]{->}(1.7,1)(1.7,1.2)
\psline[arrowscale=1.5]{->}(2.1,1)(2.1,1.2)
\psframe[linecolor=black](0,2)(2.3,2.2)
\psframe[linecolor=black](0,.4)(1.8,.6)
\end{pspicture} + \sum_{j=2}^{n} c_j  \begin{pspicture}[shift=-1.5](-.6,-.3)(3.35,2.9)
\rput[t](1,-.1){$n-1$} \rput[b](.4,1.8){$\ldots$} \psline(1,0)(1,.4)
\psline(.1,.6)(.1,2.2) \psline(.7,.6)(.7,2.2)
\psline(.9,1.8)(.9,2.2) \psline(.9,1.8)(1.1,1.6)
\psline(1.1,1.6)(1.1,1.4) \psline(1.1,1.4)(1.3,.9)
\psline(1.3,1.8)(1.1,1.6) \psline(1.3,1.8)(1.3,2.2)
\psline(1.3,.6)(1.3,.9) \psline(1.5,1.4)(1.3,.9)
\psline(1.5,1.4)(1.5,2.2) \psline(1.5,1.4)(1.7,.9)
\psline(2.2,1.4)(2,.9) \psline(2.2,1.4)(2.4,.9)
\psline(2.2,1.4)(2.2,2.2) \psline(2.4,.9)(2.6,1.4)
\psline(2.4,.6)(2.4,.9) \psline(2.6,1.4)(2.8,.9)
\psline(2.6,1.4)(2.6,2.2) \psline(2.8,.9)(3,1.4)
\psline(2.8,.9)(2.8,.6) \psline(3,1.4)(3.2,.9)
\psline(3.2,.9)(3.2,0) \psline[arrowscale=1.5]{->}(.1,1.3)(.1,1.5)
\psline[arrowscale=1.5]{->}(.7,1.3)(.7,1.5)
\psline[arrowscale=1.5]{->}(.9,1.9)(.9,2.1)
\psline[arrowscale=1.5]{->}(1.16,1.25)(1.24,1.05)
\psline[arrowscale=1.5]{->}(1.3,1.9)(1.3,2.1)
\psline[arrowscale=1.5]{->}(1.3,.65)(1.3,.85)
\psline[arrowscale=1.5]{->}(1.44,1.25)(1.36,1.05)
\psline[arrowscale=1.5]{->}(1.5,1.7)(1.5,1.9)
\psline[arrowscale=1.5]{->}(1.56,1.25)(1.64,1.05)
\psline[arrowscale=1.5]{->}(2.14,1.25)(2.06,1.05)
\psline[arrowscale=1.5]{->}(2.26,1.25)(2.34,1.05)
\psline[arrowscale=1.5]{->}(2.2,1.7)(2.2,1.9)
\psline[arrowscale=1.5]{<-}(2.46,1.05)(2.54,1.25)
\psline[arrowscale=1.5]{->}(2.4,.65)(2.4,.85)
\psline[arrowscale=1.5]{->}(2.66,1.25)(2.74,1.05)
\psline[arrowscale=1.5]{->}(2.6,1.7)(2.6,1.9)
\psline[arrowscale=1.5]{<-}(2.86,1.05)(2.94,1.25)
\psline[arrowscale=1.5]{->}(2.8,.65)(2.8,.85)
\psline[arrowscale=1.5]{<-}(3.2,.55)(3.2,.35)
\rput[b](.9,2.5){$``j"$} \rput[b](1.85,1.8){$\ldots$}
\rput[b](2.6,2.5){$``1"$} \psframe[linecolor=black](0,.4)(3,.6)
\psframe[linecolor=black](0,2.2)(2.8,2.4)
\end{pspicture}
\end{align*}If we attach the clasp of weight $(n,0)$ to
the top of each web in the right-hand side of equation~\ref{a2exp}, the
first web corresponding to the coefficient $c_1$ is nonzero because
it is a cut out from the hexagonal tiling of the plane. All other
remaining webs corresponding to $c_k$ where $k\ge 2$ are zero
because $ \begin{pspicture}[shift=-.9](0,.9)(0,.9)\begin{pspicture}[shift=.4](-.5,-1.2)(.5,.7)
\pcline(0,-.5)(-.2,-.1)\middlearrow
\pcline(0,-.5)(.2,-.1)\middlearrow \pcline(0,-.5)(0,-1)\middlearrow
\psline(-.2,.1)(-.2,.5)\psline[arrowscale=1.5]{<-}(-.2,.4)(-.2,.2)
\psline(.2,.1)(.2,.5)\psline[arrowscale=1.5]{<-}(.2,.4)(.2,.2)
\psframe[linecolor=black](-.4,-.1)(.4,.1)
\end{pspicture} \end{pspicture} =0$. Therefore, $c_1=0$.

\begin{align*} 0
=  c_i \begin{pspicture}[shift=-1.5](-.6,-.3)(3.35,2.9) \rput[t](1,-.1){$n-1$}
\rput[b](.4,1.8){$\ldots$} \psline(1,0)(1,.4) \psline(.1,.6)(.1,2.2)
\psline(.7,.6)(.7,2.2) \psline(.9,1.8)(.9,2.2)
\psline(.9,1.8)(1.1,1.6) \psline(1.1,1.6)(1.1,1.4)
\psline(1.1,1.4)(1.3,.9) \psline(1.3,1.8)(1.1,1.6)
\psline(1.3,1.8)(1.3,2.2) \psline(1.3,.6)(1.3,.9)
\psline(1.5,1.4)(1.3,.9) \psline(1.5,1.4)(1.5,2.2)
\psline(1.5,1.4)(1.7,.9) \psline(2.2,1.4)(2,.9)
\psline(2.2,1.4)(2.4,.9) \psline(2.2,1.4)(2.2,2.2)
\psline(2.4,.9)(2.6,1.4) \psline(2.4,.6)(2.4,.9)
\psline(2.6,1.4)(2.8,.9) \psline(2.6,1.4)(2.6,2.2)
\psline(2.8,.9)(3,1.4) \psline(2.8,.9)(2.8,.6)
\psline(3,1.4)(3.2,.9) \psline(3.2,.9)(3.2,0)
\psline[arrowscale=1.5]{->}(.1,1.3)(.1,1.5)
\psline[arrowscale=1.5]{->}(.7,1.3)(.7,1.5)
\psline[arrowscale=1.5]{->}(.9,1.9)(.9,2.1)
\psline[arrowscale=1.5]{->}(1.16,1.25)(1.24,1.05)
\psline[arrowscale=1.5]{->}(1.3,1.9)(1.3,2.1)
\psline[arrowscale=1.5]{->}(1.3,.65)(1.3,.85)
\psline[arrowscale=1.5]{->}(1.44,1.25)(1.36,1.05)
\psline[arrowscale=1.5]{->}(1.5,1.7)(1.5,1.9)
\psline[arrowscale=1.5]{->}(1.56,1.25)(1.64,1.05)
\psline[arrowscale=1.5]{->}(2.14,1.25)(2.06,1.05)
\psline[arrowscale=1.5]{->}(2.26,1.25)(2.34,1.05)
\psline[arrowscale=1.5]{->}(2.2,1.7)(2.2,1.9)
\psline[arrowscale=1.5]{<-}(2.46,1.05)(2.54,1.25)
\psline[arrowscale=1.5]{->}(2.4,.65)(2.4,.85)
\psline[arrowscale=1.5]{->}(2.66,1.25)(2.74,1.05)
\psline[arrowscale=1.5]{->}(2.6,1.7)(2.6,1.9)
\psline[arrowscale=1.5]{<-}(2.86,1.05)(2.94,1.25)
\psline[arrowscale=1.5]{->}(2.8,.65)(2.8,.85)
\psline[arrowscale=1.5]{<-}(3.2,.55)(3.2,.35)
\rput[b](.9,2.5){$``i"$} \rput[b](1.85,1.8){$\ldots$}
\rput[b](2.6,2.5){$``1"$} \psframe[linecolor=black](0,.4)(3,.6)
\psframe[linecolor=black](0,2.2)(1.1,2.4)
\end{pspicture} + \sum_{j=i+1}^{n} c_j \begin{pspicture}[shift=-1.4](-.6,-.3)(3.35,2.5)
\rput[t](1,-.1){$n-1$} \rput[b](.4,1.8){$\ldots$} \psline(1,0)(1,.4)
\psline(.1,.6)(.1,2.2) \psline(.7,.6)(.7,2.2)
\psline(.9,1.8)(.9,2.2) \psline(.9,1.8)(1.1,1.6)
\psline(1.1,1.6)(1.1,1.4) \psline(1.1,1.4)(1.3,.9)
\psline(1.3,1.8)(1.1,1.6) \psline(1.3,1.8)(1.3,2.2)
\psline(1.3,.6)(1.3,.9) \psline(1.5,1.4)(1.3,.9)
\psline(1.5,1.4)(1.5,2.2) \psline(1.5,1.4)(1.7,.9)
\psline(2.2,1.4)(2,.9) \psline(2.2,1.4)(2.4,.9)
\psline(2.2,1.4)(2.2,2.2) \psline(2.4,.9)(2.6,1.4)
\psline(2.4,.6)(2.4,.9) \psline(2.6,1.4)(2.8,.9)
\psline(2.6,1.4)(2.6,2.2) \psline(2.8,.9)(3,1.4)
\psline(2.8,.9)(2.8,.6) \psline(3,1.4)(3.2,.9)
\psline(3.2,.9)(3.2,0) \psline[arrowscale=1.5]{->}(.1,1.3)(.1,1.5)
\psline[arrowscale=1.5]{->}(.7,1.3)(.7,1.5)
\psline[arrowscale=1.5]{->}(.9,1.9)(.9,2.1)
\psline[arrowscale=1.5]{->}(1.16,1.25)(1.24,1.05)
\psline[arrowscale=1.5]{->}(1.3,1.9)(1.3,2.1)
\psline[arrowscale=1.5]{->}(1.3,.65)(1.3,.85)
\psline[arrowscale=1.5]{->}(1.44,1.25)(1.36,1.05)
\psline[arrowscale=1.5]{->}(1.5,1.7)(1.5,1.9)
\psline[arrowscale=1.5]{->}(1.56,1.25)(1.64,1.05)
\psline[arrowscale=1.5]{->}(2.14,1.25)(2.06,1.05)
\psline[arrowscale=1.5]{->}(2.26,1.25)(2.34,1.05)
\psline[arrowscale=1.5]{->}(2.2,1.7)(2.2,1.9)
\psline[arrowscale=1.5]{<-}(2.46,1.05)(2.54,1.25)
\psline[arrowscale=1.5]{->}(2.4,.65)(2.4,.85)
\psline[arrowscale=1.5]{->}(2.66,1.25)(2.74,1.05)
\psline[arrowscale=1.5]{->}(2.6,1.7)(2.6,1.9)
\psline[arrowscale=1.5]{<-}(2.86,1.05)(2.94,1.25)
\psline[arrowscale=1.5]{->}(2.8,.65)(2.8,.85)
\psline[arrowscale=1.5]{<-}(3.2,.55)(3.2,.35)
\rput[b](.9,2.5){$``j"$} \rput[b](1.85,1.8){$\ldots$}
\rput[b](2.6,2.5){$``1"$} \psframe[linecolor=black](0,.4)(3,.6)
\psframe[linecolor=black](0,2.2)(1.7,2.4)
\end{pspicture}
\end{align*} Inductively we assume all $c_k=0$ where $k < i$. If we attach the
clasp of weight $(n-i+1)$ to the left top of each web in the right
side of equation~\ref{a2exp}, the $i$-th web corresponding to the
coefficient $c_i$ is nonzero because it is a cut out from the
hexagonal tiling of the plane. All other remaining webs
corresponding to $c_k$ where $k\ge i+1$ are zero because the same
reason. Therefore, $c_i=0$. This completes the proof of linear
independency.

Now, we can show that the set of the webs in the right-hand side in
equation~\ref{a2exp} is a basis by counting the dimension of web
spaces. If we set $a=(n+1), b=1$, we find the dimension of the web
space of $V_{\lambda_1}^{\otimes n+1}\otimes V_{\lambda_2}\otimes
V_{(n-1)\lambda_2}$ is $n$ by Lemma~\ref{lema2dim}. Therefore, these
webs in right side of equation~\ref{a2exp} form a basis for the
single clasp expansion.

We put
$$\begin{pspicture}[shift=-1.3](-.1,-.3)(2.6,2.3)
\rput[t](1,-.1){$n$}\qline(1,0)(1,.4) \rput[c](1,1.2){$\ldots$}
\qline(.1,.6)(.1,2) \qline(.3,.6)(.3,2) \qline(.5,.6)(.5,2)
\qline(1.5,.6)(1.5,2) \qline(1.7,.6)(1.7,2) \qline(1.9,.6)(1.9,2)
\psline[arrowscale=1.5]{->}(.1,1)(.1,1.2)
\psline[arrowscale=1.5]{->}(.3,1)(.3,1.2)
\psline[arrowscale=1.5]{->}(.5,1)(.5,1.2)
\psline[arrowscale=1.5]{->}(1.5,1)(1.5,1.2)
\psline[arrowscale=1.5]{->}(1.7,1)(1.7,1.2)
\psline[arrowscale=1.5]{->}(1.9,1)(1.9,1.2)
\psframe[linecolor=black](0,.4)(2,.6)
\end{pspicture}
=  a_1\begin{pspicture}[shift=-1.3](-.1,-.3)(2.6,2.3)
\rput[t](.9,-.1){$n-1$}\qline(1,0)(1,.4) \rput[c](1,1.2){$\ldots$}
\qline(.1,.6)(.1,2) \qline(.3,.6)(.3,2) \qline(.5,.6)(.5,2)
\qline(1.5,.6)(1.5,2) \qline(1.7,.6)(1.7,2) \qline(2.1,0)(2.1,2)
\psline[arrowscale=1.5]{->}(.1,1)(.1,1.2)
\psline[arrowscale=1.5]{->}(.3,1)(.3,1.2)
\psline[arrowscale=1.5]{->}(.5,1)(.5,1.2)
\psline[arrowscale=1.5]{->}(1.5,1)(1.5,1.2)
\psline[arrowscale=1.5]{->}(1.7,1)(1.7,1.2)
\psline[arrowscale=1.5]{->}(2.1,1)(2.1,1.2)
\psframe[linecolor=black](0,.4)(1.8,.6)
\end{pspicture} + \sum_{i=2}^{n} a_i \begin{pspicture}[shift=-1.3](-.6,-.3)(3.35,2.5)
\rput[t](1,-.1){$n-1$} \rput[b](.4,1.8){$\ldots$} \psline(1,0)(1,.4)
\psline(.1,.6)(.1,2.2) \psline(.7,.6)(.7,2.2)
\psline(.9,1.8)(.9,2.2) \psline(.9,1.8)(1.1,1.6)
\psline(1.1,1.6)(1.1,1.4) \psline(1.1,1.4)(1.3,.9)
\psline(1.3,1.8)(1.1,1.6) \psline(1.3,1.8)(1.3,2.2)
\psline(1.3,.6)(1.3,.9) \psline(1.5,1.4)(1.3,.9)
\psline(1.5,1.4)(1.5,2.2) \psline(1.5,1.4)(1.7,.9)
\psline(2.2,1.4)(2,.9) \psline(2.2,1.4)(2.4,.9)
\psline(2.2,1.4)(2.2,2.2) \psline(2.4,.9)(2.6,1.4)
\psline(2.4,.6)(2.4,.9) \psline(2.6,1.4)(2.8,.9)
\psline(2.6,1.4)(2.6,2.2) \psline(2.8,.9)(3,1.4)
\psline(2.8,.9)(2.8,.6) \psline(3,1.4)(3.2,.9)
\psline(3.2,.9)(3.2,0) \psline[arrowscale=1.5]{->}(.1,1.3)(.1,1.5)
\psline[arrowscale=1.5]{->}(.7,1.3)(.7,1.5)
\psline[arrowscale=1.5]{->}(.9,1.9)(.9,2.1)
\psline[arrowscale=1.5]{->}(1.16,1.25)(1.24,1.05)
\psline[arrowscale=1.5]{->}(1.3,1.9)(1.3,2.1)
\psline[arrowscale=1.5]{->}(1.3,.65)(1.3,.85)
\psline[arrowscale=1.5]{->}(1.44,1.25)(1.36,1.05)
\psline[arrowscale=1.5]{->}(1.5,1.7)(1.5,1.9)
\psline[arrowscale=1.5]{->}(1.56,1.25)(1.64,1.05)
\psline[arrowscale=1.5]{->}(2.14,1.25)(2.06,1.05)
\psline[arrowscale=1.5]{->}(2.26,1.25)(2.34,1.05)
\psline[arrowscale=1.5]{->}(2.2,1.7)(2.2,1.9)
\psline[arrowscale=1.5]{<-}(2.46,1.05)(2.54,1.25)
\psline[arrowscale=1.5]{->}(2.4,.65)(2.4,.85)
\psline[arrowscale=1.5]{->}(2.66,1.25)(2.74,1.05)
\psline[arrowscale=1.5]{->}(2.6,1.7)(2.6,1.9)
\psline[arrowscale=1.5]{<-}(2.86,1.05)(2.94,1.25)
\psline[arrowscale=1.5]{->}(2.8,.65)(2.8,.85)
\psline[arrowscale=1.5]{<-}(3.2,.55)(3.2,.35)
\rput[b](.9,2.4){$``i"$} \rput[b](1.85,1.8){$\ldots$}
\rput[b](2.6,2.4){$``1"$} \psframe[linecolor=black](0,.4)(3,.6)
\end{pspicture},$$ with
some $a_i$, since the webs in the right hand side span the web space
which contains the web of the left hand side. If we attach a
$\begin{pspicture}[shift=.2](-.17,-.17)(.17,.17)\qline(0,0)(0,.17)
\qline(-.15,-.15)(0,0)\qline(.15,-.15)(0,0)\end{pspicture}$ on the top
of webs in equation~\ref{a2exp}, the left side of
equation~\ref{a2exp} becomes zero and all webs in the right-hand side of
equation~\ref{a2exp} become multiples of a web. Thus we get the
following $n-1$ equations.

$$a_{n-1}-[2]a_{n}=0.$$
For $i=1, 2, \ldots, n-2$,

$$a_i-[2]a_{i+1}+a_{i+2}=0.$$
From these equations, we are able to find the relations between the
coefficients $a_i$'s. By a normalization, attaching the clasp of
weight $(n,0)$ to the top of each web in the equation, we find
$a_1=1$. Then other coefficients can be found subsequently.
\end{proof}

\subsection{Single clasp expansions of a non-segregated clasp of
weight $(a,b)$}\label{sl3ab}

The most interesting case is a single clasp expansion of the clasp
of weight $(a,b)$ where $a\neq 0 \neq b$. By Lemma~\ref{lema2dim},
we know the number of webs in a single clasp expansion of the clasp
of weight $(a,b)$ is $(a+1)b$. We need a set of basis webs with a
nice rectangular order, but we can not find one in the general case.
Even if one finds such a basis, each web in the basis would have
many hexagonal faces which make it very difficult to get numerical
relations. So we start from an alternative, \emph{non-segregated}
clasp. A non-segregated clasp is obtained from the segregated clasp
by attaching a sequence of $H$'s until we get the desired shape of
edge orientations. Fortunately, there is a canonical way to find
them by putting $H$'s from the leftmost string of weight $\lambda_2$
or $-$ until it reach to the desired position. The left side of the
equation in Figure \ref{nsclaspex} is an example of a non-segregated
clasp of weight $(2,3)$. the right-hand side of the equation in Figure
\ref{nsclaspex} shows a sequence of $H$'s which illustrates how we
obtain it from the segregated clasp of weight $(2,3)$.

First of all, we can show that the non-segregated clasps are
well-defined \cite[Lemma 2.6]{kim:thesis}. One can prove that
non-segregated clasps also satisfy two properties of segregated
clasps: 1) two consecutive non-segregated clasps is equal to a
non-segregated clasp, 2) if we attach a web to a non-segregated
clasp and if it has a cut path whose weight is less than the weight
of the clasp, then it is zero \cite[Lemma 2.7]{kim:thesis}. We find
a single clasp expansion of a non-segregated clasp of weight $(a,b)$
in shown in equation~\ref{a2abexp}. Kuperberg showed that for a
fixed boundary, the interior can be filled by a cut out from the
hexagonal tiling of the plane with the given boundary
\cite{Kuperberg:spiders}. For our cases, there are two possible
fillings but we use the maximal cut out of the hexagonal tiling. We
draw examples of the case $i=6, j=5$ and the first one in
Figure~\ref{hexaex} is not a maximal cut out and the second one is
the maximal cut out which fits to the left rectangle and the last
one is the maximal cut out which fits to the right rectangle as the
number indicated in equation~\ref{a2abexp}. An example of a single
clasp expansion of a segregated clasp of weight $(2,2)$ can be found
in \cite[pp. 18]{kim:thesis}.

\begin{thm} For $a, b\ge 1$,
\begin{align}
\begin{pspicture}[shift=-1.5](-.3,-.7)(2.5,2.5) \psline(.1,.1)(.1,2)
\psline(.1,-.1)(.1,-.6) \psline(.3,.1)(.3,2) \psline(.3,-.1)(.3,-.6)
\psline(.5,.1)(.5,2) \psline(.5,-.1)(.5,-.6) \psline(1.5,.1)(1.5,2)
\psline(1.5,-.1)(1.5,-.6) \psline(1.7,.1)(1.7,2)
\psline(1.7,-.1)(1.7,-.6) \psline(1.9,.1)(1.9,2)
\psline(1.9,-.1)(1.9,-.5) \psline(1.9,-.5)(2.1,-.5)
\psline(2.1,-.5)(2.1,2) \uput{3pt}[90](1,.7){$\ldots$}
\uput{3pt}[90](.5,2.1){$\overset{b}{\overbrace{~~~~~~~}}$}
\uput{3pt}[90](1.7,2.1){$\overset{a+1}{\overbrace{~~~~~~~~~}}$}
\uput{3pt}[270](1,-.3){$\ldots$}
\psline[arrowscale=1.5]{<-}(.1,.8)(.1,1)
\psline[arrowscale=1.5]{<-}(.3,.8)(.3,1)
\psline[arrowscale=1.5]{<-}(.5,.8)(.5,2)
\psline[arrowscale=1.5]{->}(1.5,.8)(1.5,1)
\psline[arrowscale=1.5]{->}(1.7,.8)(1.7,1)
\psline[arrowscale=1.5]{->}(1.9,.8)(1.9,1)
\psline[arrowscale=1.5]{->}(2.1,.8)(2.1,1)
\psframe[linecolor=black](0,-.1)(2,.1)
\end{pspicture}
= \sum_{i=1}^{b}\sum_{j=0}^{a} \hskip .3cm
\frac{[b-i+1][b+j+1]}{[b][a+b+1]} \begin{pspicture}[shift=-1.5](-.4,-.8)(7,2.5)
\psline(.1,.1)(.1,2) \psline(.1,-.1)(.1,-.6) \psline(.3,.1)(.3,2)
\psline(.3,-.1)(.3,-.6) \psline(1.3,.1)(1.3,2)
\psline(1.3,-.1)(1.3,-.6) \psline(1.5,1)(1.5,2)
\psline(1.5,1)(1.8,1) \psline(2.2,2)(2.2,1.3)
\psline(2.2,.7)(2.2,.1) \psline(2.2,-.1)(2.2,-.6)
\psline(3,1)(3.8,1) \psline(4.6,1.3)(4.6,2)
\psline(4.6,-.1)(4.6,-.6) \psline(4.6,.7)(4.6,.1)
\psline(5,1)(5.3,1) \psline(5.3,1)(5.3,2) \psline(5.5,-.1)(5.5,-.6)
\psline(5.5,.1)(5.5,2) \psline(6.5,-.1)(6.5,-.6)
\psline(6.5,.1)(6.5,2) \psline(6.7,-.1)(6.7,-.6)
\psline(6.7,.1)(6.7,2) \uput{3pt}[90](1.5,2.1){$1$}
\uput{3pt}[90](2.2,2.1){$i-1$} \uput{3pt}[90](4.6,2.1){$j$}
\uput{3pt}[90](5.3,2.1){$1$} \uput{3pt}[0](2.2,0.3){$i-1$}
\uput{3pt}[0](4.6,0.3){$j$} \uput{3pt}[270](.8,1.3){$\ldots$}
\uput{3pt}[270](6,1.3){$\ldots$} \uput{3pt}[270](.8,-.3){$\ldots$}
\uput{3pt}[270](6,-.3){$\ldots$}
\psline[arrowscale=1.5]{<-}(.1,.8)(.1,1)
\psline[arrowscale=1.5]{<-}(.3,.8)(.3,1)
\psline[arrowscale=1.5]{<-}(1.3,.8)(1.3,1)
\psline[arrowscale=1.5]{->}(1.5,1)(1.8,1)
\psline[arrowscale=1.5]{->}(2.2,1.8)(2.2,1.6)
\psline[arrowscale=1.5]{->}(2.2,.5)(2.2,.3)
\psline[arrowscale=1.5]{->}(3.1,1)(3.3,1)
\psline[arrowscale=1.5]{->}(3.5,1)(3.7,1)
\psline[arrowscale=1.5]{->}(4.6,1.6)(4.6,1.8)
\psline[arrowscale=1.5]{<-}(4.6,.5)(4.6,.3)
\psline[arrowscale=1.5]{->}(5,1)(5.3,1)
\psline[arrowscale=1.5]{->}(5.5,.8)(5.5,1)
\psline[arrowscale=1.5]{->}(6.5,.8)(6.5,1)
\psline[arrowscale=1.5]{->}(6.7,.8)(6.7,1)
\psframe[linecolor=black](0,-.1)(6.8,.1)
\psframe[linecolor=darkgray,linewidth=1.5pt](1.8,.7)(3,1.3)
\psframe[linecolor=darkgray,linewidth=1.5pt](3.8,.7)(5,1.3)
\psline[linecolor=black,linestyle=dashed](3.4,2.2)(3.4,-.7)
\uput{1pt}[90](2.4,.7){$(1)$} \uput{1pt}[90](4.4,.7){$(2)$}
\end{pspicture}
 \label{a2abexp}
\end{align}
\label{a2abexpthm1}
\end{thm}
\begin{proof}
Let us denote the web corresponding to the coefficient
$\dfrac{[b-i+1][b+j+1]}{[b][a+b+1]}=a_{i,j}$ by $D_{i,j}$. First of
all, all these webs in the equation~\ref{a2abexp} are nonzero
because they are cut outs from the hexagonal tiling of the plane.
These webs in the right hand side of the equation form a basis
because their cardinality is the same as the dimension of the
invariant space of $V_{a\lambda_1+(b-1)\lambda_2}\otimes
V_{\lambda_1}^{\otimes a+1}\otimes V_{\lambda_2}^{\otimes b-1}$ and
they are linearly independent. Suppose that a linear combination of
webs in the right-hand side of the equation~\ref{a2abexp} is zero for
some choice of $a_{i,j}$. By attaching the clasp of weight
$(0,b-i+1)$ to the left top and $(a+1-j,0)$ on right top of webs one
can see all webs but the webs $D_{s,t}, 1\le s\le i, 0 \le t\le j$
vanish. It is clear that $a_{1,0}=0$ by attaching the clasp of
weight $(0,b)$ and $(0,a+1)$. Inductively, we can show $a_{i,j}=0$
for all $i, j$.

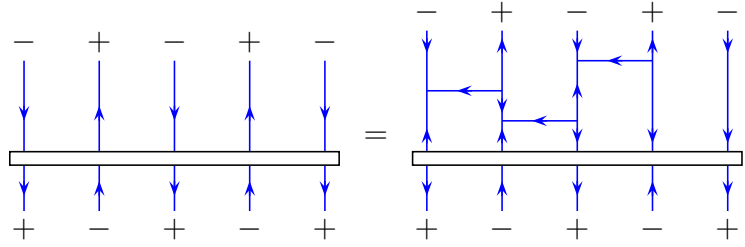
\begin{figure}
$$
\begin{pspicture}[shift=-1.3](-1.4,-1.4)(3.4,1.4) \uput{3pt}[90](0,1){$+$}
\uput{3pt}[90](-1,1){$-$} \uput{3pt}[270](-1,-1){$+$}
\uput{3pt}[90](1,1){$-$} \uput{3pt}[90](2,1){$+$}
\uput{3pt}[90](3,1){$-$} \uput{3pt}[270](0,-1){$-$}
\uput{3pt}[270](1,-1){$+$} \uput{3pt}[270](2,-1){$-$}
\uput{3pt}[270](3,-1){$+$}
\psline(0,-.2)(0,1)\psline[arrowscale=1.5]{->}(0,.2)(0,.4)
\psline(1,-.2)(1,1)\psline[arrowscale=1.5]{<-}(1,.2)(1,.4)
\psline(2,-.2)(2,1)\psline[arrowscale=1.5]{->}(2,.2)(2,.4)
\psline(3,-.2)(3,1)\psline[arrowscale=1.5]{<-}(3,.2)(3,.4)
\psline(-1,-.2)(-1,1)\psline[arrowscale=1.5]{<-}(-1,.2)(-1,.4)
\psline(-1,-.4)(-1,-1)\psline[arrowscale=1.5]{->}(-1,-.6)(-1,-.8)
\psline(0,-.4)(0,-1)\psline[arrowscale=1.5]{<-}(0,-.6)(0,-.8)
\psline(1,-.4)(1,-1)\psline[arrowscale=1.5]{->}(1,-.6)(1,-.8)
\psline(2,-.4)(2,-1)\psline[arrowscale=1.5]{<-}(2,-.6)(2,-.8)
\psline(3,-.4)(3,-1)\psline[arrowscale=1.5]{->}(3,-.6)(3,-.8)
\psframe[linecolor=black](-1.2,-.4)(3.2,-.2)
\end{pspicture}
= \begin{pspicture}[shift=-1.3](-1.4,-1.4)(3.4,1.8) \uput{3pt}[90](0,1.4){$+$}
\uput{3pt}[90](-1,1.4){$-$} \uput{3pt}[270](-1,-1){$+$}
\uput{3pt}[90](1,1.4){$-$} \uput{3pt}[90](2,1.4){$+$}
\uput{3pt}[90](3,1.4){$-$} \uput{3pt}[270](0,-1){$-$}
\uput{3pt}[270](1,-1){$+$} \uput{3pt}[270](2,-1){$-$}
\uput{3pt}[270](3,-1){$+$}
\psline(0,-.2)(0,1.4)\psline[arrowscale=1.5]{->}(0,1.1)(0,1.3)
\psline(1,-.2)(1,1.4)\psline[arrowscale=1.5]{<-}(1,1.1)(1,1.3)
\psline(2,-.2)(2,1.4)\psline[arrowscale=1.5]{->}(2,1.1)(2,1.3)
\psline(3,-.2)(3,1.4)\psline[arrowscale=1.5]{<-}(3,1.1)(3,1.3)
\psline(-1,-.2)(-1,1.4)\psline[arrowscale=1.5]{<-}(-1,1.1)(-1,1.3)
\psline(-1,-.4)(-1,-1)\psline[arrowscale=1.5]{->}(-1,-.6)(-1,-.8)
\psline(0,-.4)(0,-1)\psline[arrowscale=1.5]{<-}(0,-.6)(0,-.8)
\psline(1,-.4)(1,-1)\psline[arrowscale=1.5]{->}(1,-.6)(1,-.8)
\psline(2,-.4)(2,-1)\psline[arrowscale=1.5]{<-}(2,-.6)(2,-.8)
\psline(3,-.4)(3,-1)\psline[arrowscale=1.5]{->}(3,-.6)(3,-.8)
\psline(-1,.6)(0,.6)\psline[arrowscale=1.5]{->}(-.4,.6)(-.6,.6)
\psline(1,.2)(0,.2)\psline[arrowscale=1.5]{->}(.6,.2)(.4,.2)
\psline(1,1)(2,1)\psline[arrowscale=1.5]{->}(1.6,1)(1.4,1)
\psline[arrowscale=1.5]{->}(-1,-.1)(-1,.1)\psline[arrowscale=1.5]{->}(0,-.1)(0,.1)
\psline[arrowscale=1.5]{<-}(1,-.1)(1,.1)\psline[arrowscale=1.5]{<-}(2,-.1)(2,.1)
\psline[arrowscale=1.5]{<-}(3,-.1)(3,.1)\psline[arrowscale=1.5]{<-}(0,.3)(0,.5)
\psline[arrowscale=1.5]{->}(1,.5)(1,.7)
\psframe[linecolor=black](-1.2,-.4)(3.2,-.2)
\end{pspicture}
$$
\caption{A non-segregated clasp of weight $(2,3)$.}\label{nsclaspex}
\end{figure}

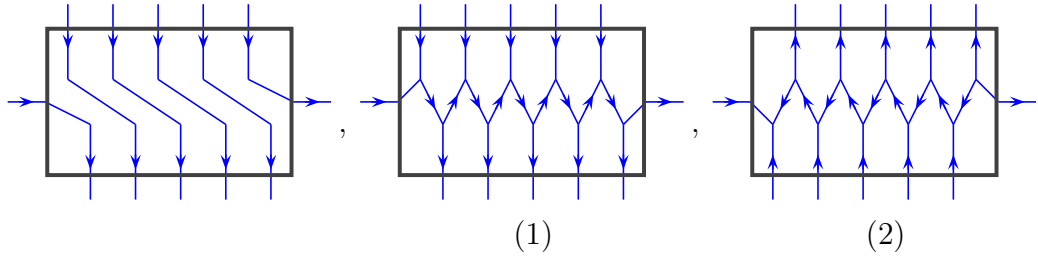
\begin{figure}
$$
\begin{pspicture}[shift=-1.7](-.6,-2.1)(3.9,1.5)
\psline(-.5,0)(0,0)\psline[arrowscale=1.5]{->}(-.35,0)(-.15,0)
\psline(3.3,0)(3.8,0)\psline[arrowscale=1.5]{->}(3.45,0)(3.65,0)
\multips(.3,.3)(.6,0){5}{\psline(0,0)(0,1)\psline[arrowscale=1.5]{->}(0,.6)(0,.4)}
\multips(.6,-1.3)(.6,0){5}{\psline(0,0)(0,1)\psline[arrowscale=1.5]{->}(0,.6)(0,.4)}
\multips(.3,.3)(.6,0){4}{\psline(0,0)(0.9,-.6)}
\psline(0,0)(.6,-.3)\psline(3.3,0)(2.7,.3)
\psframe[linecolor=darkgray,linewidth=1.5pt](0,-1)(3.3,1)
\end{pspicture}
, \begin{pspicture}[shift=-1.7](-.6,-2.1)(3.9,1.5)
\psline(-.5,0)(0,0)\psline[arrowscale=1.5]{->}(-.35,0)(-.15,0)
\psline(3.3,0)(3.8,0)\psline[arrowscale=1.5]{->}(3.45,0)(3.65,0)
\multips(.3,.3)(.6,0){5}{\psline(0,0)(0,1)\psline[arrowscale=1.5]{->}(0,.6)(0,.4)}
\multips(.6,-1.3)(.6,0){5}{\psline(0,0)(0,1)\psline[arrowscale=1.5]{->}(0,.6)(0,.4)}
\multips(.3,.3)(.6,0){5}{\psline(0,0)(0.3,-.6)\psline[arrowscale=1.5]{->}(0.1,-.2)(0.2,-.4)}
\multips(.6,-.3)(.6,0){4}{\psline(0,0)(0.3,.6)\psline[arrowscale=1.5]{->}(0.1,.2)(0.2,.4)}
\psline(0,0)(.3,.3)\psline(3.3,0)(3,-.3)
\psframe[linecolor=darkgray,linewidth=1.5pt](0,-1)(3.3,1)
\uput{3pt}[90](1.8,-2.1){$(1)$}
\end{pspicture}
, \begin{pspicture}[shift=-1.7](-.6,-2.1)(3.9,1.5)
\psline(-.5,0)(0,0)\psline[arrowscale=1.5]{->}(-.35,0)(-.15,0)
\psline(3.3,0)(3.8,0)\psline[arrowscale=1.5]{->}(3.45,0)(3.65,0)
\multips(.6,.3)(.6,0){5}{\psline(0,0)(0,1)\psline[arrowscale=1.5]{<-}(0,.6)(0,.4)}
\multips(.3,-1.3)(.6,0){5}{\psline(0,0)(0,1)\psline[arrowscale=1.5]{<-}(0,.6)(0,.4)}
\multips(.6,.3)(.6,0){4}{\psline(0,0)(0.3,-.6)\psline[arrowscale=1.5]{<-}(0.1,-.2)(0.2,-.4)}
\multips(.3,-.3)(.6,0){5}{\psline(0,0)(0.3,.6)\psline[arrowscale=1.5]{<-}(0.1,.2)(0.2,.4)}
\psline(0,0)(.3,-.3)\psline(3.3,0)(3,.3)
\psframe[linecolor=darkgray,linewidth=1.5pt](0,-1)(3.3,1)
\uput{3pt}[90](1.8,-2.1){$(2)$}
\end{pspicture}
$$
\caption{Fillings for the boxes in equation \ref{a2abexp}.}
\label{hexaex}
\end{figure}

To find $a_{i,j}$, we attach a
$\begin{pspicture}[shift=-.12](-.17,-.17)(.17,.17)\qline(0,0)(0,.17)\qline(-.15,-.15)(0,0)\qline(.15,-.15)(0,0)\end{pspicture}$
or a
$\begin{pspicture}[shift=-.12](-.17,-.17)(.17,.17)\qline(-.15,-.15)(-.15,.05)\qline(.15,-.15)(.15,.05)
\pccurve[angleA=90,angleB=90,ncurv=1](-.15,.05)(.15,.05)\end{pspicture}$
to find one exceptional and three types of equations as follow.

$$[3]a_{1,0}-[2]a_{1,1}-[2]a_{2,0}+a_{2,1}=0.$$
Type I : For $j=0,1,\ldots, a$,

$$a_{b-1,j}-[2]a_{b,j}=0.$$
Type II : For $i=1,2, \ldots, b-2$ and $j=0,1,\ldots, a$.

$$a_{i,j}-[2]a_{i+1,j}+a_{i+2,j}=0.$$
Type III : For $i=1,2, \ldots, b$ and $j=0,1,\ldots, a-2$.

$$a_{i,j}-[2]a_{i,j+1}+a_{i,j+2}=0.$$
If we set $a_{1,0}=x$, then inductively one can see that the
coefficient $a_{i,j}$ in the equation~\ref{a2abexp} is
$$\frac{[b-i+1][b+j+1]}{[b][b+1]}x.$$

One might check that these are the right coefficients. Usually we
normalize one basis web in the expansion to get a known value. But
we can not normalize for this expansion yet because it is not a
segregated clasp. Thus we use a complicate procedure in
Lemma~\ref{a2abexplem3} to find that the coefficient of $a_{1,a}$ is
$1$. Then, we find that $a_{1,0}$ is $\dfrac{[b+1]}{[a+b+1]}$ and it
completes the proof of the theorem.
\end{proof}

We find a double clasps expansion as shown in
Theorem~\ref{a2ab2expthm}, the box between two clasps is filled by
the unique maximal cut out from the hexagonal tiling with the given
boundary as we have seen in Figure~\ref{hexaex}.

\begin{thm}
For $a, b \ge 1 $,
$$
\begin{pspicture}[shift=-1.7](-.5,-1)(1,2.5) \rput[t](0,-.5){$a$}
\rput[b](0,1.9){$a$} \rput[b](.7,1.9){$b$} \rput[t](.7,-.5){$b$}
\psline(0,.6)(0,1.8)\psline[arrowscale=1.5]{->}(0,.9)(0,1.1)
\psline(0,-.4)(0,.4)\psline[arrowscale=1.5]{<-}(0,.1)(0,-.1)
\psline(.7,1.8)(.7,.6)\psline[arrowscale=1.5]{<-}(.7,.9)(.7,1.1)
\psline(.7,.4)(.7,-.4)\psline[arrowscale=1.5]{->}(.7,.1)(.7,-.1)
\psframe[linecolor=black](-.2,.4)(.9,.6)
\end{pspicture}
= \begin{pspicture}[shift=-1.7](-.5,-1)(2.2,2.5) \rput[t](0,-.5){$a$}
\rput[b](0,1.9){$a$} \rput[b](1.2,1.9){$b-1$}
\rput[t](1.2,-.5){$b-1$}
\psline(0,.6)(0,1.8)\psline[arrowscale=1.5]{->}(0,.9)(0,1.1)
\psline(0,-.4)(0,.4)\psline[arrowscale=1.5]{<-}(0,.1)(0,-.1)
\psline(1.2,1.8)(1.2,.6)\psline[arrowscale=1.5]{<-}(1.2,.9)(1.2,1.1)
\psline(1.2,.4)(1.2,-.4)\psline[arrowscale=1.5]{->}(1.2,.1)(1.2,-.1)
\psline(2,1.8)(2,-.4)\psline[arrowscale=1.5]{->}(2,.1)(2,-.1)
\psframe[linecolor=black](-.2,.4)(1.7,.6)
\end{pspicture}
+\frac{[b-1]}{[b]} \begin{pspicture}[shift=-1.7](-.8,-1)(2.5,2.5)
\rput[b](0,1.9){$a$}
\rput[bl](.4,1.9){$b-1$}\rput[t](0,-.5){$a$}\rput[tl](.4,-.5){$b-1$}
\rput[bl](0.5,.4){$b-2$}\rput[br](-.1,.4){$a$}
\psline(0,-.4)(0,0)\psline[arrowscale=1.5]{->}(0,-.3)(0,-.1)
\psline(0,.2)(0,1)\psline[arrowscale=1.5]{->}(0,.5)(0,.7)
\psline(0,1.2)(0,1.8)\psline[arrowscale=1.5]{->}(0,1.4)(0,1.6)
\psline(.4,0)(.4,-0.4)\psline[arrowscale=1.5]{->}(.4,-.1)(.4,-.3)
\psline(.4,1)(.4,.2)\psline[arrowscale=1.5]{->}(.4,.6)(.4,.4)
\psline(.4,1.8)(.4,1.2)\psline[arrowscale=1.5]{->}(.4,1.6)(.4,1.4)
\psline(2.25,1.8)(2.25,1)\psline[arrowscale=1.5]{->}(2.25,1.6)(2.25,1.4)
\psline(2.25,.2)(2.25,-.4)\psline[arrowscale=1.5]{->}(2.25,0)(2.25,-.2)
\psline(2,.45)(2,.75)\psline[arrowscale=1.5]{->}(2,.5)(2,.7)
\pccurve[angleA=270,angleB=180,ncurv=1](1.75,1)(2,.75)
\pccurve[angleA=0,angleB=270,ncurv=1](2,.75)(2.25,1)
\pccurve[angleA=90,angleB=180,ncurv=1](1.75,.2)(2,.45)
\pccurve[angleA=0,angleB=90,ncurv=1](2,.45)(2.25,.2)
\psframe[linecolor=black](-.2,0)(1.9,.2)
\psframe[linecolor=black](-.2,1)(1.9,1.2)\end{pspicture}
-\frac{[a]}{[b][a+b+1]} \begin{pspicture}[shift=-1.6](0,1.7)(0,1.7) \begin{pspicture}[shift=-1.7](-1.3,-1)(2.7,2.5)
\rput[br](0,1.9){$a-1$} \rput[br](1.75,1.9){$b-1$}
\rput[tr](0,-.6){$a-1$} \rput[tr](1.75,-.6){$b-1$}
\rput[b](0.5,1.9){$1$} \rput[b](2.25,1.9){$1$}
\rput[t](0.5,-.6){$1$} \rput[t](2.25,-.6){$1$}
\rput[br](-.1,.4){$a-1$}
\psline(0,-.5)(0,-.1)\psline[arrowscale=1.5]{->}(0,-.4)(0,-.4)
\psline(0,.1)(0,1.1)\psline[arrowscale=1.5]{->}(0,.4)(0,.6)
\psline(0,1.3)(0,1.8)\psline[arrowscale=1.5]{->}(0,1.45)(0,1.65)
\psline(.5,-.5)(.5,-.1)\psline[arrowscale=1.5]{->}(.5,-.4)(.5,-.2)
\psline(.5,1.3)(.5,1.8)\psline[arrowscale=1.5]{->}(.5,1.45)(.5,1.65)
\psline(1.5,-.1)(1.5,-0.5)\psline[arrowscale=1.5]{->}(1.5,-.2)(1.5,-.4)
\psline(1.5,1.8)(1.5,1.3)\psline[arrowscale=1.5]{->}(1.5,1.65)(1.5,1.45)
\psline(2.25,1.8)(2.25,1.1) \psline(2.25,.1)(2.25,-.5)
\pcline(.65,.75)(.5,1.1)\middlearrow
\pcline(.5,.1)(.65,.45)\middlearrow
\pcline(1.5,1.1)(1.25,.75)\middlearrow
\pcline(1.25,.45)(1.5,.1)\middlearrow
\pcline(2.25,1.1)(1.75,.65)\middlearrow
\pcline(1.75,.55)(2.25,.1)\middlearrow
\psline[linecolor=darkgray,linewidth=1.5pt](.5,.6)(.75,.85)
\psline[linecolor=darkgray,linewidth=1.5pt](.5,.6)(.75,.35)
\psline[linecolor=darkgray,linewidth=1.5pt](2,.6)(.75,.85)
\psline[linecolor=darkgray,linewidth=1.5pt](2,.6)(.75,.35)
\psframe[linecolor=black](-.2,-.1)(1.9,.1)
\psframe[linecolor=black](-.2,1.1)(1.9,1.3)\end{pspicture}\end{pspicture}
$$
\label{a2ab2expthm}
\end{thm}
\begin{proof}
It follows from Lemma \ref{a2abexplem2} and Lemma~\ref{a2abexplem3}.
\end{proof}

The expansion in equation depicted in Proposition~\ref{a2abexpprop1}
was first used to define the segregated clasp of weight
$(a,b)$~\cite{OY:quantum}. The clasps can be constructed from web
spaces \cite{Kuperberg:spiders} and these two are known to be equal.
We will apply Theorem~\ref{a2abexpthm1} to demonstrate the
effectiveness of single clasp expansions by deriving the
coefficients in Proposition~\ref{a2abexpprop1}.

\begin{prop} [\cite{OY:quantum}] A quadruple clasps expansion of the segregated clasp of
weight $(a,b)$ is
$$
\begin{pspicture}[shift=-1.4](-.5,-1)(1.3,2)
\rput[t](0,-.4){$a$}\psline(0,-.3)(0,.4)\psline[arrowscale=1.5]{->}(0,0)(0,.2)
\rput[t](.7,-.4){$b$}\psline(.7,.4)(.7,-.3)\psline[arrowscale=1.5]{->}(.7,.2)(.7,0)
\psline(0,.6)(0,1.4)\psline[arrowscale=1.5]{->}(0,.9)(0,1.1)
\psline(.7,1.4)(.7,.6)\psline[arrowscale=1.5]{->}(.7,1.1)(.7,.9)
\rput[b](0,1.5){$a$} \rput[b](.7,1.5){$b$}
\psframe[linecolor=black](-.2,.4)(.9,.6)
\end{pspicture}
= \sum_{k=0}^{\mathrm{Min}(a,b)}
(-1)^k\frac{[a]![b]![a+b-k+1]!}{[a-k]![b-k]![k]![a+b+1]!}
\begin{pspicture}[shift=-1.9](0,2)(0,2)\begin{pspicture}[shift=-1.9](-2.9,-2)(2.7,2)
\psline(-1.3,-.9)(-1.3,.9)\psline[arrowscale=1.5]{->}(-1.3,-.1)(-1.3,.1)
\rput[br](-1.5,0){$a-k$}
\psline(1.3,.9)(1.3,-.9)\psline[arrowscale=1.5]{->}(1.3,.1)(1.3,-.1)
\rput[bl](1.5,0){$b-k$}
\pccurve[angleA=270,angleB=270,ncurv=1](.7,.9)(-.7,.9)\middlearrow
\pccurve[angleA=90,angleB=90,ncurv=1](-.7,-.9)(.7,-.9)\middlearrow
\rput(0,0){$k$}
\psline(-1,-1.5)(-1,-1.1)\psline[arrowscale=1.5]{->}(-1,-1.4)(-1,-1.2)
\psline(-1,1.1)(-1,1.5)\psline[arrowscale=1.5]{->}(-1,1.2)(-1,1.4)
\psline(1,-1.1)(1,-1.5)\psline[arrowscale=1.5]{->}(1,-1.2)(1,-1.4)
\psline(1,1.5)(1,1.1)\psline[arrowscale=1.5]{->}(1,1.4)(1,1.2)
\rput[b](-1,1.7){$a$}\rput[b](1,1.7){$b$}
\rput[t](-1,-1.7){$a$}\rput[t](1,-1.7){$b$}
\psframe[linecolor=black](-1.5,-1.1)(-.5,-.9)
\psframe[linecolor=black]( 1.5,-1.1)( .5,-.9)
\psframe[linecolor=black](-1.5, 1.1)(-.5, .9)
\psframe[linecolor=black]( 1.5, 1.1)( .5, .9)
\end{pspicture} \end{pspicture}.
$$
\label{a2abexpprop1}
\end{prop}

\begin{proof}

\begin{figure}
\begin{eqnarray}&
\begin{pspicture}[shift=-1.4](-.5,-1)(1.3,2)
\rput[t](0,-.4){$a$}\psline(0,-.3)(0,.4)\psline[arrowscale=1.5]{->}(0,0)(0,.2)
\rput[t](.7,-.4){$b$}\psline(.7,.4)(.7,-.3)\psline[arrowscale=1.5]{->}(.7,.2)(.7,0)
\psline(0,.6)(0,1.4)\psline[arrowscale=1.5]{->}(0,.9)(0,1.1)
\psline(.7,1.4)(.7,.6)\psline[arrowscale=1.5]{->}(.7,1.1)(.7,.9)
\rput[b](0,1.5){$a$} \rput[b](.7,1.5){$b$}
\psframe[linecolor=black](-.2,.4)(.9,.6)
\end{pspicture}  =
\begin{pspicture}[shift=-1.9](-1.9,-2)(2.2,2)
\psline(-1,-.3)(-1,.9)\psline[arrowscale=1.5]{->}(-1,.2)(-1,.4)
\psline(-1,-.9)(-1,-.5)\psline[arrowscale=1.5]{->}(-1,-.8)(-1,-.6)
\psline(1,.9)(1,-.3)\psline[arrowscale=1.5]{->}(1,.4)(1,.2)
\psline(1,-.5)(1,-.9)\psline[arrowscale=1.5]{->}(1,-.6)(1,-.8)
\psline(1.4,.9)(1.4,-.9)\psline[arrowscale=1.5]{->}(1.4,.1)(1.4,-.1)
\psline(-1,-1.5)(-1,-1.1)\psline[arrowscale=1.5]{->}(-1,-1.4)(-1,-1.2)
\psline(-1,1.1)(-1,1.5)\psline[arrowscale=1.5]{->}(-1,1.2)(-1,1.4)
\psline(1,-1.1)(1,-1.5)\psline[arrowscale=1.5]{->}(1,-1.2)(1,-1.4)
\psline(1,1.5)(1,1.1)\psline[arrowscale=1.5]{->}(1,1.4)(1,1.2)
\rput[b](-1,1.7){$a$}\rput[b](1,1.7){$b$}
\rput[br](-1.2,.2){$a$}\rput[br](.8,.2){$b-1$}
\rput[t](-1,-1.7){$a$}\rput[t](1,-1.7){$b$}
\psframe[linecolor=black](-1.5,-1.1)(-.5,-.9)
\psframe[linecolor=black]( 1.5,-1.1)( .5,-.9)
\psframe[linecolor=black](-1.5, 1.1)(-.5, .9)
\psframe[linecolor=black]( 1.5, 1.1)( .5, .9)
\psframe[linecolor=black](-1.3,-.5)(1.3, -.3)
\end{pspicture}
-\frac{[a]}{[a+b+1]} \begin{pspicture}[shift=-1.9](0,2)(0,2)\begin{pspicture}[shift=-1.9](-2.4,-2)(2.2,2)
\psline(-1,-.3)(-1,.9)\psline[arrowscale=1.5]{->}(-1,.2)(-1,.4)
\rput[br](-1.2,0){$a-1$}
\psline(-1,-.9)(-1,-.5)\psline[arrowscale=1.5]{->}(-1,-.8)(-1,-.6)
\psline(1,.9)(1,.2)\psline[arrowscale=1.5]{->}(1,.7)(1,.5)
\rput[bl](1.2,.4){$b-1$}
\psline(1,-.5)(1,-.9)\psline[arrowscale=1.5]{->}(1,-.6)(1,-.8)
\pccurve[angleA=90,angleB=90,ncurv=1](1.2,.2)(1.4,.2)
\pccurve[angleA=225,angleB=315,ncurv=1](.7,.9)(-.7,.9)\middlearrow
\psline(1.4,.2)(1.4,-.9)\psline[arrowscale=1.5]{->}(1.4,-.3)(1.4,-.5)
\psline(-1,-1.5)(-1,-1.1)\psline[arrowscale=1.5]{->}(-1,-1.4)(-1,-1.2)
\psline(-1,1.1)(-1,1.5)\psline[arrowscale=1.5]{->}(-1,1.2)(-1,1.4)
\psline(1,-1.1)(1,-1.5)\psline[arrowscale=1.5]{->}(1,-1.2)(1,-1.4)
\psline(1,1.5)(1,1.1)\psline[arrowscale=1.5]{->}(1,1.4)(1,1.2)
\psline(.3,-.3)(.3,0)\psline[arrowscale=1.5]{->}(.3,-.25)(.3,-.05)
\psline(.5,0)(.5,-.3)\psline[arrowscale=1.5]{->}(.5,-.05)(.5,-.25)
\psline(1.2,0)(1.2,-.3)\psline[arrowscale=1.5]{->}(1.2,-.05)(1.2,-.25)
\rput[br](1.125,-.25){$\cdots$}
\rput[b](-1,1.7){$a$}\rput[b](1,1.7){$b$}
\rput[t](-1,-1.7){$a$}\rput[t](1,-1.7){$b$}
\psframe[linecolor=black](-1.5,-1.1)(-.5,-.9)
\psframe[linecolor=black]( 1.5,-1.1)( .5,-.9)
\psframe[linecolor=black](-1.5, 1.1)(-.5, .9)
\psframe[linecolor=black]( 1.5, 1.1)( .5, .9)
\psframe[linecolor=black](-1.3,-.5)(1.3, -.3)
\psframe[linecolor=darkgray,linewidth=1.2pt](.2,0)(1.3, .2)
\end{pspicture}\end{pspicture}\nonumber\\ \nonumber
&= \begin{pspicture}[shift=-1.9](-1.9,-2)(2.2,2)
\psline(-1,-.3)(-1,.9)\psline[arrowscale=1.5]{->}(-1,.2)(-1,.4)
\psline(-1,-.9)(-1,-.5)\psline[arrowscale=1.5]{->}(-1,-.8)(-1,-.6)
\psline(1,.9)(1,-.3)\psline[arrowscale=1.5]{->}(1,.4)(1,.2)
\psline(1,-.5)(1,-.9)\psline[arrowscale=1.5]{->}(1,-.6)(1,-.8)
\psline(1.4,.9)(1.4,-.9)\psline[arrowscale=1.5]{->}(1.4,.1)(1.4,-.1)
\psline(-1,-1.5)(-1,-1.1)\psline[arrowscale=1.5]{->}(-1,-1.4)(-1,-1.2)
\psline(-1,1.1)(-1,1.5)\psline[arrowscale=1.5]{->}(-1,1.2)(-1,1.4)
\psline(1,-1.1)(1,-1.5)\psline[arrowscale=1.5]{->}(1,-1.2)(1,-1.4)
\psline(1,1.5)(1,1.1)\psline[arrowscale=1.5]{->}(1,1.4)(1,1.2)
\rput[b](-1,1.7){$a$}\rput[b](1,1.7){$b$}
\rput[br](-1.2,.2){$a$}\rput[br](.8,.2){$b-1$}
\rput[t](-1,-1.7){$a$}\rput[t](1,-1.7){$b$}
\psframe[linecolor=black](-1.5,-1.1)(-.5,-.9)
\psframe[linecolor=black]( 1.5,-1.1)( .5,-.9)
\psframe[linecolor=black](-1.5, 1.1)(-.5, .9)
\psframe[linecolor=black]( 1.5, 1.1)( .5, .9)
\psframe[linecolor=black](-1.3,-.5)(1.3, -.3)
\end{pspicture} -\frac{[a][a+1]}{[a+b+1][a+b]} \begin{pspicture}[shift=-1.9](0,2)(0,2)\begin{pspicture}[shift=-1.9](-2.4,-2)(2.2,2)
\psline(-1,-.3)(-1,.9) \psline[arrowscale=1.5]{->}(-1,.2)(-1,.4)
\rput[br](-1.2,0){$a-1$} \psline(-1,-.9)(-1,-.5)
\psline[arrowscale=1.5]{->}(-1,-.8)(-1,-.6) \psline(1,.9)(1,-.3)
\psline[arrowscale=1.5]{->}(1,.7)(1,.5) \rput[bl](1.2,0){$b-1$}
\psline(1,-.5)(1,-.9) \psline[arrowscale=1.5]{->}(1,-.6)(1,-.8)
\pccurve[angleA=30,angleB=150,ncurv=1](-.7,-.9)(.7,-.9)\middlearrow
\pccurve[angleA=210,angleB=330,ncurv=1](.7,.9)(-.7,.9)\middlearrow
\rput[b](-1,1.7){$a$} \rput[b](1,1.7){$b$} \rput[t](-1,-1.7){$a$}
\rput[t](1,-1.7){$b$} \psline(-1,-1.5)(-1,-1.1)
\psline[arrowscale=1.5]{->}(-1,-1.4)(-1,-1.2)
\psline(-1,1.1)(-1,1.5) \psline[arrowscale=1.5]{->}(-1,1.2)(-1,1.4)
\psline(1,-1.1)(1,-1.5) \psline[arrowscale=1.5]{->}(1,-1.2)(1,-1.4)
\psline(1,1.5)(1,1.1) \psline[arrowscale=1.5]{->}(1,1.4)(1,1.2)
\psframe[linecolor=black](-1.5,-1.1)(-.5,-.9)
\psframe[linecolor=black]( 1.5,-1.1)( .5,-.9)
\psframe[linecolor=black](-1.5, 1.1)(-.5, .9)
\psframe[linecolor=black]( 1.5, 1.1)( .5, .9)
\psframe[linecolor=black](-1.3,-.5)(1.3, -.3)
\end{pspicture}\end{pspicture}
\end{eqnarray}
 \caption{Induction step for the proof of Proposition \ref{a2abexpprop1}. } \label{a2abexplem5} \end{figure}
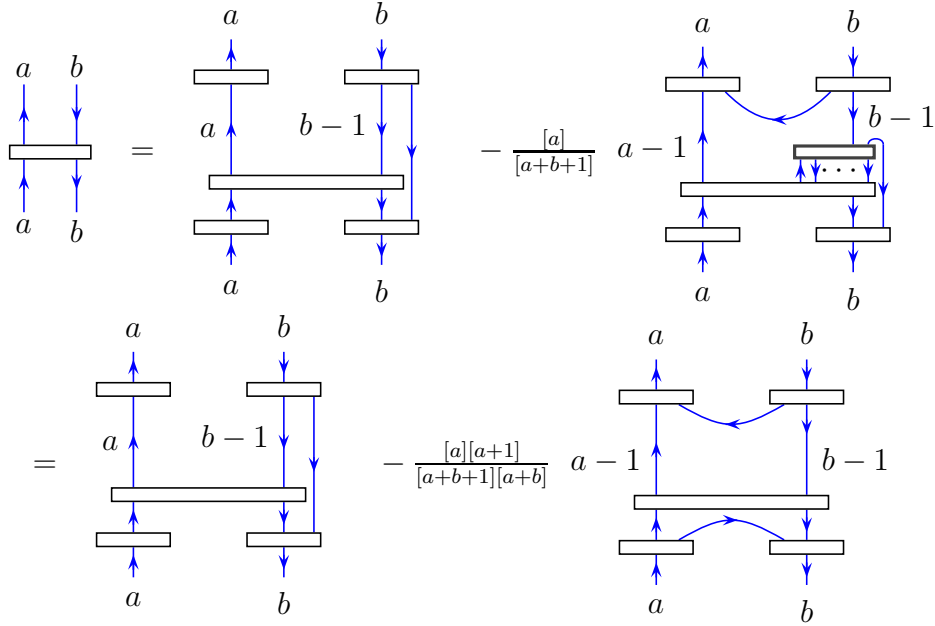

Let us denote the $k$-th term in the right-hand side of equation by
$D(k)$. We induct on $a+b$. It is clear for $a=0$ or $b=0$. If
$a\neq 0 \neq b$ then we use a segregated single clasp expansion of
weight $(a,b)$ in the middle for the first equality. Even if we do
not use the entire single clasp expansion of a segregated clasp,
once we attach clasps of weight $(a,0), (0,b)$ on the top, there are
only two nonzero webs which are webs with just one $U$-turn. One of
resulting webs has some $H$'s as in Figure~\ref{a2abexplem5} but if
we push them down to the clasp of weight $(a,b-1)$ in the middle, it
becomes a non-segregated clasp.  For the second equality we use a
non-segregated single clasp expansion at the clasp of weight
$(a,b-1)$ for which clasps of weight $(a-1,b-1)$ are located at
northeast corner. By the induction hypothesis, we have

\begin{align*}
&=\sum_{k=0}^{b-1}(-1)^{k}\frac{[a]![b-1]![a+b-k]!}{[a-k]![b-1-k]![k]![a+b]!}D(k)\\
&-\frac{[a+1][a]}{[a+b+1][a+b]}\sum_{k=0}^{b-1}(-1)^k
\frac{[a-1]![b-1]![a+b-1-k]!}{[a-1-k]![b-1-k]![k]![a+b-a]!}D(k+1)\\
&=1\cdot
D(0)+\sum_{k=1}^{b-1}((-1)^k\frac{[a]![b-1]![a+b-k]!}{[a-k]![b-1-k]![k]![a+b]!}\\
&+(-1)^k\frac{[a+1]![b-1]![a+b-k]!}{[a-k]![b-1-k]![k-1]![a+b]!})
D(k)\\
&-(-1)^{b-1}\frac{[a+1][a]}{[a+b+1][a+b]}\frac{[a-1]![b-1]![a]!}
{[a-b]![0]![b-1]![a+b-1]!} D(b)\\
&=D(0)+\sum_{k=1}^{b-1}(-1)^k\frac{[a]![b]![a+b+1-k]!}{[a-k]![b-k]![k]![a+b+1]!}
(\frac{[b-k][a+b+1]+[k][a+1]}{[b][a+b+1-k]})D(k)\\
&+(-1)^b\frac{[a]![b-1]![a+1]!}{[a-b]![0]![b-1]![a+b+1]!}D(b)\\
&=\sum_{k=0}^{b}(-1)^k\frac{[a]![b]![a+b+1-k]!}{[a-k]![b-k]![k]![a+b+1]!}D(k)
\end{align*}
\end{proof}

\section{Single Clasp Expansion for the quantum $\mathfrak{sp}(4)$}

\label{singlesp4}

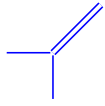
\begin{figure}
$$
\begin{pspicture}[shift=-.9](-1,-1)(1,1) \pnode(.65;180){a1}
\rput(0,0){\rnode{a2}{$$}} \pnode(.65;270){a3} \pnode(.9;45){a4}
\ncline[nodesepA=1pt,nodesepB=1pt]{a1}{a2}
\ncline[nodesepA=1pt,nodesepB=1pt]{a2}{a3}
\ncline[doubleline=true]{a2}{a4}
\end{pspicture}
$$
\caption{The generator of the quantum $\mathfrak{sp}(4)$ web space.}
\label{b2gen}
\end{figure}

The quantum $\mathfrak{sp}(4)$ webs are generated by a single web in
Figure \ref{b2gen} and a complete set of relations is given in
Figure \ref{b2rel} \cite{Kuperberg:spiders}. Again, an algebraic
proof of the existence of the clasp of the weight $(a,b)$ using the
annihilation axiom and the idempotent axiom is given in
\cite{Kuperberg:spiders}. On the other hand, one can use the double
clasps expansions in Corollary~\ref{b2n0expcor} and
Corollary~\ref{b20nexpcor} to define the clasps of the weight
$(n,0)$ and $(0,n)$. Unfortunately, we do not have any expansion
formula for the clasp of the weight $(a,b)$ where $a\neq 0\neq b$.
Using these expansions one can find Example~\ref{webexa2}. We can
define tetravalent vertices to achieve the same end as in
Figure~\ref{b2gen2}. We will use the these shapes to find a single
clasp expansion otherwise there is an ambiguity of a preferred
direction by the last relation presented in Figure~\ref{b2rel}.

First we will find single clasp expansions of clasps of weight
$(n,0)$ and $(0,n)$ and then use them to find coefficients of double
clasps expansions of clasps of weight $(n,0)$ and $(0,n)$. But we
are unable to find a single clasp expansion of the clasp of weight
$(a,b)$ where $a\neq 0 \neq b$. Remark that the cut weight is
defined slightly different way. A cut path may cut diagonally
through a tetravalent vertex, and its weight is defined as
$n\lambda_1 + (k + k')\lambda_2$, where $n$ is the number of type
``1'', single strands, that it cuts, $k$ is the number of type
``2'', double strands, that it cuts, and $k'$ is the number of
tetravalent vertices that it bisects.

\begin{figure}
\begin{eqnarray}
\begin{pspicture}[shift=-.5](-.6,-.6)(.6,.6) \pscircle(0,0){.4}
\end{pspicture}  &=   \dfrac{[6][2]}{[3]}
 \hskip 1.24cm , \hskip 1.24cm
\begin{pspicture}[shift=-.5](-.6,-.5)(.6,.5) \pscircle[doubleline=true](0,0){.4}
\end{pspicture}  &=   \dfrac{[6][5]}{[3][2]}
\nonumber \\
\begin{pspicture}[shift=-.5](-.6,-.5)(.6,.5) \psbezier(0,0)(.7,.7)(.7,-.7)(0,0)
\psline[doubleline=true](-.5,0)(0,0)
\end{pspicture}
&= 0 \hskip 1.55cm  , \hskip 1.55cm \begin{pspicture}[shift=-.5](-.8,-.5)(.8,.5)
\pcarc[arcangle=45](-.3,0)(.3,0) \pcarc[arcangle=-45](-.3,0)(.3,0)
\psline[doubleline=true](-.7,0)(-.3,0)
\psline[doubleline=true](.3,0)(.7,0)
\end{pspicture}
&=  -[2]^2\begin{pspicture}[shift=-.5](-.6,-.5)(.6,.5)
\psline[doubleline=true](-.5,0)(.5,0)
\end{pspicture} \nonumber \\
\begin{pspicture}[shift=-.8](-.9,-.9)(.9,.9) \pcarc[arcangle=-15](.4;90)(.4;210)
\pcarc[arcangle=-15](.4;210)(.4;330)
\pcarc[arcangle=-15](.4;330)(.4;90)
\psline[doubleline=true](.4;90)(.8;90)
\psline[doubleline=true](.4;210)(.8;210)
\psline[doubleline=true](.4;330)(.8;330)
\end{pspicture}
&= 0 \hskip .675cm  , \hskip .675cm \btwohvert - \btwohhoriz  &= \hh
- \vv \nonumber
\end{eqnarray}
 \caption{A complete set of relations of the quantum
$\mathfrak{sp}(4)$ web space.} \label{b2rel}
\end{figure}

\begin{exa}\label{webexa2}  The complete
expansions of the clasps of weight $(2,0)$ and $(3,0)$ are
$$
\begin{pspicture}[shift=-.6](-.5,-.7)(.5,.7) \psline(-.2,.5)(-.2,.1)
\psline(-.2,-.1)(-.2,-.5) \psline(.2,.5)(.2,.1)
\psline(.2,-.1)(.2,-.5) \psframe[linecolor=black](-.4,-.1)(.4,.1)
\end{pspicture}
 =   \begin{pspicture}[shift=-.6](-.5,-.7)(.5,.7)
\psline(-.2,.5)(-.2,-.5) \psline(.2,.5)(.2,-.5)
\end{pspicture}  +
\frac{1}{[2]^2} \begin{pspicture}[shift=-.4](-.5,-.5)(.5,.5)
\psline(-.3,.5)(.3,-.5) \psline(.3,.5)(-.3,-.5)
\end{pspicture}
+ \frac{[3][4]}{[2]^2[6]} \begin{pspicture}[shift=-.4](-.5,-.5)(.5,.5)
\psarc(0,.5){.3}{180}{360} \psarc(0,-.5){.3}{0}{180}
\end{pspicture} ,
$$
\begin{align*} \begin{pspicture}[shift=-.6](-.6,-.7)(.6,.7) \psline(-.4,.5)(-.4,.1)
\psline(-.4,-.1)(-.4,-.5) \psline(.4,.5)(.4,.1)
\psline(.4,-.1)(.4,-.5) \psline(0,.5)(0,.1) \psline(0,-.1)(0,-.5)
\psframe[linecolor=black](-.6,-.1)(.6,.1)
\end{pspicture}
 &=   \begin{pspicture}[shift=-.6](-.6,-.7)(.6,.7)
\psline(-.4,.5)(-.4,-.5) \psline(.4,.5)(.4,-.5) \psline(0,.5)(0,-.5)
\end{pspicture}  + \frac{1}{[3]} \left(
 \begin{pspicture}[shift=-.6](-.6,-.7)(.6,.7)
\psline(-.4,.5)(-.4,-.5) \psline(.4,.5)(0,-.5)
 \psline(0,.5)(.4,-.5)
\end{pspicture} +
 \begin{pspicture}[shift=-.6](-.6,-.7)(.6,.7)
\psline(.4,.5)(.4,-.5) \psline(-.4,.5)(0,-.5)
 \psline(0,.5)(-.4,-.5)
\end{pspicture}
\right)+   \frac{1}{[2]^2[3]} \left(
 \begin{pspicture}[shift=-.6](-.6,-.7)(.6,.7)
\psline(-.4,.5)(.4,-.5) \psline(0,.5)(-.4,-.5)
 \psline(.4,.5)(0,-.5)
\end{pspicture} +
 \begin{pspicture}[shift=-.6](-.6,-.7)(.6,.7)
\psline(-.4,.5)(0,-.5) \psline(0,.5)(.4,-.5)
 \psline(.4,.5)(-.4,-.5)
\end{pspicture}
\right) \\ &+  \frac{[4]^2}{[2][3][8]}\left(
\begin{pspicture}[shift=-.6](-.6,-.7)(.6,.7)  \psarc(-.2,.5){.2}{180}{0}
\psarc(.2,-.5){.2}{0}{180} \psline(.4,.5)(-.4,-.5)
\end{pspicture} + \begin{pspicture}[shift=-.6](-.6,-.7)(.6,.7)
 \psarc(.2,.5){.2}{180}{0}
\psarc(-.2,-.5){.2}{0}{180} \psline(-.4,.5)(.4,-.5)
\end{pspicture}\right) + +\frac{[2][4][6]}{[3]^2[8]} \left( \begin{pspicture}[shift=-.6](-.6,-.7)(.6,.7)
 \psarc(.2,.5){.2}{180}{0}
\psarc(.2,-.5){.2}{0}{180} \psline(-.4,.5)(-.4,-.5)
\end{pspicture} + \begin{pspicture}[shift=-.6](-.6,-.7)(.6,.7)
 \psarc(-.2,.5){.2}{180}{0}
\psarc(-.2,-.5){.2}{0}{180} \psline(.4,.5)(.4,-.5)
\end{pspicture}\right)
\\ &+ \frac{[4][6]}{[2][3]^2[8]}\left(
\begin{pspicture}[shift=-.6](-.6,-.7)(.6,.7)  \psarc(0,.5){.4}{180}{0}
\psarc(.2,-.5){.2}{0}{180} \psline(0,.5)(-.4,-.5)
\end{pspicture} + \begin{pspicture}[shift=-.6](-.6,-.7)(.6,.7)
 \psarc(.2,.5){.2}{180}{0}
\psarc(0,-.5){.4}{0}{180} \psline(-.4,.5)(0,-.5)
\end{pspicture} +
\begin{pspicture}[shift=-.6](-.6,-.7)(.6,.7)  \psarc(0,.5){.4}{180}{0}
\psarc(-.2,-.5){.2}{0}{180} \psline(0,.5)(.4,-.5)
\end{pspicture} + \begin{pspicture}[shift=-.6](-.6,-.7)(.6,.7)
 \psarc(-.2,.5){.2}{180}{0}
\psarc(0,-.5){.4}{0}{180} \psline(.4,.5)(0,-.5)
\end{pspicture}\right) +
\frac{[4][6]}{[2]^3[3]^2[8]} \begin{pspicture}[shift=-.6](-.6,-.7)(.6,.7)
\psarc(0,.5){.4}{180}{0} \psarc(0,-.5){.4}{0}{180}
\psline(0,.5)(0,-.5)
\end{pspicture} .
\end{align*}
\end{exa}

\begin{figure}
$$
\begin{pspicture}[shift=-.7](-.8,-.8)(.8,.8) \pnode(.7;45){a1} \pnode(.7;135){a2}
\pnode(.7;225){a3} \pnode(.7;315){a4} \ncline{a1}{a3}
\ncline{a2}{a4} \end{pspicture} = \begin{pspicture}[shift=-.7](-.8,-.8)(.8,.8)
\pnode(.7;45){a1} \pnode(.7;135){a2} \pnode(.7;225){a3}
\pnode(.7;315){a4} \pnode(.2;180){b1} \pnode(.2;0){b2}
\ncline{a1}{b2} \ncline{a2}{b1} \ncline{a3}{b1} \ncline{a4}{b2}
\ncline[doubleline=true]{b1}{b2}
\end{pspicture}
- \begin{pspicture}[shift=-.7](-.8,-.8)(.8,.8) \pnode(.7;45){a1}
\pnode(.7;135){a2} \pnode(.7;225){a3} \pnode(.7;315){a4}
\nccurve[angleA=210,angleB=150,ncurv=1]{a1}{a4}
\nccurve[angleA=-30,angleB=30,ncurv=1]{a2}{a3}
\end{pspicture} \hskip 1cm  ,  \hskip 1cm
\begin{pspicture}[shift=-.7](-.8,-.8)(.8,.8) \pnode(.7;45){a1} \pnode(.7;135){a2}
\pnode(.7;225){a3} \pnode(.7;315){a4} \rput(0,0){\rnode{b0}{$$}}
\ncline[doubleline=true,nodesep=1pt]{a1}{b0}
\ncline[doubleline=true,nodesep=1pt]{a2}{b0}
\ncline[doubleline=true,nodesep=1pt]{a3}{b0}
\ncline[doubleline=true,nodesep=1pt]{a4}{b0}
\end{pspicture}
 =  \begin{pspicture}[shift=-.7](-.8,-.8)(.8,.8) \pnode(.7;45){a1}
\pnode(.7;135){a2} \pnode(.7;225){a3} \pnode(.7;315){a4}
\pnode(.25;45){b1} \pnode(.25;135){b2} \pnode(.25;225){b3}
\pnode(.25;315){b4} \ncline{b1}{b2} \ncline{b4}{b1} \ncline{b3}{b2}
\ncline{b4}{b3} \ncline[doubleline=true]{a1}{b1}
\ncline[doubleline=true]{a2}{b2} \ncline[doubleline=true]{a3}{b3}
\ncline[doubleline=true]{a4}{b4}
\end{pspicture}
$$
 \caption{Tetravalent vertices.} \label{b2gen2}
\end{figure}

We demonstrate Theorem~\ref{b2n0expthm} for $n=2, 3$ by using the
presentations of clasps in Example~\ref{webexa2}. For $n=2$,
Theorem~\ref{b2n0expthm} is identical to the first formula of
Example~\ref{webexa2}. For $n=3$, we first attach the clasp of
weight $(2,0)$ to the southwest corner of each web in the second
formula of Example~\ref{webexa2}. Since $
\begin{pspicture}[shift=-.9](-.2,-.3)(1.2,1.6) \psline(.3,0)(.3,.4)
\psline(.3,.6)(.7,1.4) \psline(.7,0)(.7,.4) \psline(.7,.6)(.3,1.4)
\psframe[linecolor=black](0,.4)(1,.6)
\end{pspicture} = 0$ and $
\begin{pspicture}[shift=-.9](-.2,-.3)(1.2,1.6) \psline(.3,0)(.3,.4)
\psarc(.5,.6){.3}{0}{180} \psline(.7,0)(.7,.4)
\psframe[linecolor=black](0,.4)(1,.6)
\end{pspicture} = 0$, we find

$$ \begin{pspicture}[shift=-.8](-.5,-1.2)(.4,.6) \psframe[linecolor=black](-.6,-.7)(.2,-.5)
\qline(-.4,-.7)(-.4,-1.1) \qline(0,-.7)(0,-1.1)
\qline(.4,-.5)(.4,-1.1) \psline(-.4,.5)(-.4,.1)
\psline(-.4,-.1)(-.4,-.5) \psline(.4,.5)(.4,.1)
\psline(.4,-.1)(.4,-.5) \psline(0,.5)(0,.1) \psline(0,-.1)(0,-.5)
\psframe[linecolor=black](-.6,-.1)(.6,.1)
\end{pspicture}
=   \begin{pspicture}[shift=-.8](-.6,-1.2)(.45,.6)
\psframe[linecolor=black](-.6,-.7)(.2,-.5) \qline(-.4,-.7)(-.4,-1.1)
\qline(0,-.7)(0,-1.1) \qline(.4,-.5)(.4,-1.1)
\psline(-.4,.5)(-.4,-.5) \psline(.4,.5)(.4,-.5) \psline(0,.5)(0,-.5)
\end{pspicture}  + \frac{1}{[3]}
 \begin{pspicture}[shift=-.8](-.6,-1.2)(.45,.7)  \psframe[linecolor=black](-.6,-.7)(.2,-.5)
\qline(-.4,-.7)(-.4,-1.1) \qline(0,-.7)(0,-1.1)
\qline(.4,-.5)(.4,-1.1) \psline(-.4,.5)(-.4,-.5)
\psline(.4,.5)(0,-.5)
 \psline(0,.5)(.4,-.5)
\end{pspicture} +   \frac{1}{[2]^2[3]}
 \begin{pspicture}[shift=-.8](-.6,-1.2)(.45,.7) \psframe[linecolor=black](-.6,-.7)(.2,-.5)
\qline(-.4,-.7)(-.4,-1.1) \qline(0,-.7)(0,-1.1)
\qline(.4,-.5)(.4,-1.1) \psline(-.4,.5)(.4,-.5)
\psline(0,.5)(-.4,-.5)
 \psline(.4,.5)(0,-.5)
\end{pspicture}  + \frac{[2][4][6]}{[3]^2[8]} \begin{pspicture}[shift=-.8](-.6,-1.2)(.45,.7)
\psframe[linecolor=black](-.6,-.7)(.2,-.5) \qline(-.4,-.7)(-.4,-1.1)
\qline(0,-.7)(0,-1.1) \qline(.4,-.5)(.4,-1.1)
 \psarc(.2,.5){.2}{180}{0}
\psarc(.2,-.5){.2}{0}{180} \psline(-.4,.5)(-.4,-.5)
\end{pspicture}  + \frac{[4][6]}{[2][3]^2[8]}
\begin{pspicture}[shift=-.8](-.6,-1.2)(.45,.7)
\psframe[linecolor=black](-.6,-.7)(.2,-.5) \qline(-.4,-.7)(-.4,-1.1)
\qline(0,-.7)(0,-1.1) \qline(.4,-.5)(.4,-1.1)
\psarc(0,.5){.4}{180}{0} \psarc(.2,-.5){.2}{0}{180}
\psline(0,.5)(-.4,-.5)
\end{pspicture}
+  \frac{[4]^2}{[2][3][8]} \begin{pspicture}[shift=-.8](-.6,-1.2)(.45,.7)
\psframe[linecolor=black](-.6,-.7)(.2,-.5) \qline(-.4,-.7)(-.4,-1.1)
\qline(0,-.7)(0,-1.1) \qline(.4,-.5)(.4,-1.1)
\psarc(-.2,.5){.2}{180}{0} \psarc(.2,-.5){.2}{0}{180}
\psline(.4,.5)(-.4,-.5)
\end{pspicture}
$$

We can confirm these coefficients are the same as given in
Theorem~\ref{b2n0expthm}.

$$ \begin{pspicture}[shift=-.9](-.6,-1.2)(.6,.7)
\qline(-.4,-.5)(-.4,-1.1) \qline(0,-.5)(0,-1.1)
\qline(.4,-.5)(.4,-1.1) \psline(-.4,.5)(-.4,.1)
\psline(-.4,-.1)(-.4,-.5) \psline(.4,.5)(.4,.1)
\psline(.4,-.1)(.4,-.5) \psline(0,.5)(0,.1) \psline(0,-.1)(0,-.5)
\psframe[linecolor=black](-.6,-.1)(.6,.1)
\end{pspicture}
=   a_{01}\begin{pspicture}[shift=-.9](-.65,-1.2)(.45,.7)
\psframe[linecolor=black](-.6,-.7)(.2,-.5) \qline(-.4,-.7)(-.4,-1.1)
\qline(0,-.7)(0,-1.1) \qline(.4,-.5)(.4,-1.1)
\psline(-.4,.5)(-.4,-.5) \psline(.4,.5)(.4,-.5) \psline(0,.5)(0,-.5)
\end{pspicture}  + a_{02}
 \begin{pspicture}[shift=-.9](-.65,-1.2)(.45,.7) \psframe[linecolor=black](-.6,-.7)(.2,-.5)
\qline(-.4,-.7)(-.4,-1.1) \qline(0,-.7)(0,-1.1)
\qline(.4,-.5)(.4,-1.1) \psline(-.4,.5)(-.4,-.5)
\psline(.4,.5)(0,-.5)
 \psline(0,.5)(.4,-.5)
\end{pspicture} +   a_{03}
 \begin{pspicture}[shift=-.9](-.65,-1.2)(.45,.7) \psframe[linecolor=black](-.6,-.7)(.2,-.5)
\qline(-.4,-.7)(-.4,-1.1) \qline(0,-.7)(0,-1.1)
\qline(.4,-.5)(.4,-1.1) \psline(-.4,.5)(.4,-.5)
\psline(0,.5)(-.4,-.5)
 \psline(.4,.5)(0,-.5)
\end{pspicture}  + a_{12} \begin{pspicture}[shift=-.9](-.65,-1.2)(.45,.7)
\psframe[linecolor=black](-.6,-.7)(.2,-.5) \qline(-.4,-.7)(-.4,-1.1)
\qline(0,-.7)(0,-1.1) \qline(.4,-.5)(.4,-1.1)
 \psarc(.2,.5){.2}{180}{0}
\psarc(.2,-.5){.2}{0}{180} \psline(-.4,.5)(-.4,-.5)
\end{pspicture}  + a_{13}
\begin{pspicture}[shift=-.9](-.65,-1.2)(.45,.7)
\psframe[linecolor=black](-.6,-.7)(.2,-.5) \qline(-.4,-.7)(-.4,-1.1)
\qline(0,-.7)(0,-1.1) \qline(.4,-.5)(.4,-1.1)
\psarc(0,.5){.4}{180}{0} \psarc(.2,-.5){.2}{0}{180}
\psline(0,.5)(-.4,-.5)
\end{pspicture}
+  a_{23} \begin{pspicture}[shift=-.9](-.65,-1.2)(.45,.7)
\psframe[linecolor=black](-.6,-.7)(.2,-.5) \qline(-.4,-.7)(-.4,-1.1)
\qline(0,-.7)(0,-1.1) \qline(.4,-.5)(.4,-1.1)
\psarc(-.2,.5){.2}{180}{0} \psarc(.2,-.5){.2}{0}{180}
\psline(.4,.5)(-.4,-.5)
\end{pspicture}
$$

Now, we state a single clasp expansion of the clasp of weight
$(n,0)$.

\begin{thm} For a positive integer $n$,
\begin{eqnarray}
\begin{pspicture}[shift=-1.2](-.8,-.3)(2,2.3)
\rput[t](1,-.1){$n$}\qline(1,0)(1,.4)
\psframe[linecolor=black](0,.4)(2,.6)
\qline(1,.6)(1,2)\rput[b](1,2.1){$n$}
\end{pspicture}
= \sum_{i=0}^{n-1}\sum_{j=i+1}^{n}
[2]^{i-j+1}\frac{[n+1][n-j+1][2n-2i+2]}{[n][2n+2][n-i+1]} \hskip
.2cm \begin{pspicture}[shift=-1.3](.1,-.3)(3.55,2.5) \rput[t](1,-.1){$n-1$}
\rput(3.2,2.3){$``1"$} \rput(3.3,-.3){$``0"$} \qline(1,0)(1,.4)
\psframe[linecolor=black](0,.4)(3.1,.6) \qline(.2,.6)(.2,2)
\qline(.6,.6)(.6,2) \qline(.8,.6)(.8,2) \rput(1,2.3){$``j"$}
\rput(2.3,2.3){$``i"$} \qline(3.1,2)(3.1,1.5) \qline(3.1,1.5)(2.7,1)
\qline(2.7,1)(2.7,.6) \qline(2.5,2)(2.5,1.5) \qline(2.5,1.5)(2.1,1)
\qline(2.1,1)(2.1,.6) \qline(2.7,2)(2.7,1.5) \qline(2.7,1.5)(2.3,1)
\qline(2.3,1)(2.3,.6) \psarc(3.1,.6){.2}{0}{180}
\qline(3.3,0)(3.3,.6) \qline(1.2,1.5)(1.2,2) \qline(1.2,1.5)(1,1)
\qline(1,1)(1,.6)\qline(1.4,1.5)(1.4,2) \qline(1.4,1.5)(1.2,1)
\qline(1.2,1)(1.2,.6) \qline(2,1.5)(2,2) \qline(2,1.5)(1.8,1)
\qline(1.8,1)(1.8,.6) \qline(1.8,1.5)(1.8,2) \qline(1.8,1.5)(1.6,1)
\qline(1.6,1)(1.6,.6) \qline(.8,1.7)(.8,2)
\rput(1,2){\rnode{a1}{$$}} \rput(2.3,2){\rnode{a6}{$$}}
\nccurve[angleA=-45,angleB=225,ncurv=1]{a1}{a6}
\end{pspicture} \label{b2n0exp}
\end{eqnarray} \label{b2n0expthm}
\end{thm}

\begin{proof} By combining with the weight diagram of
$V_{\lambda_1}^{\otimes n}$ and minimal cut paths, we can find a set
of nonzero webs for single clasp expansion of a clasp of weight
$(n,0)$ as in equation~\ref{b2n0exp}. Let us denote the web
corresponding to the $i$-th in the first summation and $j$-th in the
second summation by $D_{i,j}$ and its coefficient by $a_{ij}$. First
we will show that these webs are linearly independent. Suppose that
a linear combination of the right-hand side of the equation in
Figure~\ref{a2abexp} is zero for some choice of $a_{ij}$. It is
clear that $a_{i,i+1}=0$ by attaching the clasp of weight $(n-i,0)$
to left top of webs and the clasp of weight $(i,0)$ to the right top
of each webs. By attaching the clasp of weight $(n-j+1,0)$ to left
top of webs and the clasp of weight $(i,0)$ to the right top of
webs, inductively we can show $a_{ij}=0$ for all $j\ge i+1$. By
Lemma~\ref{lemb2dim1}, we know that the dimension of the web space
of $V_{\lambda_1}^{\otimes n+1}\otimes V_{(n-1)\lambda_1}$ is
$\frac{n(n+1)}{2}$. Thus, these webs in right hand side of the
equation form a basis.

Now we are set to finds $a_{i,j}$. For equations, we remark that the
relations of webs shown in Figure~\ref{b2n0help} can be easily
obtained from the relations depicted in Figure~\ref{b2rel}. Using
these relations, we get the following $n-1$ equations by attaching a
$\begin{pspicture}[shift=.2](-.17,-.17)(.17,.17)
\qline(-.15,-.15)(-.15,.05)\qline(.15,-.15)(.15,.05)
\pccurve[angleA=90,angleB=90,ncurv=1](-.15,.05)(.15,.05)
\end{pspicture}$.
By attaching a $\begin{pspicture}[shift=.2](-.17,-.17)(.17,.17)
\psline[doubleline=true](0,0)(0,.17)\psline[nodesep=1pt](-.15,-.15)(-.025,0)
\psline(.15,-.15)(0.025 ,0)\end{pspicture}$, we have $(n-1)^2$
equations. There are two special equations and four different shapes
of equation as follows.

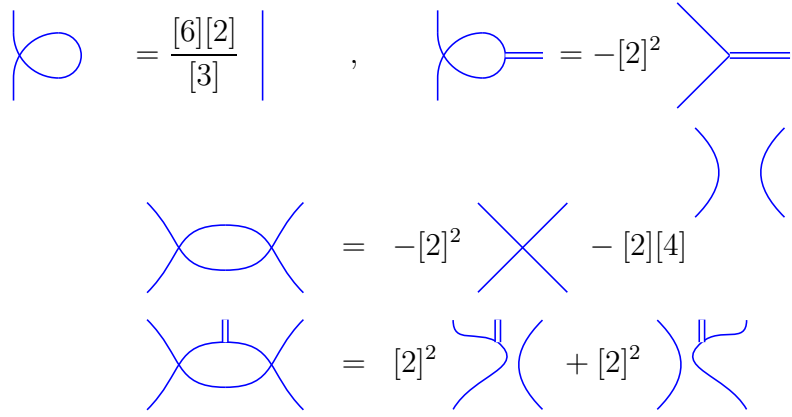
\begin{figure}
\begin{eqnarray}
\begin{pspicture}[shift=-.6](-.4,-.7)(1.2,.7) \psline(-.3,.6)(-.3,.3)
\psline(-.3,-.3)(-.3,-.6) \psarc(.3,.3){.6}{180}{270}
\psarc(.3,-.3){.6}{90}{180} \psarc(.3,0){.3}{270}{90}
\end{pspicture}
= \frac{[6][2]}{[3]} \hskip .3cm \begin{pspicture}[shift=-.6](-.1,-.7)(.1,.7)
\psline(-.1,.6)(-.1,-.6)
\end{pspicture}\hskip .5cm & , & \hskip .5cm
\begin{pspicture}[shift=-.6](-.4,-.7)(1.2,.7) \psline(-.3,.6)(-.3,.3)
\psline(-.3,-.3)(-.3,-.6) \psarc(.3,.3){.6}{180}{270}
\psarc(.3,-.3){.6}{90}{180} \psarc(.3,0){.3}{270}{90}
\psline[doubleline=true](.6,0)(1.1,0)
\end{pspicture}
= -[2]^2 \hskip .2cm \begin{pspicture}[shift=-.6](-.1,-.7)(1.5,.7)
\psline(-.1,.7)(.6,0) \psline(.6,0)(-.1,-.7)
\psline[doubleline=true](.6,0)(1.4,0)
\end{pspicture}\nonumber\\
\begin{pspicture}[shift=-.6](-1.2,-.7)(1.2,.7) \pnode(1.2;30){a1}
\pnode(1.2;150){a2} \pnode(1.2;210){a3} \pnode(1.2;330){a4}
\pnode(.6;90){b1} \pnode(.3;90){b2} \pnode(.3;270){b3}
\nccurve[angleA=315,angleB=180,ncurv=1]{a2}{b3}
\nccurve[angleA=225,angleB=0,ncurv=1]{a1}{b3}
\nccurve[angleA=135,angleB=0,ncurv=1]{a4}{b2}
\nccurve[angleA=45,angleB=180,ncurv=1]{a3}{b2}
\end{pspicture}
& = & -[2]^2\begin{pspicture}[shift=-.6](-.8,-.7)(.8,.7) \pnode(.84;45){a1}
\pnode(.84;135){a2} \pnode(.84;225){a3} \pnode(.84;315){a4}
\ncline{a2}{a4} \ncline{a3}{a1} \end{pspicture} -[2][4]
\begin{pspicture}[shift=.4](-.7,-.7)(.7,.7) \pnode(.84;45){a1}
\pnode(.84;135){a2} \pnode(.84;225){a3} \pnode(.84;315){a4}
\nccurve[angleA=135,angleB=225,ncurv=1]{a4}{a1}
\nccurve[angleA=45,angleB=-45,ncurv=1]{a3}{a2}
\end{pspicture}\nonumber\\
\begin{pspicture}[shift=-.6](-1.2,-.7)(1.2,.7) \pnode(1.2;30){a1}
\pnode(1.2;150){a2} \pnode(1.2;210){a3} \pnode(1.2;330){a4}
\pnode(.6;90){b1} \pnode(.3;90){b2} \pnode(.3;270){b3}
\nccurve[angleA=315,angleB=180,ncurv=1]{a2}{b3}
\nccurve[angleA=225,angleB=0,ncurv=1]{a1}{b3}
\nccurve[angleA=135,angleB=0,ncurv=1]{a4}{b2}
\nccurve[angleA=45,angleB=180,ncurv=1]{a3}{b2}
\ncline[doubleline=true]{b1}{b2}
\end{pspicture}
& = & [2]^2\begin{pspicture}[shift=-.6](-.8,-.7)(.8,.7) \pnode(.84;45){a1}
\pnode(.84;135){a2} \pnode(.84;225){a3} \pnode(.84;315){a4}
\pnode(.6;90){b1} \pnode(.3;90){b2} \pnode(.3;270){b3}
\nccurve[angleA=-90,angleB=150,ncurv=1]{a2}{b2}
\nccurve[angleA=60,angleB=-40,ncurv=1]{a3}{b2}
\nccurve[angleA=135,angleB=225,ncurv=1]{a4}{a1}
\ncline[doubleline=true]{b1}{b2}
\end{pspicture}  + [2]^2 \begin{pspicture}[shift=-.6](-.8,-.7)(.8,.7) \pnode(.84;45){a1}
\pnode(.84;135){a2} \pnode(.84;225){a3} \pnode(.84;315){a4}
\pnode(.6;90){b1} \pnode(.3;90){b2} \pnode(.3;270){b3}
\nccurve[angleA=315,angleB=45,ncurv=1]{a2}{a3}
\nccurve[angleA=270,angleB=40,ncurv=1]{a1}{b2}
\nccurve[angleA=120,angleB=220,ncurv=1]{a4}{b2}
\ncline[doubleline=true]{b1}{b2}
\end{pspicture} \nonumber
\end{eqnarray}
 \caption{Useful relations of webs for Theorem \ref{b2n0expthm}.  } \label{b2n0help}
\end{figure}

$$a_{n-2,n-1}+\frac{[2][6]}{[3]} a_{n-2,n} - \frac{[2][6]}{[3]} a_{n-1,n}=0,$$

$$-\frac{[2][6]}{[3]} a_{12} + \frac{[2][6]}{[3]} a_{13}
+ a_{23}+ 1+ \frac{[2][6]}{[3]} b_2 -[2][4] b_3=0.$$ Type I : For
$i=1, 2, \ldots, n-3$,

$$a_{i,i+1}+\frac{[2][6]}{[3]} a_{i,i+2} -[2][4]a_{i,i+3}  -
\frac{[2][6]}{[3]} a_{i+1,i+2} +\frac{[2][6]}{[3]}
a_{i+1,i+3}+a_{i+2,i+3}=0.$$ Type II : For $i=0, 1, \ldots, n-2$,

$$a_{i,n-1}-[2]^2a_{i,n}=0.$$
Type III : For $i= 0, 1, 2, \ldots, n-3$, $k=2, 3, \ldots, n-i-1$,

$$a_{i,n-k}-[2]^2a_{i,n-k+1}+[2]^2a_{i,n-k+2}=0.$$
Type IV : For $i=3, 4, \ldots, n$, $k=n-i+3, n-i+4, \ldots, n$,

$$[2]^2 a_{n-k,i}-[2]^2a_{n-k+1,i}+a_{n-k+2,i}=0.$$
Then we check our answer satisfies the equations and it is clear
that $a_{0,1}=1$ by a normalization. Since these webs in the
equation~\ref{b2n0exp} form a basis, the coefficients are unique.
Therefore, it completes the proof.
\end{proof}

By attaching the clasp of weigh $(n-1,0)$ on the top of all webs in
the equation presented in equation~\ref{b2n0exp}, we find the double
clasp expansion of the clasp of weight $(n,0)$.

\begin{cor} For a positive integer $n$,
$$
\begin{pspicture}[shift=-.7](-.2,-.3)(1.2,1.3) \rput[t](.5,-.1){$n$}
\qline(.5,0)(.5,.4) \qline(.5,.6)(.5,1) \rput[b](.5,1.1){$n$}
\psframe[linecolor=black](0,.4)(1,.6)
\end{pspicture}
= \begin{pspicture}[shift=-.7](-.2,-.3)(1.4,1.3) \rput[t](.5,-0.1){$n-1$}
\qline(.5,0)(.5,.4) \qline(.5,.6)(.5,1) \rput[b](.5,1.1){$n-1$}
\qline(1.3,0)(1.3,1) \psframe[linecolor=black](0,.4)(1,.6)
\end{pspicture}
+ \frac{[2n][n+1][n-1]}{[2n+2][n][n]}
\begin{pspicture}[shift=-1.2](-.2,-.3)(1.45,2.3) \rput[t](.5,-.1){$n-1$}
\qline(.5,0)(.5,.4) \psframe[linecolor=black](0,.4)(1,.6)
\qline(.25,.6)(.25,1.4) \rput[l](.35,1){$n-2$} \qline(.5,1.6)(.5,2)
\rput[b](.5,2.1){$n-1$} \psarc(1,.6){.2}{0}{180}
\qline(1.2,0)(1.2,.6) \psarc(1,1.4){.2}{180}{0}
\qline(1.2,1.4)(1.2,2) \psframe[linecolor=black](0,1.4)(1,1.6)
\end{pspicture}
+ \frac{[n-1]}{[n][2]} \begin{pspicture}[shift=-1.2](-.2,-.3)(1.45,2.3)
\rput[t](.5,-.1){$n-1$} \qline(.5,0)(.5,.4) \qline(.25,.6)(.25,1.4)
\qline(.5,1.6)(.5,2) \qline(1.2,1.4)(1.2,2) \rput[b](.5,2.1){$n-1$}
\qline(.8,.6)(1.2,1.4) \qline(1.2,0)(1.2,.6) \qline(.8,1.4)(1.2,.6)
\psframe[linecolor=black](0,.4)(1,.6)
\psframe[linecolor=black](0,1.4)(1,1.6)
\end{pspicture}
$$
 \label{b2n0expcor}
\end{cor}

Then we look at the clasp of weight $(0,n)$. The main idea for the
clasp of weight $(n,0)$ works exactly same except we replace the
basis as shown in equation~\ref{b20nexp}. For the linear
independency, every idea of the proof of Theorem~\ref{b2n0expthm}
works with the fact $ \begin{pspicture}[shift=-.8](-.2,-.3)(1.2,1.6)
\psline[doubleline=true](.3,0)(.3,.4)
\psline[doubleline=true](.3,.6)(.3,1)
\psline[doubleline=true](.7,0)(.7,.4)
\psline[doubleline=true](.7,.6)(.7,1)
\psline(.3,1.4)(.3,1)(.7,1)(.7,1.4)
\psframe[linecolor=black](0,.4)(1,.6)
\end{pspicture} = 0$.
As we did for the clasp of weight $(n,0)$, we first find the
equations as illustrated in Figure~\ref{b20nhelp} for the next step.
The same as before, we set $a_{ij}$ be the coefficient of the web of
$(i,j)$ in the summation. By attaching a
$\begin{pspicture}[shift=-.12](-.17,-.17)(.17,.17)
\psline[doubleline=true](-.15,-.15)(-.15,.05)\psline[doubleline=true](.15,-.15)(.15,.05)
\pccurve[angleA=90,angleB=90,ncurv=1,doubleline=true](-.15,.05)(.15,.05)
\end{pspicture}$
and a $\begin{pspicture}[shift=-.12](-.17,-.17)(.17,.17)
\psline[doubleline=true](-.15,-.15)(-.15,.03)
\psline[doubleline=true](.15,-.15)(.15,.03)
\qline(-.15,.03)(-.15,.17) \qline(.15,.03)(0.15,.17)
\qline(-.147,.03)(.147,.03)\end{pspicture}$ , we get the following
equations and we can solve them successively as in
Theorem~\ref{b20nexpthm}.

$$a_{n-2,n-1}-[5][2]^2a_{n-2,n}+\frac{[6][5]}{[3][2]}a_{n-1,n}=0,$$
$$-[3][2]^2a_{n-2,n}+[5]a_{n-1,n}=0.$$
Type I : For $i=0, 1, \ldots, n-3$,

$$a_{i,i+1}-[5][2]^2a_{i,i+2}+[3][2]^4a_{i,i+3}
+\frac{[6][5]}{[3][2]}a_{i+1,i+2}-[5][2]^2a_{i+1,i+3}+a_{i+2,i+3}=0$$
Type II : For $i=0, 1, \ldots, n-2$,

$$a_{i,n-1}-[4][2]a_{i,n}=0.$$
Type III : For $i=0,1,\ldots ,n-3$ and $j=i+1,i+2,\ldots, n-2$,

$$a_{i,j}-[4][2]a_{i,j+1}+[2]^4a_{i,j+2}=0.$$
Type IV : For $i=0,1,\ldots, n-3$ and $j=i+3,i+4,\ldots, n$,

$$[2]^4a_{i,j}-[4][2]a_{i+1,j}+a_{i+2,j}=0.$$
Type V : For $i=1,2,\ldots, n-2$
$$-[3][2]^2a_{i-1,i+1}+[2]^4a_{i-1,i+2}+[5]a_{i,i+1}-[3][2]^2a_{i,i+2}=0.$$

\begin{figure}
\begin{eqnarray*}
\begin{pspicture}[shift=-.6](-.4,-.7)(.4,.7) \psline(-.3,.6)(0,.3)
\psline(0,-.3)(-.3,-.6)\psline(0,.3)(0,-.3)
\psarc[doubleline=true](0,0){.3}{270}{90}
\end{pspicture}
 = [5] \hskip .2cm \begin{pspicture}[shift=-.6](-.1,-.7)(.1,.7)
\psline(0,.6)(0,-.6)
\end{pspicture} &,& \hskip .2cm \begin{pspicture}[shift=-.6](-.4,-.7)(1.2,.7)
\psline[doubleline=true](-.3,.6)(-.3,.3)
\psline[doubleline=true](-.3,-.3)(-.3,-.6)
\psarc[doubleline=true,nodesep=1pt](.3,.3){.6}{180}{207}
\psarc[doubleline=true,nodesep=1pt](.3,.3){.6}{213}{270}
\psarc[doubleline=true,nodesep=1pt](.3,-.3){.6}{90}{147}
\psarc[doubleline=true,nodesep=1pt](.3,-.3){.6}{153}{180}
\psarc[doubleline=true](.3,0){.3}{270}{90}
\end{pspicture}
= -[2]^2[5] \hskip .6cm \begin{pspicture}[shift=-.6](-.1,-.7)(.1,.7)
\psline[doubleline=true](-.1,.6)(-.1,-.6)
\end{pspicture} \nonumber\\
\begin{pspicture}[shift=-.6](-.4,-.7)(1,.7) \psline(-.3,.6)(0.02,.32)
\psline(0.02,-.32)(-.3,-.6) \psline(.02,.32)(.02,-.32)
\rput(0,.3){\rnode{a1}{$$}} \rput(.9,-.6){\rnode{a2}{$$}}
\rput(0,-.3){\rnode{a3}{$$}} \rput(.9,.6){\rnode{a4}{$$}}
\rput(.5,0){\rnode{b0}{$$}}
\nccurve[doubleline=true,angleA=30,angleB=135,ncurv=1,nodesep=1pt]{a1}{b0}
\nccurve[doubleline=true,angleA=-45,angleB=135,ncurv=1,nodesep=1pt]{b0}{a2}
\nccurve[doubleline=true,angleA=-30,angleB=225,ncurv=1,nodesep=1pt]{a3}{b0}
\nccurve[doubleline=true,angleA=45,angleB=225,ncurv=1,nodesep=1pt]{b0}{a4}
\end{pspicture}
& = &-[2][4] \begin{pspicture}[shift=-.6](-.4,-.7)(.4,.7) \psline(-.3,.6)(0,.3)
\psline(0,-.3)(-.3,-.6) \psline(0,.3)(0,-.3)
\psline[doubleline=true](0,-.3)(.3,-.6)
\psline[doubleline=true](0,.3)(.3,.6)
\end{pspicture}
-[2]^2[3]\begin{pspicture}[shift=-.6](-1,-.7)(1,.7)
\rput(.7,-.6){\rnode{a1}{$$}} \rput(.7,.6){\rnode{a2}{$$}}
\rput(-.7,-.6){\rnode{a3}{$$}} \rput(-.7,.6){\rnode{a4}{$$}}
\nccurve[doubleline=true,angleA=135,angleB=225,ncurv=1]{a1}{a2}
\nccurve[angleA=45,angleB=-45,ncurv=1]{a3}{a4}
\end{pspicture}\nonumber\\
\begin{pspicture}[shift=-.6](-1.2,-.7)(1.2,.7) \pnode(1.2;30){a1}
\pnode(1.2;150){a2} \pnode(1.2;210){a3} \pnode(1.2;330){a4}
\pnode(.6;0){b1} \pnode(.6;180){b2}
\nccurve[doubleline=true,angleA=225,angleB=45,ncurv=1,nodesep=1pt]{a1}{b1}
\nccurve[doubleline=true,angleA=135,angleB=315,ncurv=1,nodesep=1pt]{a4}{b1}
\nccurve[doubleline=true,angleA=315,angleB=135,ncurv=1,nodesep=1pt]{a2}{b2}
\nccurve[doubleline=true,angleA=45,angleB=225,ncurv=1,nodesep=1pt]{a3}{b2}
\nccurve[doubleline=true,angleA=135,angleB=45,ncurv=1,nodesep=1pt]{b1}{b2}
\nccurve[doubleline=true,angleA=225,angleB=315,ncurv=1,nodesep=1pt]{b1}{b2}
\end{pspicture}
& = & -[2][4] \begin{pspicture}[shift=-.6](-.8,-.7)(.8,.7) \pnode(.84;45){a1}
\pnode(.84;135){a2} \pnode(.84;225){a3} \pnode(.84;315){a4}
\rput(0,0){\rnode{b0}{$$}}
\ncline[doubleline=true,nodesep=1pt]{a1}{b0}
\ncline[doubleline=true,nodesep=1pt]{a2}{b0}
\ncline[doubleline=true,nodesep=1pt]{a3}{b0}
\ncline[doubleline=true,nodesep=1pt]{a4}{b0}
\end{pspicture} +[2]^4[3] \begin{pspicture}[shift=-.6](-.7,-.7)(.7,.7)
\pnode(.84;45){a1} \pnode(.84;135){a2} \pnode(.84;225){a3}
\pnode(.84;315){a4}
\nccurve[doubleline=true,angleA=135,angleB=225,ncurv=1]{a4}{a1}
\nccurve[doubleline=true,angleA=45,angleB=-45,ncurv=1]{a3}{a2}
\end{pspicture}
\end{eqnarray*}
 \caption{Useful relations of webs for Theorem \ref{b20nexpthm}. } \label{b20nhelp}
\end{figure}

\begin{thm}
For a positive integer $n$,
\begin{eqnarray} \begin{pspicture}[shift=-1.2](-.8,-.3)(2.3,2.3) \rput[t](1,-.1){$n$}
\psline[doubleline=true](1,0)(1,.4)
\psline[doubleline=true](1,.6)(1,2) \rput[b](1,2.1){$n$}
\psframe[linecolor=black](0,.4)(2,.6)
\end{pspicture}
= \sum_{i=0}^{n-1}\sum_{j=i+1}^{n}
[2]^{2(1+i-j)}\frac{[2n+1-2i][2n-2j+2]}{[2n][2n+1]} \hskip .2cm
\begin{pspicture}[shift=-1.3](-.3,-.3)(3.55,2.5) \rput[t](1,-.1){$n-1$}
\psline[doubleline=true](1,0)(1,.4)
\psline[doubleline=true](.2,.6)(.2,2)
\psline[doubleline=true](.6,.6)(.6,2)
\psline[doubleline=true](.8,.6)(.8,2) \rput(1,2.3){$``j"$}
\rput(2.3,2.3){$``i"$} \rput(3.2,2.3){$``1"$} \rput(3.3,-.3){$``0"$}
\psline[doubleline=true](3.1,2)(3.1,1.4)
\psline[doubleline=true](3.1,1.4)(2.7,1)
\psline[doubleline=true](2.7,1)(2.7,.6)
\psline[doubleline=true](2.5,1.7)(2.5,1.4)
\psline[doubleline=true](2.5,2)(2.5,1.7)
\psline[doubleline=true](2.5,1.4)(2.1,1)
\psline[doubleline=true](2.1,1)(2.1,.6)
\psline[doubleline=true](2.7,2)(2.7,1.4)
\psline[doubleline=true](2.7,1.4)(2.3,1)
\psline[doubleline=true](2.3,1)(2.3,.6)
\psline[doubleline=true](3.3,0)(3.3,.6)
\psline[doubleline=true,nodesep=1pt](1.2,1.4)(1.2,1.67)
\psline[doubleline=true,nodesep=1pt](1.2,2)(1.2,1.75)
\psline[doubleline=true](1.2,1.4)(1,1)
\psline[doubleline=true](1,1)(1,.6)
\psline[doubleline=true,nodesep=1pt](1.4,1.4)(1.4,1.67)
\psline[doubleline=true,nodesep=1pt](1.4,2)(1.4,1.75)
\psline[doubleline=true](1.4,1.4)(1.2,1)
\psline[doubleline=true](1.2,1)(1.2,.6)
\psline[doubleline=true,nodesep=1pt](2,1.4)(2,1.67)
\psline[doubleline=true,nodesep=1pt](2,2)(2,1.75)
\psline[doubleline=true](2,1.4)(1.8,1)
\psline[doubleline=true](1.8,1)(1.8,.6)
\psline[doubleline=true,nodesep=1pt](1.8,1.75)(1.8,2)
\psline[doubleline=true,nodesep=1pt](1.8,1.4)(1.8,1.67)
\psline[doubleline=true](1.8,1.4)(1.6,1)
\psline[doubleline=true](1.6,1)(1.6,.6)
\psline[doubleline=true](.8,1.7)(.8,2)
\psline[doubleline=true,nodesep=1pt](1,2)(1,1.71)(1.18,1.71)
\psline[doubleline=true,nodesep=1pt](1.22,1.71)(1.38,1.71)
\psline[doubleline=true,nodesep=1pt](1.42,1.71)(1.78,1.71)
\psline[doubleline=true,nodesep=1pt](1.82,1.71)(1.98,1.71)
\psline[doubleline=true,nodesep=1pt](2.02,1.71)(2.3,1.71)(2.3,2)
\rput(1,2){\rnode{a1}{$$}} \rput(2.3,2){\rnode{a2}{$$}}
\psarc[doubleline=true](3.1,.6){.2}{0}{180}
\psframe[linecolor=black](0,.4)(3.1,.6)
\end{pspicture} \label{b20nexp}
\end{eqnarray}
\label{b20nexpthm}
\end{thm}

By attaching the clasp of weigh $(0,n-1)$ on the top of all webs
shown in equation~\ref{b20nexp} we find the double clasps expansion
of the clasp of weight $(0,n)$.

\begin{cor} For a positive integer $n$,
$$
\begin{pspicture}[shift=-.7](-.2,-.3)(1.2,1.3)
\rput[t](.5,-.1){$n$}\psline[doubleline=true](.5,0)(.5,.4)
\psline[doubleline=true](.5,.6)(.5,1)
\psframe[linecolor=black](0,.4)(1,.6) \rput[b](.5,1.1){$n$}
\end{pspicture}
= \begin{pspicture}[shift=-.7](-.2,-.3)(1.4,1.3)
\rput[t](.5,-0.1){$n-1$}\psline[doubleline=true](.5,0)(.5,.4)
\psline[doubleline=true](.5,.6)(.5,1)\rput[b](.5,1.1){$n-1$}
\psline[doubleline=true](1.3,0)(1.3,1)
\psframe[linecolor=black](0,.4)(1,.6)
\end{pspicture}
+ \frac{[2n-1][2n-2]}{[2n+1][2]} \begin{pspicture}[shift=-1.2](-.2,-.3)(1.45,2.3)
\rput[t](.5,-.1){$n-1$} \psline[doubleline=true](.5,0)(.5,.4)
\psline[doubleline=true](.25,.6)(.25,1.4)\rput[l](.35,1){$n-2$}
\psline[doubleline=true](.5,1.6)(.5,2) \rput[b](.5,2.1){$n-1$}
\psarc[doubleline=true](1,.6){.2}{0}{180}
\psline[doubleline=true](1.2,0)(1.2,.6)
\psarc[doubleline=true](1,1.4){.2}{180}{0}
\psline[doubleline=true](1.2,1.4)(1.2,2)
\psframe[linecolor=black](0,.4)(1,.6)
\psframe[linecolor=black](0,1.4)(1,1.6)
\end{pspicture}
+ \frac{[2n-2]}{[2n][2][2]}\begin{pspicture}[shift=-1.2](-.2,-.3)(1.45,2.3)
\rput[t](.5,-.1){$n-1$} \psline[doubleline=true](.5,0)(.5,.4)
\psline[doubleline=true](.25,.6)(.25,1.4)
\psline[doubleline=true](.5,1.6)(.5,2) \rput[b](.5,2.1){$n-1$}
\psline[doubleline=true,nodesep=1pt](.8,.6)(.97,.94)
\psline[doubleline=true,nodesep=1pt](1.04,1.06)(1.2,1.4)
\psline[doubleline=true](1.2,0)(1.2,.6)
\psline[doubleline=true,nodesep=1pt](.8,1.4)(.97,1.06)
\psline[doubleline=true,nodesep=1pt](1.03,.94)(1.2,.6)
\psline[doubleline=true](1.2,1.4)(1.2,2)
\qline(.942,.94)(.942,1.06)\qline(1.06,.94)(1.06,1.06)
\psframe[linecolor=black](0,.4)(1,.6)
\psframe[linecolor=black](0,1.4)(1,1.6)
\end{pspicture}
$$
\label{b20nexpcor}
\end{cor}

\section{Applications of the quantum $\mathfrak{sl}(3)$ representation theory}\label{application}

In the section we will discuss some applications of the quantum
$\mathfrak{sl}(3)$ representation theory.

\subsection{Polynomial invariants of links}

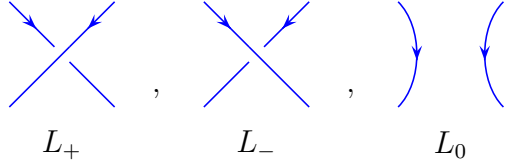
\begin{figure}
\begin{eqnarray}\nonumber \begin{pspicture}[shift=-1.2](0,-1.6)(0,1)\end{pspicture}
 \begin{pspicture}[shift=-1](0,1)(0,1)\begin{pspicture}[shift=-1.2](-1,-1.5)(1,1)
\rput(.7,.7){\rnode{a1}{$$}} \rput(-.7,.7){\rnode{a2}{$$}}
\rput(-.7,-.7){\rnode{a3}{$$}} \rput(.7,-.7){\rnode{a4}{$$}}
\rput(.1,.1){\rnode{b1}{$$}} \rput(-.1,.1){\rnode{b2}{$$}}
\rput(-.1,-.1){\rnode{b3}{$$}}
\rput(.1,-.1){\rnode{b4}{$$}}\ncline{a2}{b2}\middlearrow
\ncline{a1}{b1}\middlearrow \ncline{b4}{a4} \ncline{b1}{a3}
\rput(0,-1.2){\rnode{c4}{$L_+$}}
\end{pspicture}\end{pspicture}  \hskip .2cm, \hskip .2cm  \begin{pspicture}[shift=-1](0,1)(0,1)
\begin{pspicture}[shift=-1.2](-1,-1.5)(1,1) \rput(.7,.7){\rnode{a1}{$$}}
\rput(-.7,.7){\rnode{a2}{$$}} \rput(-.7,-.7){\rnode{a3}{$$}}
\rput(.7,-.7){\rnode{a4}{$$}} \rput(.1,.1){\rnode{b1}{$$}}
\rput(-.1,.1){\rnode{b2}{$$}} \rput(-.1,-.1){\rnode{b3}{$$}}
\rput(.1,-.1){\rnode{b4}{$$}} \ncline{a2}{b2}\middlearrow
\ncline{a1}{b1}\middlearrow \ncline{b2}{a4} \ncline{b3}{a3}
\rput(0,-1.2){\rnode{c4}{$L_-$}}
\end{pspicture}\end{pspicture}
\hskip .2cm , \hskip .2cm  \begin{pspicture}[shift=-1](0,1)(0,1)\begin{pspicture}[shift=-1.2](-1,-1.5)(1,1)
\rput(.7,.7){\rnode{a1}{$$}} \rput(-.7,.7){\rnode{a2}{$$}}
\rput(-.7,-.7){\rnode{a3}{$$}} \rput(.7,-.7){\rnode{a4}{$$}}
\rput(.1,.1){\rnode{b1}{$$}} \rput(-.1,.1){\rnode{b2}{$$}}
\rput(-.1,-.1){\rnode{b3}{$$}} \rput(.1,-.1){\rnode{b4}{$$}}
\nccurve[angleA=225,angleB=135]{a1}{a4}\middlearrow
\nccurve[angleA=315,angleB=45]{a2}{a3}\middlearrow
\rput(0,-1.2){\rnode{c4}{$L_0$}}
\end{pspicture}\end{pspicture}
\end{eqnarray}
\caption{The shape of $L_+, L_-$ and $L_0$.} \label{local}
\end{figure}

\begin{figure}
\begin{eqnarray}\nonumber \begin{pspicture}[shift=-.9](0,-1)(0,1)\end{pspicture}
\begin{pspicture}[shift=-.9](0,1)(0,1)\begin{pspicture}[shift=-.9](-1.1,-1)(1.1,1) \rput(.7,.7){\rnode{a1}{$$}}
\rput(-.7,.7){\rnode{a2}{$$}} \rput(-.7,-.7){\rnode{a3}{$$}}
\rput(.7,-.7){\rnode{a4}{$$}} \rput(.1,.1){\rnode{b1}{$$}}
\rput(-.1,.1){\rnode{b2}{$$}} \rput(-.1,-.1){\rnode{b3}{$$}}
\rput(.1,-.1){\rnode{b4}{$$}} \ncline{a2}{b2}\middlearrow
\ncline{a1}{b1}\middlearrow \ncline{b4}{a4} \ncline{b1}{a3}
\end{pspicture}\end{pspicture} &=& q^{\frac{1}{2}}
\begin{pspicture}[shift=-.9](0,1)(0,1)\begin{pspicture}[shift=-.9](-1.1,-1)(1.1,1)  \rput(.7,.7){\rnode{a1}{$$}}
\rput(-.7,.7){\rnode{a2}{$$}} \rput(-.7,-.7){\rnode{a3}{$$}}
\rput(.7,-.7){\rnode{a4}{$$}}
\nccurve[angleA=225,angleB=135]{a1}{a4}\middlearrow
\nccurve[angleA=315,angleB=45]{a2}{a3}\middlearrow
\end{pspicture}\end{pspicture} +
\begin{pspicture}[shift=-.9](0,1)(0,1)\begin{pspicture}[shift=-.9](-1.1,-1)(1.1,1)  \rput(.7,.7){\rnode{a1}{$$}}
\rput(-.7,.7){\rnode{a2}{$$}} \rput(-.7,-.7){\rnode{a3}{$$}}
\rput(.7,-.7){\rnode{a4}{$$}} \rput(0,.3){\rnode{b1}{$$}}
\rput(0,-.3){\rnode{b2}{$$}} \ncline{a1}{b1}\middlearrow
\ncline{a2}{b1}\middlearrow \ncline{b2}{b1}\middlearrow
\ncline{b2}{a3}\middlearrow \ncline{b2}{a4}\middlearrow
\end{pspicture}\end{pspicture} \\\nonumber \begin{pspicture}[shift=-.9](0,-1)(0,1)\end{pspicture}
\begin{pspicture}[shift=-.9](0,1)(0,1)\begin{pspicture}[shift=-.9](-1.1,-1)(1.1,1)  \rput(.7,.7){\rnode{a1}{$$}}
\rput(-.7,.7){\rnode{a2}{$$}} \rput(-.7,-.7){\rnode{a3}{$$}}
\rput(.7,-.7){\rnode{a4}{$$}} \rput(.1,.1){\rnode{b1}{$$}}
\rput(-.1,.1){\rnode{b2}{$$}} \rput(-.1,-.1){\rnode{b3}{$$}}
\rput(.1,-.1){\rnode{b4}{$$}} \ncline{a2}{b2}\middlearrow
\ncline{a1}{b1}\middlearrow \ncline{b2}{a4} \ncline{b3}{a3}
\end{pspicture}\end{pspicture} &=& q^{-\frac{1}{2}}
\begin{pspicture}[shift=-.9](0,1)(0,1)\begin{pspicture}[shift=-.9](-1.1,-1)(1.1,1)  \rput(.7,.7){\rnode{a1}{$$}}
\rput(-.7,.7){\rnode{a2}{$$}} \rput(-.7,-.7){\rnode{a3}{$$}}
\rput(.7,-.7){\rnode{a4}{$$}}
\nccurve[angleA=225,angleB=135]{a1}{a4}\middlearrow
\nccurve[angleA=315,angleB=45]{a2}{a3}\middlearrow
\end{pspicture} \end{pspicture}+
\begin{pspicture}[shift=-.9](0,1)(0,1)\begin{pspicture}[shift=-.9](-1.1,-1)(1.1,1)  \rput(.7,.7){\rnode{a1}{$$}}
\rput(-.7,.7){\rnode{a2}{$$}} \rput(-.7,-.7){\rnode{a3}{$$}}
\rput(.7,-.7){\rnode{a4}{$$}}  \rput(0,.3){\rnode{b1}{$$}}
\rput(0,-.3){\rnode{b2}{$$}} \ncline{a1}{b1}\middlearrow
\ncline{a2}{b1}\middlearrow \ncline{b2}{b1}\middlearrow
\ncline{b2}{a3}\middlearrow \ncline{b2}{a4}\middlearrow
\end{pspicture}\end{pspicture}
\end{eqnarray}
\caption{Expansion of crossings for $P_3(q)$.} \label{planar1}
\end{figure}
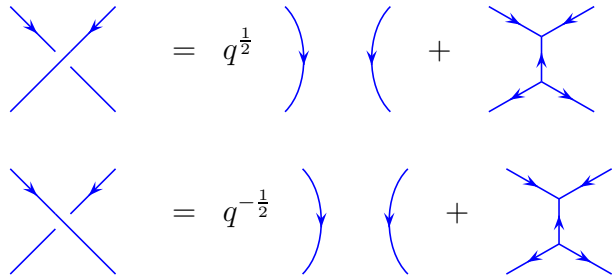

The HOMFLY polynomial $P_3(q)$ can be obtained by coloring all
components by the vector representations of the quantum
$\mathfrak{sl}(3)$ and the following skein relations

$$ P_{3}(\emptyset) =1, $$
$$ P_{3}( \begin{pspicture}[shift=-.12](-.17,-.17)(.17,.17) \pscircle(0,0){.15}
\end{pspicture} \cup D)= [3] P_{3} (D),
$$
$$
q^{\frac 32}P_{3}(L_+) - q^{-\frac 32}P_{3}(L_-) = (q^{\frac 12}-
q^{-\frac 12}) P_3(L_0),
$$
where $\emptyset$ is the empty diagram,
$\begin{pspicture}[shift=-.12](-.17,-.17)(.17,.17) \pscircle(0,0){.15}
\end{pspicture}$ is the trivial knot and $L_+, L_-$ and $L_0$ are
three diagrams which are identical except at one crossing as
illustrated in Figure~\ref{local}. On the other hand, the polynomial
$P_3(q)$ can be computed by linearly expanding each crossing into a
sum of webs as shown in Figure~\ref{planar1} then by applying
relations in Figure \ref{relations} \cite{chbili, Kuperberg:g2,
MOY:Homfly}. A benefit of using webs is that we can easily define
\emph{the colored $\mathfrak{sl}(3)$ HOMFLY polynomial} $G_3(L,\mu)$
of $L$ as follows. Let $L$ be a colored link of $l$ components say,
$L_1, L_2, \ldots, L_l$, where each component $L_i$ is colored by an
irreducible representation $V_{a_i\lambda_1 + b_i\lambda_2}$ of the
quantum $\mathfrak{sl}(3)$ and $\lambda_1, \lambda_2$ are the
fundamental weights of $\mathfrak{sl}(3)$. The coloring is denoted
by $\mu=(a_1\lambda_1 + b_1\lambda_2 , a_2\lambda_1 + b_2\lambda_2,
\ldots, a_l\lambda_1 + b_l\lambda_2)$. First we replace each
component $L_i$ by $a_i+ b_i$ copies of parallel lines and each of
$a_i$ lines is colored by the weight $\lambda_1$ and each of $b_i$
lines is colored by the weight $\lambda_2$. Then we put a clasp of
the weight $(a_i\lambda_1 + b_i\lambda_2)$ for $L_i$. If we assume
the clasps are far away from crossings, we expand each crossing as
depicted  in Figure~\ref{planar1}, then expand each clasp
inductively by Theorem~\ref{a2ab2expthm}. The value we find after
removing all faces by using the relations in Figure \ref{relations}
is \emph{the colored $\mathfrak{sl}(3)$ HOMFLY polynomial}
$G_3(L,\mu)$ of $L$. One can find the following theorem which is a
generalization of a criterion to determine the periodicity of a link
\cite{chbili, CL:period}.

\begin{thm}
Let $p$ be a positive integer and $L$ be a $p$-periodic link in
$S^3$ with the factor link $\overline{L}$. Let $\mu$ be a
$p$-periodic coloring of $L$ and $\overline{\mu}$ be the induced
coloring of $\overline{L}$. Then
$$G_3(L,\mu) \equiv G_3(\overline{L},\overline{\mu})^p \hskip ,7cm modulo \hskip .2cm
\mathcal{I}_3,$$ where $\overline L$ is the factor link and
$\mathcal{I}_3$ is the ideal of $\mathbb{Z}[q^{\pm \frac 12}]$
generated by $p$ and $[3]^p-[3]$. \label{main2}
\begin{proof}
Since the clasps are idempotents, for each component, we put $p-1$
extra clasps for each copies of components by the rotation of order
$p$. First we keep the clasps far away from the crossings. The key
idea of the proof given in \cite{chbili} is that if any expansion of
crossings occurs in the link diagram, it must be used identically
for all other $p-1$ copies of the diagram. Otherwise there will be
$p$ identical shapes by the rotation of order $p$, then it is
congruent to zero modulo $p$. By the same philosophy, if any
application of relations occurs, it must be used identically for all
other $p-1$ copies. Otherwise it is congruent to zero modulo $p$.
Once there is an unknot in the fundamental domain of the action of
order $p$, there are $p$ identical unknots by the rotation which
occurs only once in the factor link. Therefore, we get the
congruence $[3]^p-[3]$.
\end{proof}
\end{thm}

\subsection{$3j$ and $6j$ symbols for the quantum $\mathfrak{sl}(3)$
representation theory}

\begin{figure}\begin{align*}
\begin{pspicture}[shift=-1.6](-.95,-1.7)(.95,1.7) \pscircle*(-.8,0){.1}
\pscircle*(.8,0){.1}
\pccurve[angleA=60,angleB=120,ncurv=1](-.8,0)(.8,0)
\pccurve[angleA=0,angleB=180,ncurv=1](-.8,0)(.8,0)
\pccurve[angleA=-60,angleB=-120,ncurv=1](-.8,0)(.8,0)
\rput(0,1){$a$}\rput(0,.25){$b$}\rput(0,-.89){$c$}
\end{pspicture} =
\begin{pspicture}[shift=-1.6](-1.75,-1.7)(1.75,1.7) \pcline(-1.3,-.9)(-.3,-.1)
\pcline(1.3,.9)(.3,.1) \pcline(1.3,-.9)(.3,-.1)
\pcline(-1.3,.9)(-.3,.1)
\pccurve[angleA=225,angleB=315,ncurv=1](.7,.9)(-.7,.9)
\pccurve[angleA=45,angleB=135,ncurv=1](-.7,-.9)(.7,-.9)\rput(0,.8){$j$}
\rput(0,-.85){$j$} \qline(-1,-1.5)(-1,-1.1)\qline(-1,1.1)(-1,1.5)
\qline(1,-1.5)(1,-1.1) \qline(1,1.1)(1,1.5) \qline(1,-1.5)(1.7,-1.5)
\qline(1.7,-1.5)(1.7,1.5) \qline(1.7,1.5)(1,1.5)
\qline(-1,1.5)(-1.7,1.5) \qline(-1.7,1.5)(-1.7,-1.5)
\qline(-1.7,-1.5)(-1,-1.5) \rput[b](-1,.2){$i$}\rput[b](1,.2){$k$}
\psframe[linecolor=black](-1.5,-1.1)(-.5,-.9)
\psframe[linecolor=black]( 1.5,-1.1)( .5,-.9)
\psframe[linecolor=black](-1.5, 1.1)(-.5, .9)
\psframe[linecolor=black]( 1.5, 1.1)( .5, .9)
\psframe[linecolor=black](-.5,.1)(.5,-.1)
\end{pspicture}
=(-1)^{i+j+k}\frac{[i+j+k+1]![i]![j]![k]!}{[i+j]![j+k]![i+k]!}.
\end{align*}
\caption{Trihedron coefficients for $\mathfrak{sl}(2)$.}
\label{a1triexp}
\end{figure}

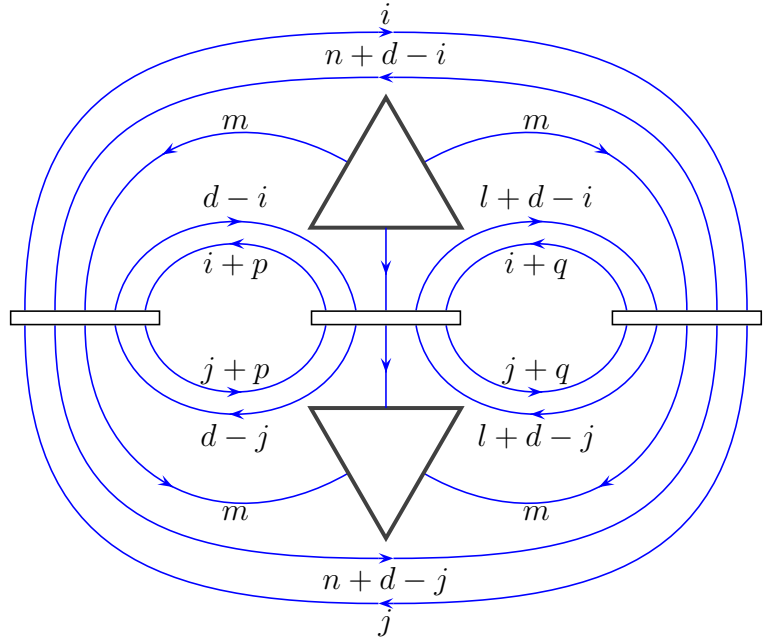
\begin{figure} $$
\begin{pspicture}[shift=-4.1](-5.1,-4.2)(5.1,4.2) \rput(0,4.05){$i$}
\rput(-2,1.6){$d-i$} \rput(0,3.5){$n+d-i$} \rput(-2,2.6){$m$}
\rput(2,2.6){$m$} \rput(-2,.7){$i+p$} \rput(2,1.6){$l+d-i$}
\rput(2,.7){$i+q$}  \rput(0,-4.05){$j$} \rput(-2,-1.6){$d-j$}
\rput(0,-3.5){$n+d-j$} \rput(-2,-2.6){$m$} \rput(2,-2.6){$m$}
\rput(-2,-.7){$j+p$} \rput(2,-1.6){$l+d-j$}
\rput(2,-.7){$j+q$}\psline[arrowscale=1.5]{->}(-.1,3.8)(.1,3.8)
\psline[arrowscale=1.5]{->}(.1,3.2)(-.1,3.2)
\pccurve[angleA=90,angleB=180,ncurv=1](-4.8,.1)(-.1,3.8)
\pccurve[angleA=90,angleB=180,ncurv=1](-4.4,.1)(-.1,3.2)
\pccurve[angleA=150,angleB=90,ncurv=1](-.5,2.066)(-4,.1)\middlearrow
\pccurve[angleA=100,angleB=80,ncurv=1](-.8,.1)(-3.2,.1)\middlearrow
\pccurve[angleA=80,angleB=100,ncurv=1](-3.6,.1)(-.4,.1)\middlearrow
\pccurve[angleA=80,angleB=100,ncurv=1](.4,.1)(3.6,.1)\middlearrow
\pccurve[angleA=100,angleB=80,ncurv=1](3.2,.1)(.8,.1)\middlearrow
\pccurve[angleA=30,angleB=90,ncurv=1](.5,2.066)(4,.1)\middlearrow
\pccurve[angleA=0,angleB=90,ncurv=1](.1,3.2)(4.4,.1)
\pccurve[angleA=0,angleB=90,ncurv=1](.1,3.8)(4.8,.1)
\pspolygon[linecolor=darkgray,linewidth=1.5pt](1,1.2)(-1,1.2)(0,2.932)
\psline[arrowscale=1.5]{->}(.1,-3.8)(-.1,-3.8)
\pcline(0,1.2)(0,.1)\middlearrow
\pccurve[angleA=-90,angleB=-180,ncurv=1](-4.8,-.1)(-.1,-3.8)
\pccurve[angleA=-90,angleB=-180,ncurv=1](-4.4,-.1)(-.1,-3.2)
\pccurve[angleA=-90,angleB=-150,ncurv=1](-4,-.1)(-.5,-2.066)\middlearrow
\pccurve[angleA=-80,angleB=-100,ncurv=1](-3.2,-.1)(-.8,-.1)\middlearrow
\pccurve[angleA=-100,angleB=-80,ncurv=1](-.4,-.1)(-3.6,-.1)\middlearrow
\pccurve[angleA=-100,angleB=-80,ncurv=1](3.6,-.1)(.4,-.1)\middlearrow
\pccurve[angleA=-80,angleB=-100,ncurv=1](.8,-.1)(3.2,-.1)\middlearrow
\pccurve[angleA=-90,angleB=-30,ncurv=1](4,-.1)(.5,-2.066)\middlearrow
\pccurve[angleA=0,angleB=-90,ncurv=1](.1,-3.2)(4.4,-.1)
\pccurve[angleA=0,angleB=-90,ncurv=1](.1,-3.8)(4.8,-.1)
\pspolygon[linecolor=darkgray,linewidth=1.5pt](1,-1.2)(-1,-1.2)(0,-2.932)
\psline[arrowscale=1.5]{->}(-.1,-3.2)(.1,-3.2)
\pcline(0,-.1)(0,-1.2)\middlearrow
\psframe[linecolor=black](-1,-.1)(1,.1)
\psframe[linecolor=black](-5,-.1)(-3,.1)
\psframe[linecolor=black](3,-.1)(5,.1)
\end{pspicture} $$
\caption{General shape of $\Theta
(a_1,b_1,a_2$,$b_2$,$a_3$,$b_3$$;i$,$j)$} \label{thetadef}
\end{figure}

$3j$ symbols and $6j$ symbols for the quantum $\mathfrak{sl}(2)$
representation theory have many significant implications in
Mathematics and Physics. $3j$ symbols are given in the equation
shown in Figure \ref{a1triexp} \cite{MV:3valent}. Its natural
generalization for the quantum $\mathfrak{sl}(3)$ representation
theory was first suggested \cite{Kuperberg:spiders} and studied
\cite{kim:theta}. Let $\lambda_1, \lambda_2$ be the fundamental
dominant weights of $\mathfrak{sl}$$(3,\mathbb{C})$. Let
$V_{a\lambda_1 +b\lambda_2}$ be an irreducible representation of
$\mathfrak{sl}$$(3,\mathbb{C})$ of highest weight $a\lambda_1
+b\lambda_2$. Now each edge of $\Theta$ is decorated by an
irreducible representation of $\mathfrak{sl}(3)$, let say
$V_{a_1\lambda_1+b_1\lambda_2}$, $V_{a_2\lambda_1+b_2\lambda_2}$ and
$V_{a_3\lambda_1+b_3\lambda_2}$ where $a_i, b_j$ are nonnegative
integers. Let $d=$ min $\{a_1,a_2,a_3,b_1,b_2,b_3\}$. If
$\dim(\Inv(V_{a_1\lambda_1+b_1\lambda_2}\otimes
V_{a_2\lambda_1+b_2\lambda_2}\otimes
V_{a_3\lambda_1+b_3\lambda_2}))$ is nonzero, in fact $d+1$, then we
say a triple of ordered pairs $((a_1, b_1)$, $(a_2, b_2)$, $(a_3,
b_3))$ is \emph{admissible}. One can show $((a_1, b_1)$, $(a_2,
b_2)$, $(a_3, b_3))$ is admissible if and only if there exist
nonnegative integers $k$, $l$, $m$, $n$, $o$, $p$, $q$ such that
$a_2=d+l+p$, $a_3=d+n+q$, $b_1=d+k+p$, $b_2=d+m+q$, $b_3=d+o$ and
$k-n=o-l=m$. For an admissible triple, we can write its trihedron
coefficients as a $(d+1)\times (d+1)$ matrix. Let us denote it by
$M_{\Theta}$ $(a_1$, $b_1$, $a_2$, $b_2$, $a_3$, $b_3)$ or
$M_{\Theta}(\lambda)$ where $\lambda = $ $a_1\lambda_1 $ $+
b_1\lambda_2$ $+ a_2\lambda_1 $ $+ b_2\lambda_2 $ $+ a_3\lambda_1 $
$+ b_3 \lambda_2$. Also we denotes its $(i,j)$ entry by
$\Theta(a_1$, $b_1$, $a_2$, $b_2$, $a_3$, $b_3$ $;$ $i$, $j)$ or
$\Theta(\lambda;i,j)$ where $0\le i,j \le d$. The trihedron shape of
$\Theta(a_1$, $b_1$, $a_2$, $b_2$, $a_3$, $b_3$ $;$ $i$, $j)$ is
given in Figure~\ref{thetadef} where the triangles are filled by cut
outs from the hexagonal tiling of the
plane~\cite{Kuperberg:spiders}. $M_{\Theta}(0$, $m+n$, $l$, $m+q$,
$n+q$, $m+l)$, $M_{\Theta}$ $(0$, $n+p$, $p+l$, $q$, $n+q$, $l)$ and
$M_{\Theta}$ $(i$, $j+k$, $k+l$, $m$, $j+m$, $j+l$ $;0$, $0)$ were
found in \cite{kim:theta}. All other cases of $3j$ symbols and $6j$
symbols are left open.

\subsection{$\mathfrak{sl}(3)$ invariants of cubic planar bipartite
graphs}

The $\mathfrak{sl}(3)$ webs are directed cubic bipartite planar
graphs together circles (no vertices) where the direction of the
edges is from one set to the other set in the bipartition. From a
given directed cubic bipartite planar graph, we remove all circles
by the relation \ref{a2defr21} and then remove the multiple edges by
the relation \ref{a2defr22} in Figure~\ref{relations}. Using a
simple application of the Euler characteristic number of a graphs in
the unit disc, we can show the existence of a rectangular
face~\cite{OY:quantum}. By inducting on the number of faces, we
provides the existence of the quantum $\mathfrak{sl}(3)$ invariants
of directed cubic bipartite planar graphs. It is fairly easy to
prove the quantum $\mathfrak{sl}(3)$ invariant does not depend on
the choice of directions in the bipartition. Thus, the quantum
$\mathfrak{sl}(3)$ invariant naturally extends to any cubic
bipartite planar graph $G$, let us denote it by $P_G(q)$. By using a
favor of graph theory, we find a classification theorem and  provide
a method to find all 3-connected cubic bipartite planar graphs which
is called \emph{prime webs} \cite{KL:graph}. As little as it is
known about the properties of the quantum invariants of links, we
know a very little how $P_G(q)$ tells us about the properties of
graphs.

For symmetries of cubic bipartite planar graph, the idea of the
Theorem~\ref{main2} and \ref{main3} works for the $\mathfrak{sl}(3)$
graph invariants with one exception. There is a critical difference
between these two invariants which is illustrated in
Theorem~\ref{symthm}.

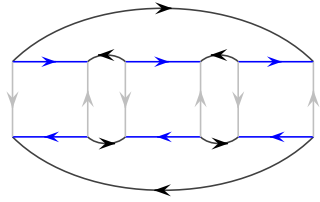
\begin{figure}
$$
\begin{pspicture}[shift=-1.1](-2.2,-.3)(2.2,2.1) \rput(2,1.5){\rnode{a1}{$$}}
\rput(1,1.5){\rnode{a2}{$$}} \rput(1,.5){\rnode{a3}{$$}}
\rput(2,.5){\rnode{a4}{$$}} \rput(.5,1.5){\rnode{b1}{$$}}
\rput(-.5,1.5){\rnode{b2}{$$}} \rput(-.5,.5){\rnode{b3}{$$}}
\rput(.5,.5){\rnode{b4}{$$}} \rput(-2,1.5){\rnode{d2}{$$}}
\rput(-1,1.5){\rnode{d1}{$$}} \rput(-1,.5){\rnode{d4}{$$}}
\rput(-2,.5){\rnode{d3}{$$}} \rput(-.5,-1.5){\rnode{c1}{$$}}
\rput(.5,-1.5){\rnode{c2}{$$}}
\nccurve[angleA=45,angleB=135,linecolor=darkgray]{d2}{a1}\mda
\nccurve[angleA=135,angleB=45,linecolor=darkgray]{a2}{b1}\mda
\nccurve[angleA=135,angleB=45,linecolor=darkgray]{b2}{d1}\mda
\nccurve[angleA=-45,angleB=-135,linecolor=darkgray]{b4}{a3}\mda
\nccurve[angleA=-45,angleB=-135,linecolor=darkgray]{d4}{b3}\mda
\nccurve[angleA=-135,angleB=-45,linecolor=darkgray]{a4}{d3}\mda
\ncline{a2}{a1}\middlearrow \ncline{a4}{a3}\middlearrow
\ncline{b2}{b1}\middlearrow \ncline{b4}{b3}\middlearrow
\ncline{d2}{d1}\middlearrow \ncline{d4}{d3}\middlearrow
\ncline[linecolor=lightgray]{b2}{b3}\mdb
\ncline[linecolor=lightgray]{b4}{b1}\mdb
\ncline[linecolor=lightgray]{a2}{a3}\mdb
\ncline[linecolor=lightgray]{a4}{a1}\mdb
\ncline[linecolor=lightgray]{d2}{d3}\mdb
\ncline[linecolor=lightgray]{d4}{d1}\mdb
\end{pspicture}
$$
\caption {Prime web $6_1$.} \label{a2graph61}
\end{figure}

\begin{thm} [\cite{KL:graph}]
Let $G$ be a planar cubic bipartite graph with the group of
symmetries $\Gamma$ of order $n$. Let $\Gamma_d$ be a subgroup of
$\Gamma$ of order $d$ such that the fundamental domain of
$G/\Gamma_d$ is not a basis web with the given boundary. Then
$$P_{G}(q) \equiv (P_{G/\Gamma_d}(q))^d \hskip ,7cm modulo \hskip .2cm
\mathcal{I}_d,$$ where $\mathcal{I}_d$ is the ideal of
$\mathbb{Z}[q^{\pm \frac 12}]$ generated by $d$ and $[3]^d-[3]$.
\label{symthm}
\end{thm}

If the fundamental domain of $G/\Gamma$ is a basis web with the
given boundary, then the main idea of the theorem no longer works
and a counterexample was found as follows \cite{KL:graph}. We look
at an example $6_1$ as shown in Figure~\ref{a2graph61}. By a help of
a machine, we can see that there does not exist an $\alpha \in
\mathbb{Z}[q^{\pm \frac 12}]$ such that
$$(\alpha)^6 \equiv [2]^4[3]+2[2]^2[3]
~~\mathrm{mod}~~\mathcal{I}_6$$ even though there do exist a
symmetry of order $6$ for $6_1$.

\subsection{Applications for the quantum
$\mathfrak{sp}(4)$ representation theory}

A quantum $\mathfrak{sp}(4)$ polynomial invariant
$G_{\mathfrak{sp}(4)}(L,\mu)$ can be defined \cite{Kuperberg:g2,
Kuperberg:spiders} where $\mu$ is a fundamental representation of
the quantum $\mathfrak{sp}(4)$. Since we have found single clasp
expansion of the clasps of weight $(a,0)$ and $(0,b)$, we can extend
$G_{\mathfrak{sp}(4)}(L,\mu)$ for $\mu$ is an irreducible
representations of weight either $(a,0)$ and $(0,b)$. If we assume a
coloring $\mu=(a,0)$ or $\mu=(0,b)$, by the same idea of the proof
of Theorem \ref{main2}, we can find the following theorem from
Corollary \ref{b2n0expcor} and \ref{b20nexpcor}.

\begin{thm}
Let $p$ be a positive integer and $L$ be a $p-$periodic link in
$S^3$ with the factor link $\overline{L}$. Let $\mu$ be a
$p$-periodic coloring of $L$ and $\overline{\mu}$ be the induced
coloring of $\overline{L}$. Then
$$G_{\mathfrak{sp}(4)}(L,\mu) \equiv G_{\mathfrak{sp}(4)}(\overline{L},\overline{\mu})^p \hskip ,7cm modulo \hskip .2cm
\mathcal{I}_{\mathfrak{sp}(4)},$$ where $\overline L$ is the factor
link and $\mathcal{I}_{\mathfrak{sp}(4)}$ is the ideal of
$\mathbb{Z}[q^{\pm \frac 12}]$ generated by $p$,
$(-\frac{[6][2]}{[3]})^p+\frac{[6][2]}{[3]}$ and
$(\frac{[6][5]}{[3][2]})^p-\frac{[6][5]}{[3][2]}$. \label{main3}
\end{thm}

In fact, Theorem \ref{main3} remains true even if $\mu$ is any
finite dimensional irreducible representation of $\mathfrak{sp}(4)$,
but we would not be able to obtain the actual polynomials because
any expansion is not known for the clasp of the weight $(a,b)$ where
$a\neq 0\neq b$.

\section{The Proof of Lemmas} \label{lemmas}

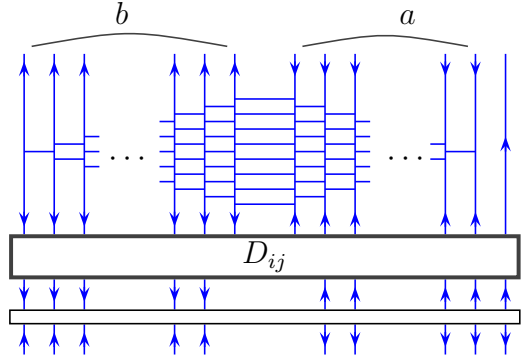
\begin{figure}
$$
\begin{pspicture}[shift=-2.4](-.5,-1.2)(7.5,3.8)
\pccurve[angleA=10,angleB=170,ncurv=1,linecolor=darkgray](.2,3.1)(2.8,3.1)
\pccurve[angleA=10,angleB=170,ncurv=1,linecolor=darkgray](3.8,3.1)(6,3.1)
\rput[b](3.3,0.1){$D_{ij}$} \rput[b](1.4,3.4){$b$}
\rput[b](5.2,3.4){$a$}
\psline(.1,.6)(.1,3)\psline[arrowscale=1.5]{->}(.1,.9)(.1,.7)
\psline[arrowscale=1.5]{->}(.1,2.7)(.1,2.9)
\psline(.1,0)(.1,-.4)\psline[arrowscale=1.5]{->}(.1,-.1)(.1,-.3)
\psline(.1,-.6)(.1,-1)\psline[arrowscale=1.5]{->}(.1,-.9)(.1,-.7)
\psline(.5,.6)(.5,3)\psline[arrowscale=1.5]{->}(.5,.9)(.5,.7)
\psline[arrowscale=1.5]{->}(.5,2.7)(.5,2.9)
\psline(.5,0)(.5,-.4)\psline[arrowscale=1.5]{->}(.5,-.1)(.5,-.3)
\psline(.5,-.6)(.5,-1)\psline[arrowscale=1.5]{->}(.5,-.9)(.5,-.7)
\psline(.9,.6)(.9,3)\psline[arrowscale=1.5]{->}(.9,.9)(.9,.7)
\psline[arrowscale=1.5]{->}(.9,2.7)(.9,2.9)
\psline(.9,0)(.9,-.4)\psline[arrowscale=1.5]{->}(.9,-.1)(.9,-.3)
\psline(.9,-.6)(.9,-1)\psline[arrowscale=1.5]{->}(.9,-.9)(.9,-.7)
\psline(2.1,.6)(2.1,3)\psline[arrowscale=1.5]{->}(2.1,.9)(2.1,.7)
\psline[arrowscale=1.5]{->}(2.1,2.7)(2.1,2.9)
\psline(2.1,0)(2.1,-.4)\psline[arrowscale=1.5]{->}(2.1,-.1)(2.1,-.3)
\psline(2.1,-.6)(2.1,-1)\psline[arrowscale=1.5]{->}(2.1,-.9)(2.1,-.7)
\psline(2.5,.6)(2.5,3)\psline[arrowscale=1.5]{->}(2.5,.9)(2.5,.7)
\psline[arrowscale=1.5]{->}(2.5,2.7)(2.5,2.9)
\psline(2.5,0)(2.5,-.4)\psline[arrowscale=1.5]{->}(2.5,-.1)(2.5,-.3)
\psline(2.5,-.6)(2.5,-1)\psline[arrowscale=1.5]{->}(2.5,-.9)(2.5,-.7)
\psline(2.9,.6)(2.9,3)\psline[arrowscale=1.5]{->}(2.9,.9)(2.9,.7)
\psline[arrowscale=1.5]{->}(2.9,2.7)(2.9,2.9)
\psline(3.7,.6)(3.7,3)\psline[arrowscale=1.5]{<-}(3.7,.9)(3.7,.7)
\psline[arrowscale=1.5]{<-}(3.7,2.7)(3.7,2.9)
\psline(4.1,.6)(4.1,3)\psline[arrowscale=1.5]{<-}(4.1,.9)(4.1,.7)
\psline[arrowscale=1.5]{<-}(4.1,2.7)(4.1,2.9)
\psline(4.1,0)(4.1,-.4)\psline[arrowscale=1.5]{<-}(4.1,-.1)(4.1,-.3)
\psline(4.1,-.6)(4.1,-1)\psline[arrowscale=1.5]{<-}(4.1,-.9)(4.1,-.7)
\psline(4.5,.6)(4.5,3)\psline[arrowscale=1.5]{<-}(4.5,.9)(4.5,.7)
\psline[arrowscale=1.5]{<-}(4.5,2.7)(4.5,2.9)
\psline(4.5,0)(4.5,-.4)\psline[arrowscale=1.5]{<-}(4.5,-.1)(4.5,-.3)
\psline(4.5,-.6)(4.5,-1)\psline[arrowscale=1.5]{<-}(4.5,-.9)(4.5,-.7)
\psline(5.7,.6)(5.7,3)\psline[arrowscale=1.5]{<-}(5.7,.9)(5.7,.7)
\psline[arrowscale=1.5]{<-}(5.7,2.7)(5.7,2.9)
\psline(5.7,0)(5.7,-.4)\psline[arrowscale=1.5]{<-}(5.7,-.1)(5.7,-.3)
\psline(5.7,-.6)(5.7,-1)\psline[arrowscale=1.5]{<-}(5.7,-.9)(5.7,-.7)
\psline(6.1,.6)(6.1,3)\psline[arrowscale=1.5]{<-}(6.1,.9)(6.1,.7)
\psline[arrowscale=1.5]{<-}(6.1,2.7)(6.1,2.9)
\psline(6.1,0)(6.1,-.4)\psline[arrowscale=1.5]{<-}(6.1,-.1)(6.1,-.3)
\psline(6.1,-.6)(6.1,-1)\psline[arrowscale=1.5]{<-}(6.1,-.9)(6.1,-.7)
\psline(6.5,.6)(6.5,3)\psline[arrowscale=1.5]{->}(6.5,1.7)(6.5,1.9)
\psline(6.5,0)(6.5,-.4)\psline[arrowscale=1.5]{<-}(6.5,-.1)(6.5,-.3)
\psline(6.5,-.6)(6.5,-1)\psline[arrowscale=1.5]{<-}(6.5,-.9)(6.5,-.7)
\qline(.1,1.7)(.5,1.7) \qline(.5,1.8)(.9,1.8) \qline(.5,1.6)(.9,1.6)
\qline(.9,1.9)(1.1,1.9) \qline(.9,1.7)(1.1,1.7)
\qline(.9,1.5)(1.1,1.5) \qline(1.9,2.1)(2.1,2.1)
\qline(1.9,1.9)(2.1,1.9) \qline(1.9,1.7)(2.1,1.7)
\qline(1.9,1.5)(2.1,1.5) \qline(1.9,1.3)(2.1,1.3)
\qline(2.1,2.2)(2.5,2.2) \qline(2.1,2)(2.5,2)
\qline(2.1,1.8)(2.5,1.8) \qline(2.1,1.6)(2.5,1.6)
\qline(2.1,1.4)(2.5,1.4) \qline(2.1,1.2)(2.5,1.2)
\qline(2.5,2.3)(2.9,2.3) \qline(2.5,2.1)(2.9,2.1)
\qline(2.5,1.9)(2.9,1.9) \qline(2.5,1.7)(2.9,1.7)
\qline(2.5,1.5)(2.9,1.5) \qline(2.5,1.3)(2.9,1.3)
\qline(2.5,1.1)(2.9,1.1) \qline(4.1,2.3)(3.7,2.3)
\qline(4.1,2.1)(3.7,2.1) \qline(4.1,1.9)(3.7,1.9)
\qline(4.1,1.7)(3.7,1.7) \qline(4.1,1.5)(3.7,1.5)
\qline(4.1,1.3)(3.7,1.3) \qline(4.1,1.1)(3.7,1.1)
\qline(4.1,2.2)(4.5,2.2) \qline(4.1,2)(4.5,2)
\qline(4.1,1.8)(4.5,1.8) \qline(4.1,1.6)(4.5,1.6)
\qline(4.1,1.4)(4.5,1.4) \qline(4.1,1.2)(4.5,1.2)
\qline(4.5,2.1)(4.7,2.1) \qline(4.5,1.9)(4.7,1.9)
\qline(4.5,1.7)(4.7,1.7) \qline(4.5,1.5)(4.7,1.5)
\qline(4.5,1.3)(4.7,1.3) \qline(5.7,1.7)(6.1,1.7)
\qline(5.5,1.8)(5.7,1.8) \qline(5.5,1.6)(5.7,1.6)
\qline(2.9,1)(3.7,1) \qline(2.9,1.2)(3.7,1.2)
\qline(2.9,1.4)(3.7,1.4) \qline(2.9,1.6)(3.7,1.6)
\qline(2.9,1.8)(3.7,1.8) \qline(2.9,2)(3.7,2)
\qline(2.9,2.2)(3.7,2.2)
\qline(2.9,2.4)(3.7,2.4)\rput[t](1.5,1.7){$\cdots$}\rput[t](5.2,1.7){$\cdots$}
\psframe[linecolor=black](-0.1,-.6)(6.7,-0.4)
\psframe[linecolor=darkgray,linewidth=1.5pt](-0.1,0)(6.7,0.6)
\end{pspicture}
$$
\caption {A sequence of H's which transforms $D_{ij}$ to a linear
combinations of webs in the single clasp expansion of segregated
clasp of weight $(a,b-1)$.} \label{a2abexp4}
\end{figure}

Let us recalled that the relation~\ref{a2defr23} in Figure
\ref{relations} is called a \emph{rectangular relation} and the
first(second) web in the right-hand side of the equality is called a
\emph{horizontal(vertical}, respectively) \emph{splitting}. The web
in the equation shown in Figure~\ref{a2abexp} corresponding to the
coefficient $a_{i,j}$ is denoted by $D_{i,j}$. After attaching $H$'s
to $D_{i,j}$ as illustrated in Figure~\ref{a2abexp4}, the resulting
web is denoted by $\tilde D_{i,j}$. We find that $\tilde D_{i,j}$
contains some elliptic faces. If we decompose each $\tilde D_{i,j}$
into a linear combination of webs which have no elliptic faces, then
the union of all these webs forms a basis. Let us prove that these
webs actually form a basis which will be denoted by $D'_{i',j'}$. As
vector spaces, this change, adding $H$'s as in Figure
\ref{a2abexp4}, induces an isomorphism between two web spaces
because its matrix representation with respect to these web bases
$\{D_{i,j}\}$ and $\{ D'_{i',j'}\}$ is an $(a+1)b\times (a+1)b$
matrix whose determinant is $\pm [2]^{ab}$ because a single $H$
contributes $\pm [2]$ depends on the choice of the direction of $H$.

\begin{figure}
$$
\begin{pspicture}[shift=-2.3](-.5,-1.2)(4,3.5) \rput[t](.8,1.8){$\cdots$}
\rput[t](2.6,1.8){$\cdots$} \pcline(.1,-0.4)(.1,3)\middlearrow
\pcline(.3,-0.4)(.3,3)\middlearrow
\pcline(1.1,-0.4)(1.1,3)\middlearrow
\pcline(1.3,-0.4)(1.3,3)\middlearrow
\pcline(1.5,-0.4)(1.5,3)\middlearrow
\pcline(1.9,3)(1.9,-0.4)\middlearrow
\pcline(2.1,3)(2.1,-0.4)\middlearrow
\pcline(2.9,3)(2.9,-0.4)\middlearrow
\pcline(3.1,3)(3.1,0.4)\middlearrow
\pcline(3.3,0.4)(3.3,3)\middlearrow
\pccurve[angleA=-90,angleB=-90,ncurv=1](3.1,0.4)(3.3,0.4)
\pcline(0.1,-1)(0.1,-.6)\middlearrow
\pcline(0.3,-1)(0.3,-.6)\middlearrow
\pcline(0.5,-1)(0.5,-.6)\middlearrow
\pcline(1.1,-1)(1.1,-.6)\middlearrow
\pcline(1.3,-1)(1.3,-.6)\middlearrow
\pcline(1.5,-1)(1.5,-.6)\middlearrow
\pcline(1.7,-.6)(1.7,-1)\middlearrow
\pcline(2.1,-.6)(2.1,-1)\middlearrow
\pcline(2.3,-.6)(2.3,-1)\middlearrow
\pcline(2.9,-.6)(2.9,-1)\middlearrow
\psframe[linecolor=black](0,3)(3,3.2)
\psframe[linecolor=black](0,-.6)(3.4,-0.4)
\psframe[linecolor=darkgray,fillstyle=solid,fillcolor=white,linewidth=1.5pt](0,.2)(3,0.6)
\end{pspicture}
, \begin{pspicture}[shift=-2.3](-.5,-1.2)(4,3.5) \rput[t](.8,1.8){$\cdots$}
\rput[t](2.4,1.8){$\cdots$} \pcline(.1,-0.4)(.1,3)\middlearrow
\pcline(.3,-0.4)(.3,3)\middlearrow
\pcline(1.1,-0.4)(1.1,3)\middlearrow
\pcline(1.3,-0.4)(1.3,3)\middlearrow
\pcline(1.5,-0.4)(1.5,3)\middlearrow
\pcline(1.9,3)(1.9,-0.4)\middlearrow
\pcline(2.1,3)(2.1,-0.4)\middlearrow
\pcline(2.7,3)(2.7,-0.4)\middlearrow
\pcline(3.3,-0.4)(3.3,3)\middlearrow
\pcline(0.1,-1)(0.1,-.6)\middlearrow
\pcline(0.3,-1)(0.3,-.6)\middlearrow
\pcline(0.5,-1)(0.5,-.6)\middlearrow
\pcline(1.1,-1)(1.1,-.6)\middlearrow
\pcline(1.3,-1)(1.3,-.6)\middlearrow
\pcline(1.5,-1)(1.5,-.6)\middlearrow
\pcline(1.7,-.6)(1.7,-1)\middlearrow
\pcline(2.1,-.6)(2.1,-1)\middlearrow
\pcline(2.3,-.6)(2.3,-1)\middlearrow
\pcline(2.9,-.6)(2.9,-1)\middlearrow \qline(2.9,2)(3.1,2)
\psline[arrows=->,arrowscale=1.5](2.9,3)(2.9,2.5)\psline(2.9,2.5)(2.9,2)
\psline[arrows=->,arrowscale=1.5](3.1,3)(3.1,2.5)\psline(3.1,2.5)(3.1,2)
\psline[arrows=->,arrowscale=1.5](3,1)(3,1.5)\psline(3,1.5)(3,2)
\qline(3,1)(3.3,1) \psframe[linecolor=black](0,-.6)(3.4,-0.4)
\psframe[linecolor=black](0,3)(3,3.2)
\psframe[linecolor=darkgray,fillstyle=solid,fillcolor=white,linewidth=1.5pt](0,.2)(2.8,0.6)
\end{pspicture} ,
\begin{pspicture}[shift=-2.3](-.5,-1.2)(4,3.5) \rput[t](.8,1.8){$\cdots$}
\rput[t](2.4,1.8){$\cdots$} \pcline(.1,-0.4)(.1,3)\middlearrow
\pcline(.3,-0.4)(.3,3)\middlearrow
\pcline(1.1,-0.4)(1.1,3)\middlearrow
\pcline(1.3,-0.4)(1.3,3)\middlearrow
\pcline(1.5,-0.4)(1.5,3)\middlearrow
\pcline(1.9,3)(1.9,-0.4)\middlearrow
\pcline(2.1,3)(2.1,-0.4)\middlearrow
\pcline(2.7,3)(2.7,-0.4)\middlearrow
\pcline(3.3,-0.4)(3.3,3)\middlearrow
\pcline(0.1,-1)(0.1,-.6)\middlearrow
\pcline(0.3,-1)(0.3,-.6)\middlearrow
\pcline(0.5,-1)(0.5,-.6)\middlearrow
\pcline(1.1,-1)(1.1,-.6)\middlearrow
\pcline(1.3,-1)(1.3,-.6)\middlearrow
\pcline(1.5,-1)(1.5,-.6)\middlearrow
\pcline(1.7,-.6)(1.7,-1)\middlearrow
\pcline(2.1,-.6)(2.1,-1)\middlearrow
\pcline(2.3,-.6)(2.3,-1)\middlearrow
\pcline(2.9,-.6)(2.9,-1)\middlearrow \qline(2.9,2)(3.1,2)
\psline[arrows=->,arrowscale=1.5](2.9,3)(2.9,2.5)\psline(2.9,2.5)(2.9,2)
\psline[arrows=->,arrowscale=1.5](3.1,3)(3.1,2.5)\psline(3.1,2.5)(3.1,2)
\psline[arrows=->,arrowscale=1.5](3,.8)(3,1.5)\psline(3,1.5)(3,2)
\qline(3,.8)(2.7,.8) \psframe[linecolor=black](0,-.6)(3.4,-0.4)
\psframe[linecolor=black](0,3)(3,3.2)
\psframe[linecolor=darkgray,fillstyle=solid,fillcolor=white,linewidth=1.5pt](0,.2)(2.8,1)
\end{pspicture}$$
 \caption{Three webs which do not vanish after attaching the clasp of weight $(a,b-1)$ to the top left side
of webs $\tilde D_{i,j}$ from the equation in
Figure~\ref{a2abexp4}.} \label{a2ab2exppo}
\end{figure}
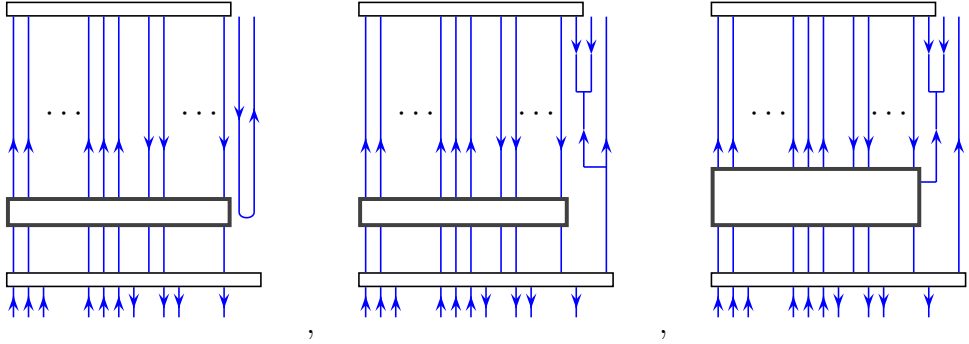

To find a single clasp expansion of the segregated clasp of weight
$(a,b)$, we have to find all linear expansions of $\tilde D_{i,j}$
into a new web basis $D'_{i',j'}$. In general this is very
complicate. Instead of using relations for linear expansions, we
look for an alternative. From $\tilde D_{i,j}$ we see that there are
$a+b+1$ nodes on top and $a+b-1$ nodes right above the clasp. A $Y$
shape in the web $D_{i,j}$ forces $\tilde D_{i,j}$ to have at least
one rectangular face. Each splitting creates another rectangular
face until it becomes a basis web (possibly using the relation
\ref{a2defr22} in Figure \ref{relations} once). If we repeatedly use
the rectangular relations as in equation~\ref{a2defr23} in Figure
\ref{relations}, we can push up $Y$'s so that there are either two
$Y$'s or one $U$ shape at the top. A \emph{stem} of a web is $a+b-1$
disjoint union of vertical lines which connect top $a+b-1$ nodes out
of $a+b+1$ nodes to the clasp of weight $(a,b-1)$ together a
$U$-turn or two $Y$'s on top. It is clear that these connecting
lines should be mutually disjoint, otherwise, we will have a cut
path with weight less than $(a,b-1)$, $i. e.,$ the web is zero.
Unfortunately some of stems do not arise naturally in the linear
expansion of $\tilde D_{i,j}$ because it may not be obtained by
removing elliptic faces. If a stem appears, we call it an \emph{
admissible stem}. For single clasp expansions, finding all these
admissible stems will be more difficult than linear expansions by
relations. But for double clasp expansions of segregated clasps of
weight $(a,b)$ there are only few possible admissible stems whose
coefficients are nonzero.

\begin{lem}
After attaching the clasp of weight $(a,b-1)$ to the top left side
of webs $\tilde D_{i,j}$ from the equation in Figure~\ref{a2abexp4},
the only non-vanishing webs are those three webs as depicted in
Figure~\ref{a2ab2exppo}. \label{a2abexplem2}
\end{lem}

\begin{proof}
It is possible to have two adjacent $Y$'s which appear in the second
and third webs in Figure~\ref{a2ab2exppo} but a $U$-turn can appear
in only two places because of the orientation of edges. If we attach
the clasp of weight $(a,b-1)$ to the northwest corner of the
resulting web and if there is a $U$ or a $Y$ shape just below the
clasp of weight $(a,b-1)$, the web becomes zero. Therefore only
these three webs do not vanish.
\end{proof}

In the following lemma, we find all $\tilde D_{i,j}$'s which can be
transformed to each of the web in Figure~\ref{a2ab2exppo}.

\begin{lem}
Only $\tilde D_{1,a}(\tilde D_{2,a})$ can be transformed to the
first(second, respectively) web in Figure~\ref{a2ab2exppo}. Only the
three webs, $\tilde D_{1,a-1}, \tilde D_{1,a}$ and $\tilde
D_{2,a-1}$ can be transformed to the last web. Moreover, all of
these transformations use only rectangular relations as in
equation~\ref{a2defr23} except the transformation from $\tilde
D_{1,a-1}$ to the third web uses the relation~\ref{a2defr22} in
Figure \ref{relations} exactly once. \label{a2abexplem3}
\end{lem}
\begin{proof}
For the first web shown in Figure~\ref{a2ab2exppo}, it is fairly
easy to see that we need to look at $\tilde D_{i,a}$, for $i=1, 2,
\ldots, b$, otherwise the last two strings can not be changed to the
first web presented in Figure~\ref{a2ab2exppo} with a $U$-turn. Now
we look at the $D_{i,a}$ where $i>1$ as illustrated in
Figure~\ref{a2abexp5-1}. Since we picked where the $U$ turn appears
already, only possible disjoint lines are given as thick and lightly
shaded lines but we can not finish to have a stem because the darkly
shaded string from the left top can not be connected to the bottom
clasp without being zero, $i. e.,$ if we connect the tick line to
clasp, there will be either
$\begin{pspicture}[shift=.2](-.17,-.17)(.17,.17)\qline(0,0)(0,.17)\qline(-.15,-.15)(0,0)\qline(.15,-.15)(0,0)\end{pspicture}$
or a
$\begin{pspicture}[shift=.2](-.17,-.17)(.17,.17)\qline(-.15,-.15)(-.15,.05)\qline(.15,-.15)(.15,.05)
\pccurve[angleA=90,angleB=90,ncurv=1](-.15,.05)(.15,.05)\end{pspicture}$
right above of the clasp of weight $(a,b-1)$.

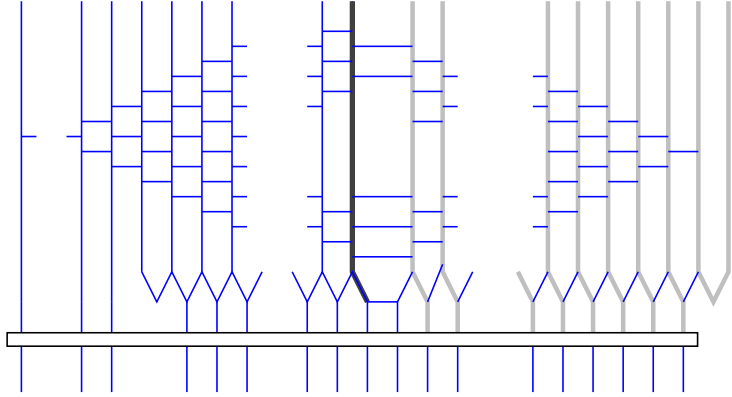
\begin{figure}
$$
\begin{pspicture}[shift=-2.7](-.2,-.3)(10.1,5.1)
\qline(0.0,0.6)(0.0,1.0)\qline(0.8,0.6)(0.8,1.0)\qline(1.2,0.6)(1.2,1.0)
\qline(2.2,0.6)(2.2,1.0)\qline(2.6,0.6)(2.6,1.0)\qline(3.0,0.6)(3.0,1.0)
\qline(3.8,0.6)(3.8,1.0)\qline(4.2,0.6)(4.2,1.0)\qline(4.6,0.6)(4.6,1.0)
\qline(5.0,0.6)(5.0,1.0)\qline(0.0,0.4)(0.0,-0.2)\qline(0.8,0.4)(0.8,-0.2)
\qline(1.2,0.4)(1.2,-0.2)\qline(2.2,0.4)(2.2,-0.2)\qline(2.6,0.4)(2.6,-0.2)
\qline(3.0,0.4)(3.0,-0.2)\qline(3.8,0.4)(3.8,-0.2)\qline(4.2,0.4)(4.2,-0.2)
\qline(4.6,0.4)(4.6,-0.2)\qline(5.0,0.4)(5.0,-0.2)\qline(5.4,0.4)(5.4,-0.2)
\qline(5.8,0.4)(5.8,-0.2)\qline(6.8,0.4)(6.8,-0.2)\qline(7.2,0.4)(7.2,-0.2)
\qline(7.6,0.4)(7.6,-0.2)\qline(8.0,0.4)(8.0,-0.2)\qline(8.4,0.4)(8.4,-0.2)
\qline(8.8,0.4)(8.8,-0.2)\qline(0.0,1.0)(0.,5.0)\qline(0.8,1.0)(0.8,5.0)
\qline(1.2,1.0)(1.2,5.0)\qline(1.6,1.4)(1.6,5.0)\qline(2.0,1.4)(2.0,5.0)
\qline(2.4,1.4)(2.4,5.0)\qline(2.8,1.4)(2.8,5.0)\qline(4.0,1.4)(4.0,5.0)
\qline(4.4,1.4)(4.4,5.0)
\psline[linecolor=lightgray,linewidth=1.7pt](6.6,1.4)(6.8,1)(6.8,.6)
\psline[linecolor=lightgray,linewidth=1.7pt](5.4,.6)(5.4,1)(5.2,1.4)(5.2,5.0)
\psline[linecolor=lightgray,linewidth=1.7pt](5.8,.6)(5.8,1)(5.6,1.4)(5.6,5.0)
\psline[linecolor=lightgray,linewidth=1.7pt](7.2,.6)(7.2,1)(7.0,1.4)(7.0,5.0)
\psline[linecolor=lightgray,linewidth=1.7pt](7.6,.6)(7.6,1)(7.4,1.4)(7.4,5.0)
\psline[linecolor=lightgray,linewidth=1.7pt](8,.6)(8,1)(7.8,1.4)(7.8,5.0)
\psline[linecolor=lightgray,linewidth=1.7pt](8.4,.6)(8.4,1)(8.2,1.4)(8.2,5.0)
\psline[linecolor=lightgray,linewidth=1.7pt](8.8,.6)(8.8,1)(8.6,1.4)(8.6,5.0)
\psline[linecolor=lightgray,linewidth=1.7pt](9.4,5)(9.4,1.4)(9.2,1)(9.0,1.4)(9.0,5.0)
\psline[linecolor=darkgray,linewidth=2pt](4.4,5)(4.4,1.4)(4.6,1)
\psline(1.6,1.4)(1.8,1)(2,1.4)(2.2,1)(2.4,1.4)(2.6,1)(2.8,1.4)(3,1)(3.2,1.4)
\psline(3.6,1.4)(3.8,1)(4,1.4)(4.2,1)(4.4,1.4)(4.6,1)(5,1)(5.2,1.4)
\qline(5.4,1)(5.6,1.5)\qline(5.8,1)(6,1.4)
\qline(6.8,1)(7,1.4)\qline(7.2,1)(7.4,1.4)\qline(7.6,1)(7.8,1.4)
\qline(8,1)(8.2,1.4)\qline(8.4,1)(8.6,1.4)\qline(8.8,1)(9,1.4)
\qline(0,3.2)(.2,3.2)
\qline(0.6,3.2)(0.8,3.2)\qline(0.8,3.0)(1.2,3.0)\qline(0.8,3.4)(1.2,3.4)
\qline(1.2,2.8)(1.6,2.8)\qline(1.2,3.2)(1.6,3.2)\qline(1.2,3.6)(1.6,3.6)
\qline(1.6,2.6)(2.0,2.6)\qline(1.6,3.0)(2.0,3.0)\qline(1.6,3.4)(2.0,3.4)
\qline(1.6,3.8)(2.0,3.8)\qline(2.0,2.4)(2.4,2.4)\qline(2.0,2.8)(2.4,2.8)
\qline(2.0,3.2)(2.4,3.2)\qline(2.0,3.6)(2.4,3.6)\qline(2.0,4.0)(2.4,4.0)
\qline(2.4,2.2)(2.8,2.2)\qline(2.4,2.6)(2.8,2.6)\qline(2.4,3.0)(2.8,3.0)
\qline(2.4,3.4)(2.8,3.4)\qline(2.4,3.8)(2.8,3.8)\qline(2.4,4.2)(2.8,4.2)
\qline(2.8,2.0)(3.0,2.0)\qline(2.8,2.4)(3.0,2.4)\qline(2.8,2.8)(3.0,2.8)
\qline(2.8,3.2)(3.0,3.2)\qline(2.8,3.6)(3.0,3.6)\qline(2.8,4.0)(3.0,4.0)
\qline(2.8,4.4)(3.0,4.4)\qline(3.8,2.0)(4.0,2.0)\qline(3.8,2.4)(4.0,2.4)
\qline(3.8,3.6)(4.0,3.6)\qline(3.8,4.0)(4.0,4.0)\qline(3.8,4.4)(4.0,4.4)
\qline(4.0,1.8)(4.4,1.8)\qline(4.0,2.2)(4.4,2.2)\qline(4.0,3.8)(4.4,3.8)
\qline(4.0,4.2)(4.4,4.2)\qline(4.0,4.6)(4.4,4.6)\qline(4.4,1.6)(5.2,1.6)
\qline(5.2,2.0)(4.4,2.0)\qline(5.2,2.4)(4.4,2.4)\qline(5.2,4.0)(4.4,4.0)
\qline(5.2,4.4)(4.4,4.4)\qline(5.2,1.8)(5.6,1.8)\qline(5.2,2.2)(5.6,2.2)
\qline(5.2,3.8)(5.6,3.8)\qline(5.2,4.2)(5.6,4.2)\qline(5.2,3.4)(5.6,3.4)
\qline(5.6,2.0)(5.8,2.0)\qline(5.6,2.4)(5.8,2.4)\qline(5.6,4.0)(5.8,4.0)
\qline(5.6,3.6)(5.8,3.6)\qline(7.0,2.0)(6.8,2.0)\qline(6.8,2.4)(7.0,2.4)
\qline(7.0,4.0)(6.8,4.0)\qline(6.8,3.6)(7.0,3.6)\qline(7.0,2.2)(7.4,2.2)
\qline(7.0,2.6)(7.4,2.6)\qline(7.0,3.0)(7.4,3.0)\qline(7.0,3.4)(7.4,3.4)
\qline(7.0,3.8)(7.4,3.8)\qline(7.4,2.4)(7.8,2.4)\qline(7.4,2.8)(7.8,2.8)
\qline(7.4,3.2)(7.8,3.2)\qline(7.4,3.6)(7.8,3.6)\qline(7.8,2.6)(8.2,2.6)
\qline(7.8,3.0)(8.2,3.0)\qline(7.8,3.4)(8.2,3.4)\qline(8.2,2.8)(8.6,2.8)
\qline(8.2,3.2)(8.6,3.2)\qline(8.6,3)(9.0,3)
\psframe[linecolor=black](-.2,.4)(9,.6)
\end{pspicture}
$$   \caption{$D_{i,a}$ where $i>1$.} \label{a2abexp5-1}\end{figure}

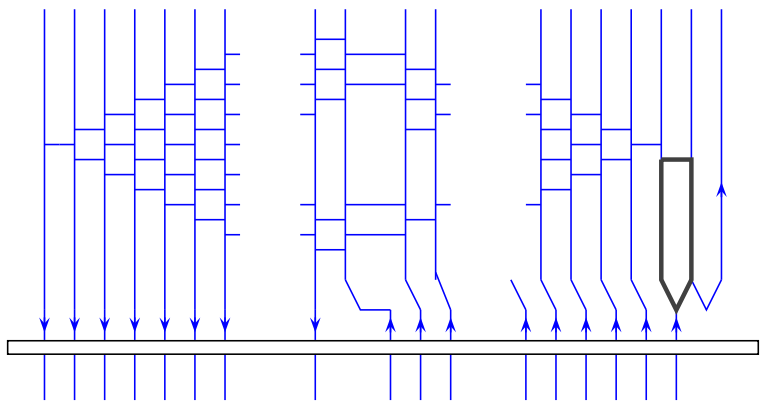
\begin{figure}
$$
\begin{pspicture}[shift=-2.7](-.2,-.3)(10.1,5.1)
\qline(0.4,0.6)(0.4,1.0)\psline[arrowscale=1.5]{->}(0.4,.9)(0.4,.7)
\qline(0.8,0.6)(0.8,1.0)\psline[arrowscale=1.5]{->}(0.8,.9)(0.8,.7)
\qline(1.2,0.6)(1.2,1.0)\psline[arrowscale=1.5]{->}(1.2,.9)(1.2,.7)
\qline(1.6,0.6)(1.6,1.0)\psline[arrowscale=1.5]{->}(1.6,.9)(1.6,.7)
\qline(2,0.6)(2,1.0)\psline[arrowscale=1.5]{->}(2,.9)(2,.7)
\qline(2.4,0.6)(2.4,1.0)\psline[arrowscale=1.5]{->}(2.4,.9)(2.4,.7)
\qline(2.8,0.6)(2.8,1.0)\psline[arrowscale=1.5]{->}(2.8,.9)(2.8,.7)
\qline(4.0,0.6)(4.0,1.0)\psline[arrowscale=1.5]{->}(4.0,.9)(4.0,.7)
\qline(5.0,0.6)(5.0,1.0)\psline[arrowscale=1.5]{<-}(5.0,.9)(5.0,.7)
\qline(5.4,0.6)(5.4,1.0)\psline[arrowscale=1.5]{<-}(5.4,.9)(5.4,.7)
\qline(5.8,0.6)(5.8,1.0)\psline[arrowscale=1.5]{<-}(5.8,.9)(5.8,.7)
\qline(6.8,0.6)(6.8,1.0)\psline[arrowscale=1.5]{<-}(6.8,.9)(6.8,.7)
\qline(7.2,0.6)(7.2,1.0)\psline[arrowscale=1.5]{<-}(7.2,.9)(7.2,.7)
\qline(7.6,0.6)(7.6,1.0)\psline[arrowscale=1.5]{<-}(7.6,.9)(7.6,.7)
\qline(8.0,0.6)(8.0,1.0)\psline[arrowscale=1.5]{<-}(8,.9)(8,.7)
\qline(8.4,0.6)(8.4,1.0)\psline[arrowscale=1.5]{<-}(8.4,.9)(8.4,.7)
\qline(8.8,0.6)(8.8,1.0)\psline[arrowscale=1.5]{<-}(8.8,.9)(8.8,.7)
\qline(0.4,0.4)(0.4,-0.2) \qline(0.8,0.4)(0.8,-0.2)
\qline(1.2,0.4)(1.2,-0.2) \qline(1.6,0.4)(1.6,-0.2)
\qline(2,0.4)(2,-0.2) \qline(2.4,0.4)(2.4,-0.2)
\qline(2.8,0.4)(2.8,-0.2) \qline(4,0.4)(4,-0.2)
\qline(5.0,0.4)(5.0,-0.2) \qline(5.4,0.4)(5.4,-0.2)
\qline(5.8,0.4)(5.8,-0.2) \qline(6.8,0.4)(6.8,-0.2)
\qline(7.2,0.4)(7.2,-0.2) \qline(7.6,0.4)(7.6,-0.2)
\qline(8.0,0.4)(8.0,-0.2) \qline(8.4,0.4)(8.4,-0.2)
\qline(8.8,0.4)(8.8,-0.2)
\qline(0.4,1.0)(0.4,5.0)\qline(0.8,1.0)(0.8,5.0)\qline(1.2,1.0)(1.2,5.0)
\qline(1.6,1.0)(1.6,5.0)\qline(2.0,1.0)(2.0,5.0)\qline(2.4,1.0)(2.4,5.0)
\qline(2.8,1.0)(2.8,5.0)\qline(4.0,1.0)(4.0,5.0)\qline(4.4,1.4)(4.4,5.0)
\qline(5.2,1.4)(5.2,5.0)\qline(5.6,1.4)(5.6,5.0)\qline(7.0,1.4)(7.0,5.0)
\qline(7.4,1.4)(7.4,5.0)\qline(7.8,1.4)(7.8,5.0)\qline(8.2,1.4)(8.2,5.0)
\qline(9.0,1.4)(9.0,5.0)\qline(8.6,1.4)(8.6,5.0)
\psline(4.4,1.4)(4.6,1.0)(5.0,1.0)\qline(5.2,1.4)(5.4,1.0)
\qline(5.6,1.5)(5.8,1.0)\qline(6.6,1.4)(6.8,1.0)\qline(7.0,1.4)(7.2,1.0)
\qline(7.4,1.4)(7.6,1.0)\qline(7.8,1.4)(8,1)\qline(8.2,1.4)(8.4,1.0)
\qline(8.6,1.4)(8.8,1)\psline(9,1.4)(9.2,1)(9.4,1.4)\qline(0.4,3.2)(.6,3.2)
\qline(0.6,3.2)(0.8,3.2)\qline(0.8,3.0)(1.2,3.0)\qline(0.8,3.4)(1.2,3.4)
\qline(1.2,2.8)(1.6,2.8)\qline(1.2,3.2)(1.6,3.2)\qline(1.2,3.6)(1.6,3.6)
\qline(1.6,2.6)(2.0,2.6)\qline(1.6,3.0)(2.0,3.0)\qline(1.6,3.4)(2.0,3.4)
\qline(1.6,3.8)(2.0,3.8)\qline(2.0,2.4)(2.4,2.4)\qline(2.0,2.8)(2.4,2.8)
\qline(2.0,3.2)(2.4,3.2)\qline(2.0,3.6)(2.4,3.6)\qline(2.0,4.0)(2.4,4.0)
\qline(2.4,2.2)(2.8,2.2)\qline(2.4,2.6)(2.8,2.6)\qline(2.4,3.0)(2.8,3.0)
\qline(2.4,3.4)(2.8,3.4)\qline(2.4,3.8)(2.8,3.8)\qline(2.4,4.2)(2.8,4.2)
\qline(2.8,2.0)(3.0,2.0)\qline(2.8,2.4)(3.0,2.4)\qline(2.8,2.8)(3.0,2.8)
\qline(2.8,3.2)(3.0,3.2)\qline(2.8,3.6)(3.0,3.6)\qline(2.8,4.0)(3.0,4.0)
\qline(2.8,4.4)(3.0,4.4)\qline(3.8,2.0)(4.0,2.0)\qline(3.8,2.4)(4.0,2.4)
\qline(3.8,3.6)(4.0,3.6)\qline(3.8,4.0)(4.0,4.0)\qline(3.8,4.4)(4.0,4.4)
\qline(4.0,1.8)(4.4,1.8)\qline(4.0,2.2)(4.4,2.2)\qline(4.0,3.8)(4.4,3.8)
\qline(4.0,4.2)(4.4,4.2)\qline(4.0,4.6)(4.4,4.6)\qline(5.2,2.0)(4.4,2.0)
\qline(5.2,2.4)(4.4,2.4)\qline(5.2,4.0)(4.4,4.0)\qline(5.2,4.4)(4.4,4.4)
\qline(5.2,2.2)(5.6,2.2)\qline(5.2,3.8)(5.6,3.8)\qline(5.2,4.2)(5.6,4.2)
\qline(5.2,3.4)(5.6,3.4)\qline(5.6,2.4)(5.8,2.4)\qline(5.6,4.0)(5.8,4.0)
\qline(5.6,3.6)(5.8,3.6)\qline(6.8,2.4)(7.0,2.4)\qline(7.0,4.0)(6.8,4.0)
\qline(6.8,3.6)(7.0,3.6)\qline(7.0,2.6)(7.4,2.6)\qline(7.0,3.0)(7.4,3.0)
\qline(7.0,3.4)(7.4,3.4)\qline(7.0,3.8)(7.4,3.8)\qline(7.4,2.8)(7.8,2.8)
\qline(7.4,3.2)(7.8,3.2)\qline(7.4,3.6)(7.8,3.6)\qline(7.8,3.0)(8.2,3.0)
\qline(7.8,3.4)(8.2,3.4)\qline(8.2,3.2)(8.6,3.2)\qline(9.4,1.4)(9.4,5)
\psline[arrowscale=1.5]{->}(9.4,2.5)(9.4,2.7)
\psline[linecolor=darkgray,linewidth=1.7pt](8.6,3)(9,3)(9,1.4)(8.8,1)(8.6,1.4)(8.6,3)
\psframe[linecolor=black](-.1,.4)(9.9,.6)
\end{pspicture}
$$  \caption{The web $\tilde D_{1,a}$. } \label{a2abexp6}
\end{figure}

So only nonzero admissible stems should be obtained from $\tilde
D_{1,a}$. As we explained before, one can see that there is a
rectangular face in the web $\tilde D_{1,a}$. Since the horizontal
splitting makes it zero, we have to split vertically. For the
resulting web, this process created one rectangular face at right
topside of previous place. We have to split vertically and the
process are repeated until the last step, both splits do not vanish.
The web in the last step is drawn in Figure~\ref{a2abexp6} with the
rectangular face, darkly shaded. The vertical split gives us the
first web in Figure \ref{a2ab2exppo} and the horizontal split gives
the third web in Figure \ref{a2ab2exppo}.

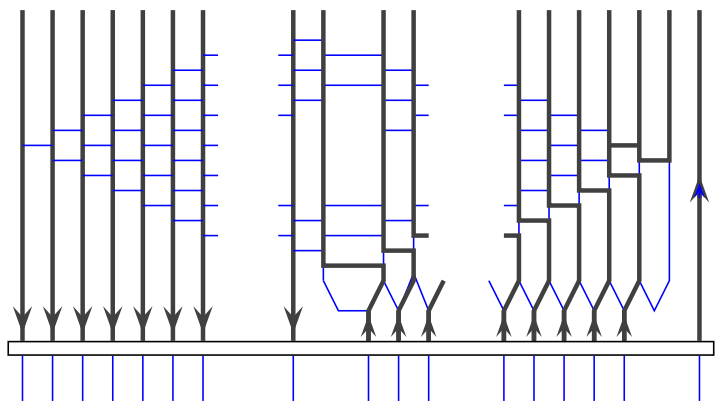
\begin{figure}
$$
\begin{pspicture}[shift=-2.7](-.2,-.3)(10.1,5.1)
\psline[linecolor=darkgray,linewidth=1.7pt](0.4,0.6)(0.4,5.0)
\psline[arrowscale=1.5,linecolor=darkgray,linewidth=1.7pt]{->}(0.4,.9)(0.4,.7)
\psline[linecolor=darkgray,linewidth=1.7pt](0.8,0.6)(0.8,5.0)
\psline[arrowscale=1.5,linecolor=darkgray,linewidth=1.7pt]{->}(0.8,.9)(0.8,.7)
\psline[linecolor=darkgray,linewidth=1.7pt](1.2,0.6)(1.2,5.0)
\psline[arrowscale=1.5,linecolor=darkgray,linewidth=1.7pt]{->}(1.2,.9)(1.2,.7)
\psline[linecolor=darkgray,linewidth=1.7pt](1.6,0.6)(1.6,5.0)
\psline[arrowscale=1.5,linecolor=darkgray,linewidth=1.7pt]{->}(1.6,.9)(1.6,.7)
\psline[linecolor=darkgray,linewidth=1.7pt](2,0.6)(2,5.0)
\psline[arrowscale=1.5,linecolor=darkgray,linewidth=1.7pt]{->}(2,.9)(2,.7)
\psline[linecolor=darkgray,linewidth=1.7pt](2.4,0.6)(2.4,5.0)
\psline[arrowscale=1.5,linecolor=darkgray,linewidth=1.7pt]{->}(2.4,.9)(2.4,.7)
\psline[linecolor=darkgray,linewidth=1.7pt](2.8,0.6)(2.8,5.0)
\psline[arrowscale=1.5,linecolor=darkgray,linewidth=1.7pt]{->}(2.8,.9)(2.8,.7)
\psline[linecolor=darkgray,linewidth=1.7pt](4.0,0.6)(4.0,5.0)
\psline[arrowscale=1.5,linecolor=darkgray,linewidth=1.7pt]{->}(4.0,.9)(4.0,.7)
\psline[linecolor=darkgray,linewidth=1.7pt](5.0,0.6)(5.0,1.0)
\psline[arrowscale=1.2,linecolor=darkgray,linewidth=1.7pt]{<-}(5.0,.94)(5.0,.8)
\psline[linecolor=darkgray,linewidth=1.7pt](5.4,0.6)(5.4,1.0)
\psline[arrowscale=1.2,linecolor=darkgray,linewidth=1.7pt]{<-}(5.4,.94)(5.4,.8)
\psline[linecolor=darkgray,linewidth=1.7pt](5.8,0.6)(5.8,1.0)
\psline[arrowscale=1.2,linecolor=darkgray,linewidth=1.7pt]{<-}(5.8,.94)(5.8,.8)
\psline[linecolor=darkgray,linewidth=1.7pt](6.8,0.6)(6.8,1.0)
\psline[arrowscale=1.2,linecolor=darkgray,linewidth=1.7pt]{<-}(6.8,.94)(6.8,.8)
\psline[linecolor=darkgray,linewidth=1.7pt](7.2,0.6)(7.2,1.0)
\psline[arrowscale=1.2,linecolor=darkgray,linewidth=1.7pt]{<-}(7.2,.94)(7.2,.8)
\psline[linecolor=darkgray,linewidth=1.7pt](7.6,0.6)(7.6,1.0)
\psline[arrowscale=1.2,linecolor=darkgray,linewidth=1.7pt]{<-}(7.6,.94)(7.6,.8)
\psline[linecolor=darkgray,linewidth=1.7pt](8.0,0.6)(8.0,1.0)
\psline[arrowscale=1.2,linecolor=darkgray,linewidth=1.7pt]{<-}(8,.94)(8,.8)
\psline[linecolor=darkgray,linewidth=1.7pt](8.4,0.6)(8.4,1.0)
\psline[arrowscale=1.2,linecolor=darkgray,linewidth=1.7pt]{<-}(8.4,.94)(8.4,.8)
\psline[linecolor=darkgray,linewidth=1.7pt](9.4,0.6)(9.4,5.0)
\psline[arrowscale=1.5,linecolor=darkgray,linewidth=1.7pt]{<-}(9.4,2.8)(9.4,2.5)
\qline(0.4,0.4)(0.4,-0.2) \qline(0.8,0.4)(0.8,-0.2)
\qline(1.2,0.4)(1.2,-0.2) \qline(1.6,0.4)(1.6,-0.2)
\qline(2,0.4)(2,-0.2) \qline(2.4,0.4)(2.4,-0.2)
\qline(2.8,0.4)(2.8,-0.2) \qline(4,0.4)(4,-0.2)
\qline(5.0,0.4)(5.0,-0.2) \qline(5.4,0.4)(5.4,-0.2)
\qline(5.8,0.4)(5.8,-0.2) \qline(6.8,0.4)(6.8,-0.2)
\qline(7.2,0.4)(7.2,-0.2) \qline(7.6,0.4)(7.6,-0.2)
\qline(8.0,0.4)(8.0,-0.2) \qline(8.4,0.4)(8.4,-0.2)
\qline(9.4,0.4)(9.4,-0.2) \qline(4.4,1.4)(4.4,1.6)
\qline(5.2,1.6)(5.2,1.8)\qline(5.6,1.8)(5.6,2.0)\qline(7.0,2)(7.0,2.2)
\qline(7.4,2.4)(7.4,2.2)\qline(7.8,2.4)(7.8,2.6)\qline(8.2,2.6)(8.2,2.8)
\qline(8.6,2.8)(8.6,3.0)\psline(9,3)(9,1.4)(8.8,1)(8.6,1.4)
\psline(4.4,1.4)(4.6,1.0)(5.0,1.0)(5.2,1.4)(5.4,1.0)(5.6,1.5)(5.8,1.0)
\qline(6.6,1.4)(6.8,1)\qline(7,1.4)(7.2,1)\qline(7.4,1.4)(7.6,1)\qline(7.8,1.4)(8,1)
\qline(8.2,1.4)(8.4,1)
\qline(0.4,3.2)(0.8,3.2)\qline(0.8,3.0)(1.2,3.0)\qline(0.8,3.4)(1.2,3.4)
\qline(1.2,2.8)(1.6,2.8)\qline(1.2,3.2)(1.6,3.2)\qline(1.2,3.6)(1.6,3.6)
\qline(1.6,2.6)(2.0,2.6)\qline(1.6,3.0)(2.0,3.0)\qline(1.6,3.4)(2.0,3.4)
\qline(1.6,3.8)(2.0,3.8)\qline(2.0,2.4)(2.4,2.4)\qline(2.0,2.8)(2.4,2.8)
\qline(2.0,3.2)(2.4,3.2)\qline(2.0,3.6)(2.4,3.6)\qline(2.0,4.0)(2.4,4.0)
\qline(2.4,2.2)(2.8,2.2)\qline(2.4,2.6)(2.8,2.6)\qline(2.4,3.0)(2.8,3.0)
\qline(2.4,3.4)(2.8,3.4)\qline(2.4,3.8)(2.8,3.8)\qline(2.4,4.2)(2.8,4.2)
\qline(2.8,2.0)(3.0,2.0)\qline(2.8,2.4)(3.0,2.4)\qline(2.8,2.8)(3.0,2.8)
\qline(2.8,3.2)(3.0,3.2)\qline(2.8,3.6)(3.0,3.6)\qline(2.8,4.0)(3.0,4.0)
\qline(2.8,4.4)(3.0,4.4)\qline(3.8,2.0)(4.0,2.0)\qline(3.8,2.4)(4.0,2.4)
\qline(3.8,3.6)(4.0,3.6)\qline(3.8,4.0)(4.0,4.0)\qline(3.8,4.4)(4.0,4.4)
\qline(4.0,1.8)(4.4,1.8)\qline(4.0,2.2)(4.4,2.2)\qline(4.0,3.8)(4.4,3.8)
\qline(4.0,4.2)(4.4,4.2)\qline(4.0,4.6)(4.4,4.6)\qline(5.2,2.0)(4.4,2.0)
\qline(5.2,2.4)(4.4,2.4)\qline(5.2,4.0)(4.4,4.0)\qline(5.2,4.4)(4.4,4.4)
\qline(5.2,2.2)(5.6,2.2)\qline(5.2,3.8)(5.6,3.8)\qline(5.2,4.2)(5.6,4.2)
\qline(5.2,3.4)(5.6,3.4)\qline(5.6,2.4)(5.8,2.4)\qline(5.6,4.0)(5.8,4.0)
\qline(5.6,3.6)(5.8,3.6)\qline(6.8,2.4)(7.0,2.4)\qline(7.0,4.0)(6.8,4.0)
\qline(6.8,3.6)(7.0,3.6)\qline(7.0,2.6)(7.4,2.6)\qline(7.0,3.0)(7.4,3.0)
\qline(7.0,3.4)(7.4,3.4)\qline(7.0,3.8)(7.4,3.8)\qline(7.4,2.8)(7.8,2.8)
\qline(7.4,3.2)(7.8,3.2)\qline(7.4,3.6)(7.8,3.6)\qline(7.8,3.0)(8.2,3.0)
\qline(7.8,3.4)(8.2,3.4)
\psline[linecolor=darkgray,linewidth=1.7pt](4.4,5)(4.4,1.6)(5.2,1.6)(5.2,1.4)(5.0,1)(5.0,.6)
\psline[linecolor=darkgray,linewidth=1.7pt](5.2,5)(5.2,1.8)(5.6,1.8)(5.6,1.4)(5.4,1)(5.4,.6)
\psline[linecolor=darkgray,linewidth=1.7pt](5.6,5)(5.6,2.0)(5.8,2.0)
\psline[linecolor=darkgray,linewidth=1.7pt](6,1.4)(5.8,1)(5.8,.6)
\psline[linecolor=darkgray,linewidth=1.7pt](6.8,2)(7,2)(7,1.4)(6.8,1)(6.8,.6)
\psline[linecolor=darkgray,linewidth=1.7pt](7,5)(7,2.2)(7.4,2.2)(7.4,1.4)(7.2,1)(7.2,.6)
\psline[linecolor=darkgray,linewidth=1.7pt](7.4,5)(7.4,2.4)(7.8,2.4)(7.8,1.4)(7.6,1)(7.6,.6)
\psline[linecolor=darkgray,linewidth=1.7pt](7.8,5)(7.8,2.6)(8.2,2.6)(8.2,1.4)(8,1)(8,.6)
\psline[linecolor=darkgray,linewidth=1.7pt](8.2,5)(8.2,2.8)(8.6,2.8)(8.6,1.4)(8.4,1)(8.4,.6)
\psline[linecolor=darkgray,linewidth=1.7pt](8.6,5)(8.6,3)(9,3)(9,5)
\psline[linecolor=darkgray,linewidth=1.7pt](8.2,3.2)(8.6,3.2)
\psline[arrowscale=1.5]{->}(9.4,2.5)(9.4,2.7)
\psframe[linecolor=black](.2,.4)(9.6,.6)
\end{pspicture}
$$ \caption{The nonzero admissible stem for $\tilde D_{1,a-1}$. }
\label{a2abexp7}\end{figure}

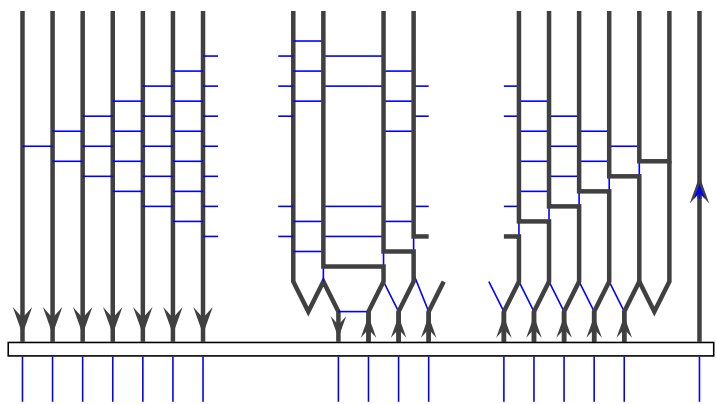
\begin{figure}
$$
\begin{pspicture}[shift=-2.7](-.2,-.3)(10.1,5.1)
\psline[linecolor=darkgray,linewidth=1.7pt](0.4,0.6)(0.4,5.0)
\psline[arrowscale=1.5,linecolor=darkgray,linewidth=1.7pt]{->}(0.4,.9)(0.4,.7)
\psline[linecolor=darkgray,linewidth=1.7pt](0.8,0.6)(0.8,5.0)
\psline[arrowscale=1.5,linecolor=darkgray,linewidth=1.7pt]{->}(0.8,.9)(0.8,.7)
\psline[linecolor=darkgray,linewidth=1.7pt](1.2,0.6)(1.2,5.0)
\psline[arrowscale=1.5,linecolor=darkgray,linewidth=1.7pt]{->}(1.2,.9)(1.2,.7)
\psline[linecolor=darkgray,linewidth=1.7pt](1.6,0.6)(1.6,5.0)
\psline[arrowscale=1.5,linecolor=darkgray,linewidth=1.7pt]{->}(1.6,.9)(1.6,.7)
\psline[linecolor=darkgray,linewidth=1.7pt](2,0.6)(2,5.0)
\psline[arrowscale=1.5,linecolor=darkgray,linewidth=1.7pt]{->}(2,.9)(2,.7)
\psline[linecolor=darkgray,linewidth=1.7pt](2.4,0.6)(2.4,5.0)
\psline[arrowscale=1.5,linecolor=darkgray,linewidth=1.7pt]{->}(2.4,.9)(2.4,.7)
\psline[linecolor=darkgray,linewidth=1.7pt](2.8,0.6)(2.8,5.0)
\psline[arrowscale=1.5,linecolor=darkgray,linewidth=1.7pt]{->}(2.8,.9)(2.8,.7)
\psline[linecolor=darkgray,linewidth=1.7pt](4,5)(4,1.4)(4.2,1)(4.4,1.4)(4.6,1)(4.6,.6)
\psline[arrowscale=1.2,linecolor=darkgray,linewidth=1.7pt]{->}(4.6,.8)(4.6,.65)
\psline[linecolor=darkgray,linewidth=1.7pt](5.0,0.6)(5.0,1.0)
\psline[arrowscale=1.2,linecolor=darkgray,linewidth=1.7pt]{<-}(5.0,.94)(5.0,.8)
\psline[linecolor=darkgray,linewidth=1.7pt](5.4,0.6)(5.4,1.0)
\psline[arrowscale=1.2,linecolor=darkgray,linewidth=1.7pt]{<-}(5.4,.94)(5.4,.8)
\psline[linecolor=darkgray,linewidth=1.7pt](5.8,0.6)(5.8,1.0)
\psline[arrowscale=1.2,linecolor=darkgray,linewidth=1.7pt]{<-}(5.8,.94)(5.8,.8)
\psline[linecolor=darkgray,linewidth=1.7pt](6.8,0.6)(6.8,1.0)
\psline[arrowscale=1.2,linecolor=darkgray,linewidth=1.7pt]{<-}(6.8,.94)(6.8,.8)
\psline[linecolor=darkgray,linewidth=1.7pt](7.2,0.6)(7.2,1.0)
\psline[arrowscale=1.2,linecolor=darkgray,linewidth=1.7pt]{<-}(7.2,.94)(7.2,.8)
\psline[linecolor=darkgray,linewidth=1.7pt](7.6,0.6)(7.6,1.0)
\psline[arrowscale=1.2,linecolor=darkgray,linewidth=1.7pt]{<-}(7.6,.94)(7.6,.8)
\psline[linecolor=darkgray,linewidth=1.7pt](8.0,0.6)(8.0,1.0)
\psline[arrowscale=1.2,linecolor=darkgray,linewidth=1.7pt]{<-}(8,.94)(8,.8)
\psline[linecolor=darkgray,linewidth=1.7pt](8.4,0.6)(8.4,1.0)
\psline[arrowscale=1.2,linecolor=darkgray,linewidth=1.7pt]{<-}(8.4,.94)(8.4,.8)
\psline[linecolor=darkgray,linewidth=1.7pt](9.4,0.6)(9.4,5.0)
\psline[arrowscale=1.5,linecolor=darkgray,linewidth=1.7pt]{<-}(9.4,2.8)(9.4,2.5)
\qline(0.4,0.4)(0.4,-0.2) \qline(0.8,0.4)(0.8,-0.2)
\qline(1.2,0.4)(1.2,-0.2) \qline(1.6,0.4)(1.6,-0.2)
\qline(2,0.4)(2,-0.2) \qline(2.4,0.4)(2.4,-0.2)
\qline(2.8,0.4)(2.8,-0.2) \qline(4.6,0.4)(4.6,-0.2)
\qline(5.0,0.4)(5.0,-0.2) \qline(5.4,0.4)(5.4,-0.2)
\qline(5.8,0.4)(5.8,-0.2) \qline(6.8,0.4)(6.8,-0.2)
\qline(7.2,0.4)(7.2,-0.2) \qline(7.6,0.4)(7.6,-0.2)
\qline(8.0,0.4)(8.0,-0.2) \qline(8.4,0.4)(8.4,-0.2)
\qline(9.4,0.4)(9.4,-0.2) \qline(4.4,1.4)(4.4,1.6)
\qline(5.2,1.6)(5.2,1.8)\qline(5.6,1.8)(5.6,2.0)\qline(7.0,2)(7.0,2.2)
\qline(7.4,2.4)(7.4,2.2)\qline(7.8,2.4)(7.8,2.6)\qline(8.2,2.6)(8.2,2.8)
\qline(8.6,2.8)(8.6,3.0)
\qline(4.6,1.0)(5.0,1.0)\qline(5.2,1.4)(5.4,1.0)\qline(5.6,1.5)(5.8,1.0)
\qline(6.6,1.4)(6.8,1)\qline(7,1.4)(7.2,1)\qline(7.4,1.4)(7.6,1)\qline(7.8,1.4)(8,1)
\qline(8.2,1.4)(8.4,1)
\qline(0.4,3.2)(0.8,3.2)\qline(0.8,3.0)(1.2,3.0)\qline(0.8,3.4)(1.2,3.4)
\qline(1.2,2.8)(1.6,2.8)\qline(1.2,3.2)(1.6,3.2)\qline(1.2,3.6)(1.6,3.6)
\qline(1.6,2.6)(2.0,2.6)\qline(1.6,3.0)(2.0,3.0)\qline(1.6,3.4)(2.0,3.4)
\qline(1.6,3.8)(2.0,3.8)\qline(2.0,2.4)(2.4,2.4)\qline(2.0,2.8)(2.4,2.8)
\qline(2.0,3.2)(2.4,3.2)\qline(2.0,3.6)(2.4,3.6)\qline(2.0,4.0)(2.4,4.0)
\qline(2.4,2.2)(2.8,2.2)\qline(2.4,2.6)(2.8,2.6)\qline(2.4,3.0)(2.8,3.0)
\qline(2.4,3.4)(2.8,3.4)\qline(2.4,3.8)(2.8,3.8)\qline(2.4,4.2)(2.8,4.2)
\qline(2.8,2.0)(3.0,2.0)\qline(2.8,2.4)(3.0,2.4)\qline(2.8,2.8)(3.0,2.8)
\qline(2.8,3.2)(3.0,3.2)\qline(2.8,3.6)(3.0,3.6)\qline(2.8,4.0)(3.0,4.0)
\qline(2.8,4.4)(3.0,4.4)\qline(3.8,2.0)(4.0,2.0)\qline(3.8,2.4)(4.0,2.4)
\qline(3.8,3.6)(4.0,3.6)\qline(3.8,4.0)(4.0,4.0)\qline(3.8,4.4)(4.0,4.4)
\qline(4.0,1.8)(4.4,1.8)\qline(4.0,2.2)(4.4,2.2)\qline(4.0,3.8)(4.4,3.8)
\qline(4.0,4.2)(4.4,4.2)\qline(4.0,4.6)(4.4,4.6)\qline(5.2,2.0)(4.4,2.0)
\qline(5.2,2.4)(4.4,2.4)\qline(5.2,4.0)(4.4,4.0)\qline(5.2,4.4)(4.4,4.4)
\qline(5.2,2.2)(5.6,2.2)\qline(5.2,3.8)(5.6,3.8)\qline(5.2,4.2)(5.6,4.2)
\qline(5.2,3.4)(5.6,3.4)\qline(5.6,2.4)(5.8,2.4)\qline(5.6,4.0)(5.8,4.0)
\qline(5.6,3.6)(5.8,3.6)\qline(6.8,2.4)(7.0,2.4)\qline(7.0,4.0)(6.8,4.0)
\qline(6.8,3.6)(7.0,3.6)\qline(7.0,2.6)(7.4,2.6)\qline(7.0,3.0)(7.4,3.0)
\qline(7.0,3.4)(7.4,3.4)\qline(7.0,3.8)(7.4,3.8)\qline(7.4,2.8)(7.8,2.8)
\qline(7.4,3.2)(7.8,3.2)\qline(7.4,3.6)(7.8,3.6)\qline(7.8,3.0)(8.2,3.0)
\qline(7.8,3.4)(8.2,3.4)\qline(8.2,3.2)(8.6,3.2)
\psline[linecolor=darkgray,linewidth=1.7pt](4.4,5)(4.4,1.6)(5.2,1.6)(5.2,1.4)(5.0,1)(5.0,.6)
\psline[linecolor=darkgray,linewidth=1.7pt](5.2,5)(5.2,1.8)(5.6,1.8)(5.6,1.4)(5.4,1)(5.4,.6)
\psline[linecolor=darkgray,linewidth=1.7pt](5.6,5)(5.6,2.0)(5.8,2.0)
\psline[linecolor=darkgray,linewidth=1.7pt](6,1.4)(5.8,1)(5.8,.6)
\psline[linecolor=darkgray,linewidth=1.7pt](6.8,2)(7,2)(7,1.4)(6.8,1)(6.8,.6)
\psline[linecolor=darkgray,linewidth=1.7pt](7,5)(7,2.2)(7.4,2.2)(7.4,1.4)(7.2,1)(7.2,.6)
\psline[linecolor=darkgray,linewidth=1.7pt](7.4,5)(7.4,2.4)(7.8,2.4)(7.8,1.4)(7.6,1)(7.6,.6)
\psline[linecolor=darkgray,linewidth=1.7pt](7.8,5)(7.8,2.6)(8.2,2.6)(8.2,1.4)(8,1)(8,.6)
\psline[linecolor=darkgray,linewidth=1.7pt](8.2,5)(8.2,2.8)(8.6,2.8)(8.6,1.4)(8.4,1)(8.4,.6)
\psline[linecolor=darkgray,linewidth=1.7pt](8.6,5)(8.6,3)(9,3)(9,5)
\psline[linecolor=darkgray,linewidth=1.7pt](9,3)(9,1.4)(8.8,1)(8.6,1.4)
\psline[arrowscale=1.5]{->}(9.4,2.5)(9.4,2.7)
\psframe[linecolor=black](.2,.4)(9.6,.6)
\end{pspicture}
$$  \caption{The nonzero admissible stem for $\tilde D_{2,a-1}$.}
\label{a2abexp7-1}\end{figure}

A similar argument works for the second web illustrated in
Figure~\ref{a2ab2exppo}. The third web depicted in
Figure~\ref{a2ab2exppo} is a little subtle. First one can see that
none of $\tilde D_{i,j}$ can be transformed if either $i>2$ or
$j<a-1$. Thus, we only need to check $\tilde D_{1,a-1}, \tilde
D_{1,a}, \tilde D_{2,a-1}$ and $\tilde D_{2,a}$ but we already know
about $ \tilde D_{1,a}, \tilde D_{2,a}$. Figure~\ref{a2abexp7} shows
the nonzero admissible stem for $ \tilde D_{1,a-1}$. As usual, we
draw a stem as a union of thick and darkly shaded lines. Note that
we have used relation ~\ref{a2defr23} in Figure \ref{relations}
exactly once which contributes $-[2]$. The Figure~\ref{a2abexp7-1}
shows the nonzero admissible stem for $\tilde D_{2,a-1}$. This
completes the proof of the lemma.
\end{proof}

\bibliographystyle{amsplain}

\end{document}